\documentclass[10pt]{amsart}
\usepackage{amssymb}
\usepackage{amsmath}
\usepackage{fancyhdr}
\usepackage[british]{babel}
\usepackage{geometry}
\usepackage{enumitem}
\usepackage{algpseudocode}
\usepackage{dsfont}
\usepackage{centernot}
\usepackage{xstring}
\usepackage{colortbl}
\usepackage[section]{algorithm}
\usepackage{graphicx}
\usepackage[font={footnotesize}]{caption}
\usepackage[usenames,dvipsnames,table]{xcolor}
\usepackage{pstricks,tikz}
\usepackage[h]{esvect}
\usepackage[
		bookmarksopen=true,
		bookmarksopenlevel=1,
		colorlinks=true,
		linkcolor=darkblue,
        linktoc=page,
		citecolor=darkblue,
]{hyperref}

\title[]{A rainbow blow-up lemma}
\date{\today}
\author[S.~Glock and F.~Joos]{Stefan Glock and Felix Joos}

\thanks{The research leading to these results was partially supported by the European Research Council
under the European Union's Seventh Framework Programme (FP/2007--2013) / ERC Grant
Agreement no. 306349 (S.~Glock) and the DFG, grant no.~JO 1457/1-1 (F.~Joos).
}

\newtheorem{theorem}[algorithm]{Theorem}
\newtheorem{prop}[algorithm]{Proposition}
\newtheorem{lemma}[algorithm]{Lemma}
\newtheorem{cor}[algorithm]{Corollary}
\newtheorem{fact}[algorithm]{Fact}

\theoremstyle{definition}

\newtheorem{defin}[algorithm]{Definition}

\newtheoremstyle{claimstyle}{5pt}{5pt}{\em}{5pt}{\em}{:}{5pt}{}
\theoremstyle{claimstyle}
\newtheorem{claim}{Claim}

\newtheoremstyle{stepstyle}{10pt}{5pt}{\em}{0pt}{\em}{:}{5pt}{}
\theoremstyle{stepstyle}
\newtheorem{step}{Step}

\numberwithin{equation}{section}

\definecolor{darkblue}{rgb}{0,0,0.5}

\def\noproof{{\unskip\nobreak\hfill\penalty50\hskip2em\hbox{}\nobreak\hfill%
       $\square$\parfillskip=0pt\finalhyphendemerits=0\par}\goodbreak}
\def\endproof{\noproof\bigskip}

\def\noclaimproof{{\unskip\nobreak\hfill\penalty50\hskip2em\hbox{}\nobreak\hfill%
       $-$\parfillskip=0pt\finalhyphendemerits=0\par}\goodbreak}

\def\endclaimproof{\noclaimproof\medskip}

\newdimen\margin
\def\textno#1&#2\par{
   \margin=\hsize
   \advance\margin by -4\parindent
          \setbox1=\hbox{\sl#1}
   \ifdim\wd1 < \margin
      $$\box1\eqno#2$$
   \else
      \bigbreak
      \hbox to \hsize{\indent$\vcenter{\advance\hsize by -3\parindent
      \it\noindent#1}\hfil#2$}
      \bigbreak
   \fi}

\def\proof{\removelastskip\penalty55\medskip\noindent\setcounter{claim}{0}\setcounter{step}{0}{\bf Proof. }} % in each main proof, claim counter set back

\def\lateproof#1{\removelastskip\penalty55\medskip\noindent\setcounter{claim}{0}\setcounter{step}{0}{\bf Proof of #1. }} % in each main proof, claim counter set back

\def\claimproof{\removelastskip\penalty55\medskip\noindent{\em Proof of claim: }}

\begin{document}

\newcommand{\new}[1]{\textcolor{red}{#1}}
\def\COMMENT#1{}
\def\TASK#1{}

\newcommand{\todo}[1]{\begin{center}\textbf{to do:} #1 \end{center}}

\def\eps{{\varepsilon}}
\newcommand{\ex}{\mathbb{E}}
\newcommand{\pr}{\mathbb{P}}
\newcommand{\cB}{\mathcal{B}}
\newcommand{\cA}{\mathcal{A}}
\newcommand{\cE}{\mathcal{E}}
\newcommand{\cS}{\mathcal{S}}
\newcommand{\cF}{\mathcal{F}}
\newcommand{\cG}{\mathcal{G}}
\newcommand{\bL}{\mathbb{L}}
\newcommand{\bF}{\mathbb{F}}
\newcommand{\bZ}{\mathbb{Z}}
\newcommand{\cH}{\mathcal{H}}
\newcommand{\cC}{\mathcal{C}}
\newcommand{\cM}{\mathcal{M}}
\newcommand{\bN}{\mathbb{N}}
\newcommand{\bR}{\mathbb{R}}
\def\O{\mathcal{O}}
\newcommand{\cP}{\mathcal{P}}
\newcommand{\cQ}{\mathcal{Q}}
\newcommand{\cR}{\mathcal{R}}
\newcommand{\cJ}{\mathcal{J}}
\newcommand{\cL}{\mathcal{L}}
\newcommand{\cK}{\mathcal{K}}
\newcommand{\cD}{\mathcal{D}}
\newcommand{\cI}{\mathcal{I}}
\newcommand{\cV}{\mathcal{V}}
\newcommand{\cT}{\mathcal{T}}
\newcommand{\cU}{\mathcal{U}}
\newcommand{\cX}{\mathcal{X}}
\newcommand{\cZ}{\mathcal{Z}}
\newcommand{\1}{{\bf 1}_{n\not\equiv \delta}}
\newcommand{\eul}{{\rm e}}
\newcommand{\Erd}{Erd\H{o}s}
\newcommand{\cupdot}{\mathbin{\mathaccent\cdot\cup}}
\newcommand{\whp}{whp }
\newcommand{\bX}{\mathcal{X}}
\newcommand{\bV}{\mathcal{V}}

\newcommand{\doublesquig}{%
  \mathrel{%
    \vcenter{\offinterlineskip
      \ialign{##\cr$\rightsquigarrow$\cr\noalign{\kern-1.5pt}$\rightsquigarrow$\cr}%
    }%
  }%
}

\newcommand{\defn}{\emph}

\newcommand\restrict[1]{\raisebox{-.5ex}{$|$}_{#1}}

\newcommand{\prob}[1]{\mathrm{\mathbb{P}}\left[#1\right]}
\newcommand{\expn}[1]{\mathrm{\mathbb{E}}\left[#1\right]}
\def\gnp{G_{n,p}}
\def\G{\mathcal{G}}
\def\lflr{\left\lfloor}
\def\rflr{\right\rfloor}
\def\lcl{\left\lceil}
\def\rcl{\right\rceil}

\newcommand{\qbinom}[2]{\binom{#1}{#2}_{\!q}}
\newcommand{\binomdim}[2]{\binom{#1}{#2}_{\!\dim}}

\newcommand{\grass}{\mathrm{Gr}}

\newcommand{\brackets}[1]{\left(#1\right)}
\def\sm{\setminus}
\newcommand{\Set}[1]{\{#1\}}
\newcommand{\set}[2]{\{#1\,:\;#2\}}
\newcommand{\krq}[2]{K^{(#1)}_{#2}}
\newcommand{\ind}[1]{$\mathbf{S}(#1)$}
\newcommand{\indcov}[1]{$(\#)_{#1}$}
\def\In{\subseteq}

\begin{abstract}  \noindent
We prove a rainbow version of the blow-up lemma of Koml\'os, S\'ark\"ozy and Szemer\'edi for $\mu n$-bounded edge colourings.
This enables the systematic study of rainbow embeddings of bounded degree spanning subgraphs.
As one application, we show how our blow-up lemma can be used to transfer the bandwidth theorem of B\"ottcher, Schacht and Taraz to the rainbow setting.
It can also be employed as a tool beyond the setting of $\mu n$-bounded edge colourings.
Kim, K\"uhn, Kupavskii and Osthus exploit this to prove several rainbow decomposition results.
%approximate versions of conjectures of Brualdi--Hollingsworth, Kaneko--Kano-Suzuki and Constantine.
%that in particular every properly coloured complete graph admits an approximate decomposition into rainbow almost Hamilton cycles yielding an approximate solution to a conjecture of Brualdi and Hollingsworth.
Our proof methods include the strategy of an alternative proof of the blow-up lemma given by R\"odl and Ruci\'nski, the switching method, and the partial resampling algorithm developed by Harris and Srinivasan.
\end{abstract}

\maketitle

\section{Introduction} \label{sec:intro}

A subgraph $H$ of an edge-coloured graph $G$ is called \defn{rainbow} if all its edges have different colours. 
Rainbow colourings appear for example in canonical Ramsey theory, and many open problems in combinatorics such as the Ryser--Brualdi--Stein conjecture on partial transversals in Latin squares and the graceful labelling conjecture can be phrased as rainbow subgraph problems. 
The central question is under which conditions on $G$ and its edge colouring a rainbow copy of $H$ is guaranteed.
Here, $H$ is usually a spanning subgraph such as a perfect matching~\cite{CPS:17,CP:17,HS:08,KY:17,pokrovskiy:16}, Hamilton cycle~\cite{AFR:95,APS:17,BM:17,BPS:17,CPS:17,ENR:83,FR:93,HT:86}, spanning tree~\cite{BLM:17,BPS:17,FK:08,PS:17}, or a general bounded degree graph~\cite{BKP:12,KSV:17,SV:17}.
Closely related questions concern properly coloured subgraphs and rainbow decompositions, which we shall discuss briefly in Section~\ref{subsec:related stuff}.

Clearly, a necessary condition for the existence of a rainbow copy of $H$ in $G$ is that $H$ is at least a subgraph of~$G$. 
Thus, the best one can hope for is to find a rainbow copy of $H$ in~$G$ `whenever' $H$ is a subgraph of~$G$.
The blow-up lemma of Koml\'os, S\'ark\"ozy and Szemer\'edi is a powerful tool to find spanning subgraphs, which, since its invention roughly 20~years ago, has significantly shaped the landscape of extremal combinatorics~\cite{BST:09,KSS:98,KSS:98a,KSS:01,KO:09,KO:13}.
It is tailored to be used after an application of Szemer\'edi's regularity lemma and roughly says that super-regular pairs behave like complete bipartite graphs in terms of embedding bounded degree subgraphs. 
%Koml\'os, S\'ark\"ozy and Szemer\'edi used it in a series of papers to prove several major open problems at the time. Since then, the blow-up lemma has become one of the standard tools in extremal combinatorics
%and has been used for many further breakthroughs in the area (e.g.~\cite{BST:09,KO:09,KO:13}).
In the present paper, we prove a rainbow blow-up lemma.
As one application, we transfer the bandwidth theorem of B\"ottcher, Schacht and Taraz~\cite{BST:09} to the rainbow setting. It would be interesting to find out whether other results can be transferred in a similar way.

In many of the classical rainbow problems, the host graph $G$ is complete.
For instance, Erd\H{o}s and Stein asked for the maximal $k$
such that any $k$-bounded edge colouring of $K_n$ contains a rainbow Hamilton cycle (cf.~\cite{ENR:83}). Here, an edge colouring is \defn{$k$-bounded} if each colour appears on at most $k$ edges.
After several subsequent improvements,\COMMENT{(see~e.g.~\cite{ENR:83,FR:93,HT:86})} Albert, Frieze, and Reed~\cite{AFR:95} showed that $k=\Omega(n)$, i.e.~there exists a constant $\mu>0$ such that for any $\mu n$-bounded edge colouring of~$K_n$, there exists a rainbow Hamilton cycle. Note that this is best possible up to the value of the constant~$\mu$. 
A natural generalization is to ask for general rainbow (spanning) subgraphs. For example, Frieze and Krivelevich~\cite{FK:08} showed that there exists some $\mu>0$ such that any almost spanning tree with bounded degree is contained as a rainbow copy in $K_n$ for any $\mu n$-bounded edge colouring of~$K_n$.
This was greatly improved by B\"ottcher, Kohayakawa, and Procacci~\cite{BKP:12}, who showed the following very general result. Given any $n/(51\Delta^2)$-bounded edge colouring of $K_n$ and any graph $H$ on $n$ vertices with $\Delta(H)\le \Delta$, one can find a rainbow copy of~$H$.
Their proof is based on the Lopsided Lov\'asz local lemma as well as the framework of Lu and Sz\'ekely~\cite{LS:07} for random injections. Using these tools, they show that a random injection $V(H)\to V(K_n)$ yields with positive probability a rainbow copy of~$H$. Kam\v{c}ev, Sudakov, and Volec~\cite{KSV:17} recently extended this result 
to the setting where $G$ is complete multipartite, and Sudakov and Volec~\cite{SV:17} considered the case when the number of cherries in $H$ is restricted (instead of the maximum degree).

There is a major stumbling stone if one wants to consider the above problem for incomplete host graphs,
say, for example quasi-random graphs with density $d$ for some arbitrarily small (fixed) $d>0$.
If $G$ is complete, then any injection $V(H)\to V(G)$ yields a valid embedding of $H$ (similarly for the multipartite setting). However, this is not the case for general host graphs $G$, where a random injection yields a valid embedding with exponentially small probability. Restricting the probability space to the `valid' injections does not seem to work with the Lu--Sz\'ekely framework, as the latter relies on the perfect symmetry of the setup. 
\COMMENT{Thus, the Lu--Sz\'ekely framework does not seem to be the right tool to study general host graphs~$G$.}

Some recent results on rainbow subgraphs in incomplete host graphs with $\mu n$-bounded edge colourings were obtained using the so-called `switching method'. 
For example, Coulson and Perarnau~\cite{CP:17} found rainbow perfect matchings in Dirac bipartite graphs, improving an approximate result of \cite{CPS:17}. The crucial property is that given a perfect matching $M$ and an edge $e\in M$ (which is in conflict with another edge in~$M$), there are many ways of `switching $e$ out of $M$' to obtain a new perfect matching which does not contain~$e$.
As another example, Coulson, Keevash, Perarnau, and Yepremyan~\cite{CKPY:18} show the existence of a rainbow $F$-factor in a graph $G$ 
whenever $\delta(G)\geq (\delta_F+o(1))n$, where $\delta_F$ is the minimum degree threshold for the existence of an $F$-factor (cf.~\cite{KO:09}). 
(Here $F$ is an arbitrary fixed hypergraph as their results also apply to hypergraphs.)\COMMENT{Divisibility}
However, the switching method seems to be limited to `simply structured' spanning graphs~$H$ with rich symmetry properties.

We are motivated by the following question:
\begin{quote}
\emph{Given a (dense) graph $G$ on $n$ vertices with a $\mu n$-bounded edge colouring and a (bounded degree) graph $H$ on $n$ vertices, is there a rainbow copy of $H$ in~$G$?}
\end{quote}
By proving a rainbow blow-up lemma, we provide a tool which allows for the systematic study of this question,
profiting from various techniques and methods which have been developed in the non-coloured setting.
In particular, we give affirmative answers to the above question if $G$ is quasi-random (see Corollary~\ref{cor:quasirandom}) or has sufficiently high minimum degree (see Section~\ref{sec:apps}).
%In the latter case, the `threshold' we obtain is the same as in the non-coloured setting, thus showing that a rainbow copy of $H$ exists `as soon as' $H$ is a subgraph of~$G$.
We remark that the constant $\mu>0$ we obtain is very small.

Nevertheless, our rainbow blow-up lemma has applications beyond this setting, in at least the following two aspects.
Firstly, it can still be applied even if the number of available colours is only slightly larger than the number of edges in the desired subgraph.
Secondly, one can even obtain approximate decompositions, for example into Hamilton cycles and $H$-factors.
This has recently been demonstrated by Kim, K\"uhn, Kupavskii and Osthus (see Section~\ref{subsec:related stuff} for further details).

We will discuss the blow-up lemma in more detail in the next subsection. As mentioned above, commonly used techniques like the Lu--Sz\'ekely framework and the switching method do not seem capable of dealing with quasi-random host graphs $G$ and/or general graphs $H$. (Here, the idea would have been to use the original blow-up lemma as a `blackbox' result.)
Another natural question is whether the proof of the blow-up lemma can be adapted to work in the rainbow setting. 
In the original proof due to Koml\'os, S\'ark\"ozy and Szemer\'edi, the vertices of $H$ are embedded one by one using a randomized algorithm, until all but a small fraction of the vertices are embedded, and the embedding is then completed using the K\"onig-Hall theorem. Note that this approach is extremely vulnerable in the rainbow setting, as already a constant number of choices can render the algorithm unsuccessful (as a vertex may be incident to only a constant number of different colours), which seems rather hopeless.

To overcome these obstacles, several new ideas are needed.
As an underlying strategy, we employ the alternative proof of the blow-up lemma given by R\"odl and Ruci\'nski,
and combine it with various techniques such as the partial resampling algorithm, the switching method, and a parallelization of the embedding procedure.
We will provide a detailed proof overview in Section~\ref{sec:sketch}. 
Sections~\ref{sec:colour splitting} and~\ref{sec:main proof} contain the main proof.
In Section~\ref{sec:apps}, we demonstrate the applicability of our rainbow blow-up lemma.
%Most notably, we will show how it can be used to prove a rainbow version of the celebrated bandwidth theorem of B\"ottcher, Schacht and Taraz~\cite{BST:09}, which we discuss in more detail in Section~\ref{subsec:intro bandwidth}.

\subsection{Rainbow blow-up lemma}\label{subsec:blow-up intro}

%In this paper we lift the blow-up lemma to a higher level
%and show that we can still find any bounded degree spanning subgraph even if we impose strong restrictions on this subgraph.
In order to state a simplified version of our rainbow blow-up lemma,
we need some more terminology.
For $k\in \bN$, we write $[k]_0:=[k]\cup \Set{0}=\Set{0,1,\dots,k}$. We say that $(H,G,R,(X_i)_{i\in[r]_0},(V_i)_{i\in[r]_0})$ is a \defn{blow-up instance} if the following hold:
\begin{itemize}
\item $H$ and $G$ are graphs, $(X_i)_{i\in[r]_0}$ is a partition of $V(H)$ into independent sets, $(V_i)_{i\in[r]_0}$ is a partition of $V(G)$, and $|X_i|=|V_i|$ for all $i\in[r]_0$;
\item $R$ is a graph on $[r]$ such that for all distinct $i,j\in [r]$, the graph $H[X_i,X_j]$ is empty if $ij\notin E(R)$.
\end{itemize}
Here, $X_0$ and $V_0$ are so-called `exceptional sets'. For simplicity, we assume in this subsection that they are empty. 

For a graph $G$ and two disjoint subsets $S,T\In V(G)$, denote by $e_G(S,T)$ the number of edges of $G$ with one endpoint in $S$ and the other one in $T$, and define $$d_G(S,T):=\frac{e_G(S,T)}{|S||T|}$$ as the \defn{density} of the pair $S,T$ in $G$.
We say that the bipartite graph $G[V_1,V_2]$ is \defn{lower $(\eps,d)$-super-regular}
if
\begin{itemize}
\item for all $S\In V_1$ and $T\In V_2$ with $|S|\ge \eps|V_1|$, $|T|\ge \eps |V_2|$, we have $d_G(S,T)\ge d-\eps$;
\item for all $i\in [2]$ and $v\in V_i$, we have $|N_G(v)\cap V_{3-i}|\ge (d-\eps)|V_{3-i}|$.
\end{itemize}

We say that the blow-up instance $(H,G,R,(X_i)_{i\in[r]_0},(V_i)_{i\in[r]_0})$ is \defn{lower $(\eps,d)$-super-regular} if for all $ij\in E(R)$, the bipartite graph $G[V_i,V_j]$ is lower $(\eps,d)$-super-regular.

We now state a simplified version of the rainbow blow-up lemma. The full statement can be found in Lemma~\ref{lem:blow up} and also allows exceptional vertices and candidate sets. Moreover, it does not only apply to rainbow embeddings, but to slightly more general conflict-free embeddings.

\begin{lemma}[Rainbow blow-up lemma---simplified]\label{lem:blow-up simple}
For all $d,\Delta,r$,
there exist $\eps=\eps(d,\Delta)$, $\mu=\mu(d,\Delta,r)>0$ and an $n_0\in \mathbb{N}$
such that the following holds for all $n\geq n_0$.
%Let $1/n \ll \mu , \eps \ll d,1/\Delta$ and $\mu \ll 1/r$.

Let $(H,G,R,(X_i)_{i\in[r]},(V_i)_{i\in[r]})$ be a lower $(\eps,d)$-super-regular blow-up instance.
Assume further that
%$\Delta(H),\Delta(R)\le \Delta$ and $|V_i|=(1\pm \eps)n/r$ for all $i\in[r]$.
\begin{enumerate}[label={\rm (\roman*)}]
	\item $\Delta(H), \Delta(R)\leq \Delta$,
	\item $|V_i|=(1\pm \eps)n/r$ for all $i\in[r]$.
\end{enumerate}
Then, given any $\mu n$-bounded edge colouring of $G$,
there exists a rainbow embedding of $H$ into $G$ (where $X_i$ is mapped to $V_i$ for all $i\in [r]$).
\end{lemma}

Observe that $\eps$ does not depend on $r$, which is crucial in applications. Originally the blow-up lemma was formulated with $R$ being the clique on at most $\Delta$ vertices.
In essentially all applications it is applied to many clique blow-ups iteratively.
However, in order to have a useful rainbow blow-up lemma,
we cannot apply it independently to two or more clique blow-ups as
in both parts the blow-up lemma may use the same colour.
We will discuss this issue further in Section~\ref{sec:sketch}.

\subsection{A rainbow bandwidth theorem and quasirandom host graphs} \label{subsec:intro bandwidth}

One of the fundamental results of extremal combinatorics is the Erd\H{o}s--Stone theorem, stating that for a fixed graph $H$, every large graph $G$ on $n$ vertices with average degree at least $(1-1/(\chi(H)-1)+o(1))n$ contains $H$ as a subgraph.
The bandwidth theorem can be viewed as an analogue of the Erd\H{o}s--Stone theorem when $H$ is a spanning subgraph of~$G$. Clearly, in this setting, one has to replace the average degree condition by a minimum degree condition in order to obtain sensible results. 
Let $\ell:=\chi(H)$ and assume that $H$ has bounded degree.
A long line of research confirmed for various cases that $\delta(G)\geq (1-1/\ell+o(1))n$ suffices to find $H$ as a spanning subgraph in $G$, for instance when $H$ is a spanning tree, a (power of a) Hamilton cycle, or a clique factor.
A conjecture of Bollob\'as and Koml\'os, which became known as the bandwidth conjecture, attempted to
generalize all the mentioned results by claiming that $\delta(G)\geq (1-1/\ell+o(1))n$ suffices whenever $H$ has not too strong expansion properties (in this case measured by the parameter bandwidth). The \emph{bandwidth} of a graph $H$ with vertex set $[n]$ is defined as
$\min_{\pi}\max \{|\pi(i)-\pi(j)|\colon ij\in E(H)\}$ where the minimum ranges over all permutations on $[n]$.

B\"ottcher, Schacht and Taraz~\cite{BST:09} proved the bandwidth conjecture roughly ten years ago. We refer to their paper for more information on the history of the conjecture. Graph classes with appropriately small bandwidth include for example powers of Hamilton cycles, $F$-factors and bounded degree trees and planar graphs. We also remark that the minimum degree bound given by the bandwidth theorem is not necessarily the optimal threshold (cf.~Section~\ref{subsec:bandwidth}). 
Here, we extend the bandwidth theorem to the rainbow setting. 
We will prove the following theorem in Section~\ref{subsec:bandwidth} 
by replacing the usual blow-up lemma with our rainbow blow-up lemma in the approach of~\cite{BST:09}.

%Our blow-up lemma can be used to prove several of the aforementioned results for bounded degree graphs
%(up to worse dependencies of the parameters)
%and we can also quite easily extend many seminal results in extremal graph theory to the rainbow setting.
%One particular example is the bandwidth theorem which was proved by B\"ottcher, Schacht and Taraz~\cite{BST:09}
%and was priorly known as the bandwidth conjecture of Bollob\'as and Komlos.

\begin{theorem}\label{thm:bandwsimple}
For all $\delta,\Delta,\ell$,
there are $\beta,\mu>0$ and an $n_0\in \mathbb{N}$
such that the following holds for all $n\geq n_0$.
%Suppose $1/n\ll \mu,\beta \ll \delta,1/\Delta,1/\ell$.
Suppose $G$ is a graph on $n$ vertices with $\delta(G)\geq (1-1/\ell+\delta)n$.
Suppose $H$ is a graph on $n$ vertices with $\Delta(H)\leq \Delta$, bandwidth at most $\beta n$, and $\chi(H)\le \ell$.
Then, given any $\mu n$-bounded edge colouring of $G$, the graph $G$ contains a rainbow copy of $H$.
\end{theorem}

%The bandwidth theorem generalized many previously known results
%about the minimum degree threshold for the containment of specific families of bounded degree graphs.
%Among other families this includes trees, planar graphs, $H$-factors and powers of Hamilton cycles.

%\medskip

A particular example which is covered by the bandwidth theorem is the case when $H$ is a tree.
This was known as Bollob\'as' conjecture and was proved by Koml\'os, S\'ark\"ozy and Szemer\'edi~\cite{KSS:95}.
In Section~\ref{sec:appstree}, we present an easy proof of this result in the rainbow setting.
Our main motivation is to give the reader a straightforward example of how to apply our rainbow blow-up lemma. (The proof in~\cite{KSS:95} is without the blow-up lemma and thus much longer.)

\medskip

Another straightforward corollary of our rainbow blow-up lemma is the existence of rainbow bounded degree graphs in quasi-random host graphs.
Let us say a graph $G$ on $n$ vertices is \emph{$(\eps,d)$-dense}
if for all disjoint $S,T\subseteq V(G)$ with $|S|,|T|\geq \eps n$,
we have $e_G(S,T)\geq d|S||T|$.

\begin{cor} \label{cor:quasirandom}
%Suppose $1/n\ll \mu, \eps \ll d,1/\Delta$.
For all $d,\Delta$,
there are $\eps,\mu>0$ and an $n_0\in \mathbb{N}$
such that the following holds for all $n\geq n_0$.
Suppose $G,H$ are graphs on $n$ vertices,
$\delta(G)\geq dn$, the graph $G$ is $(\eps,d)$-dense, and $\Delta(H)\leq \Delta$.
Then, given any $\mu n$-bounded edge colouring of $G$,
the graph $G$ contains a rainbow copy of $H$.
\end{cor}
\COMMENT{
\proof
By the Hajnal--Szemer\'edi-theorem,
there is a partition $(X_1,\ldots,X_{\Delta+1})$ of $V(H)$ into independent sets such that their sizes differ by at most $1$.
Observe that with probability at least $1/2$
a random partition $(V_1,\ldots,V_{\Delta+1})$ of $V(G)$ such that $|V_i|=|X_i|$ for all $i\in [\Delta+1]$
induces lower $(2\eps,d/2)$-super-regular graphs $G[V_j,V_{j'}]$ for all distinct $j,j'\in [\Delta+1]$ by some standard Chernoff-type inequality.
Let $(V_1,\ldots,V_{\Delta+1})$ be such a partition.
Now the application of Lemma~\ref{lem:blow-up simple} with $R=K_{\Delta+1}$ completes the proof.
\endproof}

This is the first result that considers spanning rainbow embeddings in graphs of arbitrarily small (but fixed) density.

\subsection{Further applications and related problems} \label{subsec:related stuff}

An edge colouring is \defn{locally $k$-bounded} if each colour appears at most $k$ times on the edges which are incident to a single vertex. (For instance, a locally $1$-bounded colouring is proper.) For many of the aforementioned results, there is a corresponding result where the colouring is not $k$-bounded, but only locally $k$-bounded, and instead of aiming for a rainbow copy of $H$, one aims for a properly coloured copy of $H$. Our proof of the rainbow blow-up lemma transfers to this setting as well (with $k=\mu n$). In fact, it is then much easier, because there are only local conflicts and no global conflicts,
and a simple modification of the original proofs already yields such a result.

Note that every locally $O(1)$-bounded edge colouring is (globally) $O(n)$-bounded. For example, every proper edge colouring is $n/2$-bounded. 
Quite a few open problems are formulated in this setting. 
For example, the Ryser--Brualdi--Stein conjecture asks for a rainbow matching of size $n-1$ in every properly $n$-edge-coloured~$K_{n,n}$, and Andersen~\cite{andersen:89} conjectured that there is a rainbow path of length $n-2$ in every properly edge-coloured~$K_n$.
The currently best approximate result for the first problem is due to Hatami and Shor~\cite{HS:08}. 
An approximate result for the second problem was achieved in a breakthrough result by Alon, Pokrovskiy and Sudakov~\cite{APS:17} and slightly improved by Balogh and Molla~\cite{BM:17}. One may also ask for any other rainbow subgraph with at most $n$ edges, for instance a particular (almost) spanning tree or an (almost) spanning collection of cycles. (Observe that a proper edge-colouring may only use $n$ colours.
Hence it only makes sense to ask for rainbow subgraphs with at most $n$ edges.) 

Recently, Montgomery, Pokrovskiy and Sudakov~\cite{MPS:18a} showed that $K_n$ equipped with a locally $k$-bounded edge colouring
contains any tree with at most $(1-o(1))n/k$ vertices as a rainbow subgraph.
This also implies an approximate version of Ringel's conjecture on decompositions of $K_n$ into a given tree.

Instead of finding a single rainbow subgraph, it is also natural to ask for decompositions into rainbow subgraphs. For instance, there are three conjectures due to Brualdi--Hollingsworth, Kaneko--Kano--Suzuki, and Constantine regarding the decomposition of properly edge-coloured complete graphs into rainbow (isomorphic) spanning trees. Montgomery, Pokrovskiy and Sudakov~\cite{MPS:18b} obtain strong results which prove these conjectures approximately (for not necessarily isomorphic trees).
Independently, Kim, K\"uhn, Kupavskii and Osthus~\cite{KKKO:18} use our blow-up lemma to prove several related results. In their setting, the given colouring does not need to be proper, but can be locally $n/\log^5 n$-bounded, say. On the other hand, the global condition is slightly more restrictive than in a proper colouring.
As one example, they show the following.

\begin{theorem}[Kim, K\"uhn, Kupavskii and Osthus~\cite{KKKO:18}]
For any $\eps > 0$, there exists $n_0\in \bN$ such that for all $n\geq  n_0$, any $(1 - \eps)n/2$-bounded, locally $n/\log^5 n$-bounded edge-colouring of $K_n$ contains $(1/2- \eps)n$ edge-disjoint rainbow Hamilton cycles.
\end{theorem}
This also establishes approximate solutions to the above tree decomposition conjectures. 
Analogous results hold for approximate decompositions into rainbow $H$-factors~\cite{KKKO:18}.

%Instead of finding a single rainbow subgraph, it is also natural to ask for decompositions into rainbow subgraphs. 
%Beside several conjectures on the existence of a single rainbow subgraph,
%there are also three conjectures due to Brualdi--Hollingsworth, Kaneko--Kano--Suzuki,
%and Constantine regarding the decomposition of the complete graph into rainbow (isomorphic) spanning trees. 
%Kim, K\"uhn, Kupavskii and Osthus use our rainbow blow-up lemma and a hypergraph matching result to prove approximate solutions to these conjectures. 
%Montgomery, Pokrovskiy and Sudakov~\cite{MPS:18b} obtained independently similar results.\footnote{
%To be more precise, they prove that every properly (i.e.~$1$-bounded) edge-coloured $K_n$ contains $(1/2- o(1))n$ (not necessarily isomorphic) trees.}
%
%\begin{theorem}[Kim, K\"uhn, Kupavskii and Osthus~\cite{KKKO:18}]
%For any $\eps > 0$, there exists $n_0\in \bN$ such that for all $n\geq  n_0$, any $(1 - \eps)n/2$-bounded, locally $n/\log^5 n$-bounded edge-colouring of $K_n$ contains $(1/2- \eps)n$ edge-disjoint rainbow Hamilton cycles.
%\end{theorem}
%Analogous results also hold for approximate decompositions into rainbow $H$-factors. 

The following observation indicates why our rainbow blow-up lemma can be applied in this setting.
Suppose for simplicity that we consider a properly edge-coloured complete graph.
If one chooses a random subset $U$ of vertices of size $\mu n$, then with high probability, the colouring restricted to this subset is $\mu |U|$-bounded,\COMMENT{expected number of edges in $U$ of fixed colour is at most $(|U|/n)^2 n/2 =\mu |U|/2 $, appear independently because colouring is proper.} and the blow-up lemma can be applied to $U$ to complete a partial embedding (constructed outside $U$ using different methods, e.g.~via a hypergraph matching argument).

Rainbow decomposition results may in fact hold in a more general setting.
Kim, K\"uhn, Osthus and Tyomkyn~\cite{KKOT:ta} proved a blow-up lemma for approximate decompositions into arbitrary (bounded degree) graphs. 
In combination with our result and methods there is hope for a rainbow blow-up lemma for approximate decompositions.

\section{Proof overview} \label{sec:sketch}

Our starting point is the alternative proof of the blow-up lemma given by R\"odl and Ruci\'nski~\cite{RR:99}. We briefly sketch their approach, before considering the rainbow setting.

Suppose that $(H,G,R,(X_i)_{i\in[r]},(V_i)_{i\in[r]})$ is a blow-up instance.
Assume that $\Delta(H)\le \Delta$ and that $G[V_i,V_j]$ is super-regular\COMMENT{not yet defined} for all $1\le i<j\le r$, i.e.~$R=K_r$.
We want to find an embedding $\phi$ of $H$ into $G$ (where $X_i$ is mapped onto $V_i$).
The strategy is to embed $H$ in $r$ rounds, where in round $i$, all vertices of $X_i$ are embedded simultaneously into~$V_i$. Clearly, this corresponds to a perfect matching between $X_i$ and~$V_i$.
Yet not every perfect matching yields a valid embedding. They keep track of this by defining a (dynamic) `candidacy graph' $A^j$ for each pair of clusters $(X_j,V_j)$. This candidacy graph depends on the partial embedding up to round $i$.
More precisely, for $j>i$, let $A_i^j$ denote the candidacy graph for $(X_j,V_j)$ after $X_1,\dots,X_i$ have been embedded, where $xv\in E(A_i^j)$ if and only if $v$ is still a candidate to be the image of~$x$.
In the beginning, all candidacy graphs are complete bipartite graphs, but during the embedding, they become gradually sparser. 
The main idea is to maintain super-regularity of the candidacy graphs.
More precisely, given an embedding of the clusters $X_1,\dots,X_i$, we assume that the candidacy graphs $A_{i}^{i+1},\dots,A_i^{r}$ are super-regular (which is true for $i=0$), and now have to find a perfect matching $\sigma$ of $A_i^{i+1}$ (which defines an embedding of $X_{i+1}$ into $V_{i+1}$) in such a way that the updated candidacy graphs $A_{i+1}^{i+2},\dots,A_{i+1}^{r}$ are again super-regular.
Fortunately, it turns out that if $\sigma$ is chosen randomly among all perfect matchings of $A_i^{i+1}$, then this desired property holds with high probability.
The key tool to show this is the so-called Four graphs lemma from \cite{RR:99} (see~Lemma~\ref{lem:four graphs}).
For this to work, we actually need a stronger assumption on~$H$, that is, we require that the clusters $X_i$ are not only independent in~$H$, but $2$-independent. However, using the Hajnal--Szemer\'edi theorem, it is not hard to reduce the general case to this setting. This will also be the first step in our proof (see Section~\ref{subsec:matching split}).

Now, assume that $c\colon E(G)\to C$ is a $\mu n$-bounded edge colouring of $G$. We might try to proceed as above. Suppose that we have already found a rainbow embedding $\phi_i$ of $H[X_1\cup \dots \cup X_i]$ into $G[V_1\cup \dots \cup V_i]$ and that the candidacy graphs $A_{i}^{i+1},\dots,A_i^{r}$ are super-regular. Clearly, the definition of candidacy graphs needs to be adapted to the rainbow setting.
Consider $x\in X_j$ and $v\in V_j$ for $j>i$. Let $F_{xv}$ denote the set of edges $\phi_i(y)v$ where $y\in N_H(x) \cap (X_1\cup \dots \cup X_i)$.
Previously, $v$ was a candidate for $x$ if and only if $F_{xv}\In E(G)$, where $F_{xv}$ is the set of edges in $G$ which are used if $x$ is embedded at $v$. Now, we need in addition that the edges of $F_{xv}$ have mutually distinct colours, and do not have a colour which is already used by $\phi_i$.
Suppose that this is the case and $A_{i}^{i+1}$ is the candidacy graph for $(X_{i+1},V_{i+1})$.
In the uncoloured setting, any perfect matching $\sigma$ of $A_{i}^{i+1}$ yields a valid embedding of $H[X_1\cup \dots \cup X_{i+1}]$ into $G[V_1\cup \dots \cup V_{i+1}]$ (and almost all of them leave the updated candidacy graphs super-regular).
In contrast to this, by far not every perfect matching $\sigma$ of $A_{i}^{i+1}$ yields a rainbow embedding of $H[X_1\cup \dots \cup X_{i+1}]$.
For this to be the case, we need that $c(F_{x\sigma(x)})\cap c(F_{x'\sigma(x')})=\emptyset$ for all distinct $x,x'\in X_{i+1}$.
This can be modelled by defining a suitable conflict system on the edges of $A_{i}^{i+1}$, where two edges $xv,x'v'$ \defn{conflict} if $c(F_{xv})\cap c(F_{x'v'})\neq \emptyset$.
Observe that a conflict-free perfect matching would yield a rainbow embedding of $H[X_1\cup \dots \cup X_{i+1}]$.
Crucially, since $c$ is $\mu n$-bounded, every edge conflicts with at most $\Delta(H)\mu n$ other edges.\COMMENT{Here already assume $2$-independence} Using the switching method, we can show that a uniformly random perfect matching is conflict-free with positive probability.
This step uses a recent result of Coulson and Perarnau~\cite{CP:17}.

For the embedding procedure to work, we have to make sure that the updated candidacy graphs are again super-regular. However, this might easily be wrong, as we might accidentally isolate a vertex. For instance, vertex $v\in V_{i+2}$ might only see $\mu^{-2}$ colours, say, but all of those colours might already be used in the embedding $\phi_{i+1}$. Then vertex $v$ cannot be a candidate for any $x\in X_{i+2}$ and thus the embedding is stuck.

To overcome this issue, we reserve an exclusive set of colours for each embedding round in the beginning. For instance, we might partition the colour set $C$ into sets $(C_{ij})_{1\le i<j\le r}$, and instead of embedding $H$ into $G$, we embed $H$ into the graph $G^\ast$ which between clusters $V_i,V_j$ only contains those edges of $G[V_i,V_j]$ whose colour is contained in $C_{ij}$. This decouples much of the colour dependency between the embedding rounds and allows us to analyse the updated candidacy graphs as before using the Four graphs lemma.
Of course, this strategy only works if the pairs $G^\ast[V_i,V_j]$ are super-regular. A natural approach to show the existence of a suitable partition of $C$ is to partition $C$ randomly. Simplified, the problem we face here is as follows: Consider the complete bipartite graph $K_{n,n}$ with any $\mu n$-bounded edge colouring. If we `activate' each colour independently with probability $1/2$, what is the probability that the activated edges induce a $(\sqrt{\mu},1/2)$-super-regular pair, say?
This seems to be an interesting problem in its own right, and we are not aware of any coverage in the literature. Heuristically one expects this to be a rare event of exponentially small probability. However, the Lov\'asz local lemma does not seem flexible enough to solve this problem.
We will use the partial resampling algorithm introduced by Harris and Srinivasan~\cite{HS:18} to deal with this problem (see Section~\ref{sec:colour splitting}).

One remaining issue arises from the applicability of the blow-up lemma. In the original version of the blow-up lemma, the regularity parameter~$\eps$ needs to be small as a function of the number of clusters~$r$. However, in a regularity partition obtained via an application of Szemer\'edi's regularity lemma, the number of clusters~$r$ is large as a function of~$\eps$. Usually, this is dealt with roughly as follows: one finds a vertex partition of the reduced graph $R$ into $r/k$ cliques of order $k$, where $k$ is `small', and then applies the blow-up lemma to each clique independently. 
This approach is not feasible in the rainbow setting, as the union of rainbow subgraphs does not necessarily yield a rainbow subgraph.
Also, splitting the colour set in the beginning in order to reserve an exclusive set of colours for each application of the blow-up lemma does not solve this problem, as the density of the regular pairs is divided by $r/k$ through the splitting and will thus be too small compared to~$\eps$. We thus prove a `global' version of the blow-up lemma which allows the number of clusters $r$ to be large, but requires that $\Delta(R)\le \Delta$.
Such a version of the blow-up lemma was first proved by Csaba~\cite{csaba:07} following the proof in~\cite{KSS:97}.
The proof of R\"odl and Ruci\'nski~\cite{RR:99} can also be adapted to work in this setting (see e.g.~\cite{KKOT:ta}).
However, we have to be careful with the colour splitting, for the very reason discussed above:
if we reserved an exclusive colour set for each pair $V_i,V_j$ with $ij\in E(R)$,
the density of the regular pairs would be far too low.
This can be dealt with by `parallelizing' the embedding (for different reasons, this has also been done in~\cite{KKOT:ta}).
We split $V(R)$ into $2$-independent sets $J_1,\dots,J_T$, where we might take $T=\Delta(R)^2+1$.
We then reserve an exclusive colour set $C_{t_1t_2}$ for each pair $1\le t_1<t_2\le T$.
The sparsified graph $G^\ast$ contains between clusters $V_i,V_j$ only those edges whose colour is contained in $C_{t_1t_2}$, where $i\in J_{t_1}$ and $j\in J_{t_2}$.
We now embed $H$ in just $T$ rounds, where in round $t\in [T]$, we embed all clusters $X_i$ with $i\in J_t$ simultaneously. Consequently, the perfect matching $\sigma$ is not drawn as a uniform perfect matching of one candidacy graph, but is the union of $|J_t|$ such perfect matchings. Fortunately, this still works nicely with the switching method and the Four graphs lemma.

The resulting proof can be found in Section~\ref{subsec:main proof}.

\section{Notation}

We only consider finite, undirected and simple graphs. For a graph $G$, we let $V(G)$ and $E(G)$ denote the vertex set and edge set, respectively. Moreover, for a vertex $v$, let $d_G(v)$ denote the degree of $v$, and let $d_G(v,A)$ denote the number of neighbours of $v$ in $A\In V(G)$. 
If we deal with several graphs simultaneously, we say \defn{$G$-neighbour} to clarify which graph we refer to.
As usual, $\delta(G)$ and $\Delta(G)$ denote the minimum and maximum degree of $G$, respectively.
For disjoint subsets $A,B\In V(G)$, 
let $G[A,B]$ denote the bipartite subgraph of $G$ between $A$ and $B$
and $G[A]$ the subgraph in $G$ induced by $A$.
We let $G^2$ denote the square of $G$, i.e.~the graph obtained from $G$ by adding edges between vertices which have a common neighbour in~$G$.
A subset $X\In V(G)$ is $2$-independent if it is independent in~$G^2$.

%Let $G$ be a graph and $c\colon E(G)\to C$ an edge colouring of $G$. For a colour $i$, we let $G_i$ denote the subgraph of $G$ with vertex set $V(G)$ and edge set $c^{-1}(i)$, and for a subset $C'\In C$ of colours, we write $G_{C'}:=\bigcup_{i\in C'}G_i$. If $G$ is clear from the context, we might write $d_i(v),d_{C'}(v)$ instead of $d_{G_i}(v),d_{G_{C'}}(v)$, respectively.
Let $G$ be a graph.
Given a set $C$, 
a function $c\colon E(G)\to 2^C$ is called an \emph{edge set colouring} of~$G$. A colour $\alpha \in C$ \emph{appears} on an edge $e$ if $\alpha\in c(e)$. 
We say that $c$ is \defn{$k$-bounded} if each colour appears on at most $k$ edges. 
Moreover, we say that $c$ is \defn{$(k,\Delta)$-bounded} if it is $k$-bounded and $|c(e)|\le \Delta$ for all $e\in E(G)$.
For $C'\In C$, let $G_{C'}$ denote the subgraph of $G$ with vertex set $V(G)$ and edge set $\set{e\in E(G)}{c(e)\In C'}$. 
 We let $d_G^{\alpha}(v)$ denote the number of edges incident to $v$ on which $\alpha$ appears.
If $G$ is clear from the context, we may simply write $d_{C'}(v)$ and $d^\alpha(v)$ instead of $d_{G_{C'}}(v)$ and $d_G^{\alpha}(v)$, respectively.

Given $\eps>0$ and $d\in[0,1]$, a bipartite graph $G$ with vertex classes $V_1,V_2$ is called \defn{$(\eps,d)$-regular} if for all pairs $S\In V_1$ and $T\In V_2$ with $|S|\ge \eps|V_1|$, $|T|\ge \eps |V_2|$,
we have $|d_G(S,T)-d|\le \eps$.
We simply say that $G$ is \defn{$\eps$-regular} if $G$ is $(\eps,d_G(V_1,V_2))$-regular.
Moreover, an $(\eps,d)$-regular graph $G$ is called \defn{$(\eps,d)$-super-regular} if for all $v\in V_1$, we have $d_G(v)=(d\pm \eps)|V_2|$ and for all $v\in V_2$ we have $d_G(v)=(d\pm \eps)|V_1|$.
Clearly, if $G$ is $(\eps,d)$-super-regular, it is also lower $(\eps,d)$-super-regular.

%For any $n\in \mathbb{N}$, we define $[n]=\{1,\ldots,n\}$ and $[n]_0:=[n]\cup \{0\}$.
For a finite set $S$ and $i\in \mathbb{N}$,
we write $\binom{S}{i}$ for the set of all subsets of $S$ of size $i$
and $2^S$ for the powerset of $S$. 
For a set $\Set{i,j}$, we sometimes simply write $ij$.
For $a,b,c\in \mathbb{R}$,
we write $a=b\pm c$ whenever $a\in [b-c,b+c]$.
For $a,b,c\in (0,1]$,
we sometimes write $a\ll b \ll c$ in our statements meaning that there are increasing functions $f,g:(0,1]\to (0,1]$
such that whenever $a\leq f(b)$ and $b \leq g(c)$,
then the subsequent result holds.

\section{Colour splitting} \label{sec:colour splitting}

The goal of this section is to prove the following lemma, which will allow us to reserve exclusive colour sets for our embedding rounds in the proof of the rainbow blow-up lemma.
Our main tools are the partial resampling algorithm introduced by Harris and Srinivasan~\cite{HS:18} and McDiarmid's inequality.

\begin{lemma} \label{lem:separate colours}
Let $1/n \ll \mu \ll 1/r,\eps, 1/t,d,1/\Delta$. Suppose that $G$ is a graph with vertex partition $V(G)=V_1\cupdot\dots\cupdot V_r$ and $R_S\In R$ are graphs on $[r]$ such that
\begin{enumerate}[label={\rm (\roman*)}]
	\item $n\le |V_i|\le 2n$ for all $i\in[r]$;
	\item $G[V_i,V_j]$ is (lower) $(\eps,d)$-regular for all $ij\in E(R)$;
	\item $G[V_i,V_j]$ is (lower) $(\eps,d)$-super-regular for all $ij\in E(R_S)$.
\end{enumerate}
Let $c\colon E(G)\to 2^C$ be a $(\mu n,\Delta)$-bounded edge set colouring of $G$ with $|C|\le \mu^{-1.3}n$.\COMMENT{deleted: Let $(E_\ell)_{\ell\in [t]}$ be any partition of $E(R)$. Doesn't work with the bandwidth application. Thus, we split every pair into many colour exclusive subpairs. Keep in mind}
Then there exists a partition $(C_\ell)_{\ell\in [t]}$ of $C$ such that for all $\ell \in [t]$ and all $ij \in E(R)$, there is a subgraph $G_\ell^{ij}\In G_{C_\ell}[V_i,V_j]$ such that $G_\ell^{ij}$ is (lower) $(2\eps,d/t^\Delta)$-regular, and (lower) $(2\eps,d/t^\Delta)$-super-regular if $ij \in E(R_S)$.
\end{lemma}

We remark that if $c$ is uniform, i.e.~$c(e)=\Delta$ for all $e\in E(G)$, for some $\Delta\in \bN$, then we can take $G_\ell^{ij}= G_{C_\ell}[V_i,V_j]$. In particular, this applies if $c$ is a normal edge colouring. In the general case, we introduce `dummy colours' to make $c$ uniform.

In the proof of Lemma~\ref{lem:separate colours},
we will partition $C$ randomly and then show that with positive probability, all the obtained subpairs are (super-)regular.
In order to motivate our approach, we briefly discuss the following simplified variant of the problem (which will be the running example throughout this section): Consider $K_{n,n}$ with any $\mu n$-bounded edge colouring. If we `activate' each colour independently with probability $1/2$, what is the probability that the activated edges induce a super-regular pair?
(An edge is activated if its colour is activated.)
Using McDiarmid's inequality, it is not hard to see that with high probability, the activated subgraph is indeed $(\eps,1/2)$-regular.
We first recall McDiarmid's inequality and then prove a more elaborated version of what we just indicated.

\begin{theorem}[McDiarmid's inequality, see~\cite{mcdiarmid:89}\COMMENT{Lemma~1.2}] \label{thm:McDiarmid}
Suppose $X_1,\dots,X_m$ are independent Bernoulli random variables and suppose $w_1,\dots,w_m\in [0,W]$.
Suppose $X$ is a real-valued random variable determined by $X_1,\dots,X_m$ such that changing the outcome of $X_i$ changes $X$ by at most $w_i$ for all $i\in [m]$.
Then, for all $t>0$, we have $$\prob{|X-\expn{X}|\ge t} \le 2 \eul^{-\frac{2t^2}{W\sum_{i=1}^m w_i}}.$$
\end{theorem}

\begin{prop}\label{prop:regularity slicing}
Let $1/n \ll \mu \ll  \eps , p,d,1/\Delta$.
Suppose that $G$ is a (lower) $(\eps,d)$-regular graph between $V,V'$ with $\eps n\le |V|,|V'|\le 2n$ and let $c\colon E(G)\to \binom{C}{\Delta}$ be a $\mu n$-bounded edge set colouring of $G$.
Suppose we choose a random subset $C'$ of $C$ by including every colour $\alpha \in C$ independently with probability $p$.
Then with probability at least $1-\eul^{-\mu^{-0.9} n}$,
the graph $G_{C'}$ is (lower) $(2\eps,p^{\Delta}d)$-regular.
\end{prop}

\proof
We only consider the case when $G$ is $(\eps,d)$-regular.
If $G$ is lower\COMMENT{deleted brackets} $(\eps,d)$-regular, the proof is similar.
Consider sets $S\In V$, $T \In V'$ with $|S|\ge \eps |V|$ and $|T|\ge \eps |V'|$.
Let $X(S,T)=e_{G_{C'}}(S,T)$ denote the number of edges $e$ between $S$ and $T$
such that $c(e)\subseteq C'$.
Clearly, $\expn{X(S,T)}=p^\Delta e_G(S,T) = p^\Delta d |S||T| \pm \eps |S||T|$.
Moreover, $X(S,T)$ is a random variable where a single choice of whether $\alpha\in C'$ changes $X(S,T)$ by at most $\mu n$.
Hence, by Theorem~\ref{thm:McDiarmid},
\begin{align*}
\prob{|X(S,T)-\expn{X(S,T)}|> \eps |S||T|}
\le 2\eul^{-\frac{2(\eps |S||T|)^2}{\mu n \Delta e(G)}}
\le 2\eul^{-\frac{\eps^{10}}{2\mu \Delta} n}.
\end{align*}
Thus, a union bound shows that the activated subgraph is not $(2\eps,p^\Delta d)$-regular with probability at most $2^{4n}2\eul^{-\frac{\eps^{10}}{2\mu \Delta}n}\le \eul^{-\mu^{-0.9} n}$.
\endproof

However, in order to show that the activated subgraph is super-regular, we have to control the degree of every vertex as well. Whilst for sets $S,T$ as above, the activation or deactivation of one colour has only a small effect, this is very different for the degree of a single vertex.
To appreciate this issue, note that there might only be $\mu^{-1}$ colours present at a particular vertex at all.
Hence,
if we suppose the `activated' degrees of the vertices are mutually independent,
then the probability that all `activated' degrees are in the desired range decays exponentially with~$n$.
We note that under the additional assumption that the given colouring is locally $o(n/\log n)$-bounded, we could apply McDiarmid's inequality for each vertex and conclude with a union bound that our desired properties hold with probability close to~$1$.
We remark that in \cite{APS:17}, the authors considered a problem of a similar kind. They prove that if the colours of a proper edge colouring of $K_n$ are activated with probability~$1/2$, say, then with high probability, the activated subgraph has very good expansion properties.

Our efforts in the rest of this section focus on overcoming this local barrier. As a consequence, we cannot hope that a random partition satisfies the desired properties with high probability. The approach we pursue is motivated by the Lov\'asz local lemma, but the Lov\'asz local lemma itself does not seem flexible enough for our purposes. We will discuss this in more detail in the next subsection.

In order to motivate the next result, which allows us to split the degree of each vertex evenly among the colour sets, we note that (the simplified version of) this problem can be phrased as a `column-sparse assignment-packing problem'. Let $G$ be a graph on $n$ vertices and $c\colon E(G)\to C$ a $\mu n$-bounded edge colouring. Our aim is to find a subset $C'$ of `activated' colours such that for every vertex $v\in V(G)$, we have $d_{C'}(v)=d_G(v)/2\pm 2\mu n$, say.
Define the matrix $\hat{A}$ whose rows are indexed by $V(G)$ and whose columns are indexed by $C$, where $\hat{A}_{v,i}:=d^i(v)$ is the number of edges at $v$ coloured~$i$.
Let $$D:=\max_{i\in C}\sum_{v\in V(G)}d^i(v)$$ denote the maximum column sum of $\hat{A}$.
Let $\hat{c}$ denote the vector indexed by $V(G)$ with $\hat{c}_v:=d_G(v)/2$. Define 
\begin{align*}
	\text{$A:=
	\left(
	\begin{array}{c}
	\hat{A}\\
	-\hat{A}
	\end{array}
	\right)$
	and $c:=
	\left(
	\begin{array}{c}
	\hat{c}\\
	-\hat{c}
	\end{array}
\right)$.}
\end{align*}
Let $x$ be a vector indexed by~$C$.
Clearly, the linear system $Ax=c$ has the solution $x_i^\ast=1/2$ for all $i\in C$. The rounding theorem of Karp et al.~\cite{KLRTVV:87} implies the existence of an integral vector $x$ such that $|x_i-x^\ast_i|<1$ for all $i\in C$ and such that $(\hat{A}x)_v= (\hat{A} x^\ast)_v \pm D$ for all $v\in V(G)$. Clearly, $x$ is a $0/1$-vector.
Let $C'$ be the subset of $C$ which represents the $1$-coordinates of $x$. Then we have $d_{C'}(v)=(\hat{A}x)_v=\hat{c}_v\pm D=d_G(v)/2\pm D$ for all $v\in V(G)$. Crucially, since $c$ is $\mu n$-bounded, we have that $D\le 2\mu n$, as desired.

This is good news for our problem. However, in order to prove Lemma~\ref{lem:separate colours}, we need that a random partition of $C$ satisfies the desired degree properties with high enough probability so that after applying Proposition~\ref{prop:regularity slicing}, we can use a union bound to conclude that the sliced pairs are super-regular.
We will deal with this in the next two subsections.
The following corollary is the result of these efforts.
With this tool at hand, we can deduce that the sliced pairs do not only inherit regularity, but also super-regularity, all we need to prove Lemma~\ref{lem:separate colours}.

\begin{cor} \label{cor:degree split advanced}
Suppose $1/n\ll \mu \ll 1/t,d,1/\Delta,1/r$ and let $\kappa:=4\sqrt{t^\Delta d^{-1}\mu}$.
Let $G$ be a graph on $n$ vertices and let $c\colon E(G)\to \binom{C}{\Delta}$ be $\mu n$-bounded.
Suppose for every vertex $v\in V(G)$, we are given a partition $S^v_1,\ldots,S_{m_v}^v$ of $N_G(v)$ into at most $r$ sets of size at least $drn$ each.\COMMENT{Allows us to use same $\kappa$ as in Lemma~\ref{lem:degree split advanced}}
Suppose we split $C$ randomly into $C_1,\dots,C_t$, by independently assigning $\alpha \in C$ to $C_\ell$ with probability $1/t$.
Then with probability at least $\frac{\kappa}{3}(1+\kappa)^{-|C|}$, the following holds:
$d_{C_\ell}(v,S_i^v)= t^{-\Delta} |S_i^v| \pm  (2t)^\Delta \kappa r n$
for all $v\in V(G)$, $i\in [m_v]$, and $\ell\in [t]$.
\end{cor}

\lateproof{Lemma~\ref{lem:separate colours}}
We only consider the case when $G[V_i,V_j]$ is in fact $(\eps,d)$-(super)-regular.
If $G[V_i,V_j]$ is only lower $(\eps,d)$-(super)-regular, the proof is similar.
If $2\eps \ge d$, there is nothing to prove, so let us assume that $\eps\le d/2$.

In order to apply Proposition~\ref{prop:regularity slicing} and Corollary~\ref{cor:degree split advanced}, we need the edge set colouring to be uniform. We thus add a set of new `dummy' colours $C'$ with $|C'|= \lceil \mu^{-1}\Delta n \rceil$. Let $$C^\ast:=C \cupdot C'$$ be the extended colour set and note that $|C^\ast|\le 2\mu^{-1.3}n$.
Let $c'\colon E(G) \to \binom{C'}{\Delta}$ be any $\mu n$-bounded edge set colouring. Clearly, $c'$ exists since $\frac{\Delta e(G)}{\mu n} \le |C'|$.

Now, let $c^*\colon E(G) \to \binom{C^\ast}{\Delta}$ be such that
\begin{align}
	\text{$c^* (e) \In c(e) \cup c'(e)$ and $c(e) \In c^* (e)$ for all $e\in E(G)$.} \label{filled colour sets}
\end{align}
Clearly, $c^*$ is $\mu n$-bounded.

Split $C^\ast$ randomly into $C_1^*,\dots,C_t^*$, by independently assigning $\alpha\in C^\ast$ to $C_\ell^*$ with probability $1/t$.
By Proposition~\ref{prop:regularity slicing},
we conclude that with probability at least $1-tr^2\eul^{-\mu^{-0.9} n}$,
the graph $G_\ell^{ij}:=G_{C_\ell^*}[V_i,V_j]$ is $(2\eps,d/t^\Delta)$-regular for all $\ell \in [t]$ and all $ij \in E(R)$.

In order to apply Corollary~\ref{cor:degree split advanced}, define $G^\ast:=\bigcup_{ij\in E(R_S)}G[V_i,V_j]$. Consider $i\in[r]$ and $v\in V_i$. For $j\in N_{R_S}(i)$, define $S^v_j:=N_G(v)\cap V_j$. Note that $(S^v_j)_{j\in N_{R_S}(i)}$ partitions $N_{G^\ast}(v)$, and that $|S^v_j|= (d\pm \eps)|V_j|$ since $G[V_i,V_j]$ is $(\eps,d)$-super-regular. In particular, $|S^v_j|\ge d|V_j|/2 $.

By Corollary~\ref{cor:degree split advanced} (with $G^\ast$, $d/4r^2$ playing the roles of $G,d$), we have with probability at least $\mu (1+ \sqrt{\mu})^{-|C^\ast|} $ that
\begin{align}
d_{G^{\ast}_{C_\ell^*}}(v,V_j)= t^{-\Delta}(d\pm \eps)|V_j| \pm  \mu^{1/3} n  \label{splitted degree}
\end{align}
for all $ij\in E(R_S)$, $v\in V_i$, and $\ell\in [t]$.

Since $|C^*|\le 2\mu^{-1.3}n$, we have\COMMENT{with $(1+ \sqrt{\mu})^{-|C^*|}\ge \eul^{\mu^{1/2}(-2\mu^{-1.3}n)}=\eul^{-2\mu^{-0.8}n}$}
\begin{align*}
	\mu(1+ \sqrt{\mu})^{-|C^*|}-tr^2\eul^{-\mu^{-0.9} n}
	\geq \mu \eul^{-2\mu^{-0.8} n }-tr^2\eul^{-\mu^{-0.9} n}
	>0.
\end{align*}
Therefore, there exists a partition $(C_\ell^*)_{\ell\in[t]}$ of $C^\ast$ with the above properties. From \eqref{splitted degree}, we can infer that $G_\ell^{ij}$ is $(2\eps,d/t^\Delta)$-super-regular for all $\ell \in [t]$ and all $ij\in E(R_S)$.
Finally, define $C_\ell:=C\cap C_\ell^*$ for all $\ell\in [t]$. Clearly, $(C_\ell)_{\ell\in[t]}$ is a partition of $C$, and we have $G_\ell^{ij}\In G_{C_\ell}[V_i,V_j]$ by~\eqref{filled colour sets} for all $\ell \in [t]$ and all $ij \in E(R)$.
\endproof

\subsection{The partial resampling algorithm}
The proof of Corollary~\ref{cor:degree split advanced} is based on the partial resampling algorithm introduced by Harris and Srinivasan, which is a relative of the Moser-Tardos resampling algorithm. These algorithms are set in the so-called `variable model' of the Lov\'asz local lemma.
Let $X_1,\dots,X_N$ be independent random variables, and suppose we want to avoid a set of bad events $\cB$.
For every $A\in \cB$, let $vbl(A)$ be the set of variables which determine $A$.
This naturally defines a dependency graph $\Gamma$ on $\cB$ as follows: $AB \in E(\Gamma)$ if $vbl(A)\cap vbl(B)\neq\emptyset$.
The (symmetric) Lov\'asz local lemma implies that if each bad event has probability at most $p$ and $4p\Delta(\Gamma)\le 1$, then there exists an assignment of the variables $X_1,\dots, X_N$ which avoids all bad events.
In their well-known paper, Moser and Tardos~\cite{MT:10} provided an algorithmic version thereof. Their so-called resampling algorithm is as simple as it could be: Start with a random assignment of the variables $X_1,\dots,X_N$. As long as some bad event $A\in \cB$ is still true, resample all the variables in $vbl(A)$ (according to their respective distribution).
Clearly, if the algorithm terminates, then this yields an assignment of the variables which avoids all bad events. Moser and Tardos showed that, under the assumptions of the Lov\'asz local lemma, the algorithm terminates with probability~$1$.

In our (simplified) setting, let $X_i$ denote the Bernoulli random variable indicating whether colour $i\in C$ is activated or not. For each vertex $v\in V(G)$, let $A_v$ be the bad event that $\sum_{i\in C}d_G^i(v)\mathds{1}_{\Set{X_i=1}} \neq d_G(v)/2\pm \sqrt{\mu}n$, say.
Note that $A_v$ depends on all colours which are present at~$v$. This may lead to a very dense dependency graph, possibly far too dense to apply the Lov\'asz local lemma.

Nonetheless, a lot of these dependencies might be very weak. For example, given a vertex $v$ and colours $i,j\in C$, where $d^i(v)$ is much larger than $d^j(v)$, the dependency of $A_v$ on $X_i$ should intuitively be much more significant than its dependency on~$X_j$.
A variant of the resampling algorithm of Moser and Tardos tailored towards this problem could thus vaguely look as follows: Start with a random assignment of the variables $(X_i)_{i\in C}$. As long as some bad event $A_v$ is still true, choose randomly a subset $C'$ of colours, where a colour is more likely to be chosen if many edges at $v$ are coloured with this colour. Now only resample $(X_i)_{i\in C'}$.

This `partial resampling algorithm' falls into a very general framework, which was introduced by Harris and Srinivasan~\cite{HS:13b,HS:13a,HS:18}. They eliminated the need of a dependency graph by introducing `fractional hitting sets' instead, which capture `how dependent' a bad event $A$ is on a variable $X_i$ (or more generally, a set of variables).

We now introduce this framework, following their exposition in~\cite{HS:18}.
Let $X_1,\dots,X_N$ be independent random variables, where each $X_i$ has a finite set $L_i$ of possible assignments. For $i\in [N]$ and $j\in L_i$, let $p_{i,j}=\prob{X_i=j}$. Hence, $\sum_{j\in L_i}p_{i,j}=1$. Let $(\Omega,\mathbb{P})$ be the respective product probability space.

An ordered pair $(i,j)$ with $i\in [N]$ and $j\in L_i$ is referred to as an \defn{element}, and $\cX$ denotes the set of all elements.
An \defn{atomic event} is a set $Y\In \cX$ such that if $(i,j),(i,j')\in Y$, then $j=j'$. Let $\cA$ denote the family of all atomic events.
A set $\cB\In \cA$ is called a \defn{complex event}.
An \defn{assignment} is an atomic event $A$ such that $|A|=N$.
We say that an assignment $A$ \defn{avoids $\cB\In \cA$} if there is no $B\in \cB$ with $B\In A$.

These definitions naturally relate to the more standard notation for events. For an element $\omega=(j_1,\dots,j_N)$ of the probability space $\Omega$, define $A_{\omega}:=\Set{(1,j_1),\dots,(N,j_N)}$. Then $A_{\omega}$ is an assignment. For an event $\cE\In \Omega$, define $\cB_{\cE}:=\set{A_{\omega}}{\omega \in \cE}$. Clearly, we have $\omega\in \cE$ if and only if $A_\omega$ does not avoid $\cB_{\cE}$. Note that our complex events $\cB_{\cE}$ only contain assignments. In general, it is also possible to deal with complex events which contain atomic events $A$ with $|A|<N$. 

Given a function $\lambda\colon \cX \to \bR$, we use the following shorthand notation: for an element $(i,j)\in \cX$, we write $\lambda_{i,j}:=\lambda((i,j))$, and for a set $Y\In \cX$, we write $$\lambda^Y:=\prod_{(i,j)\in Y}\lambda_{i,j}.$$
For $Y\In \cX$ and $i\in [N]$, we write $Y\sim i$ if there exists some $j\in L_i$ such that $(i,j)\in Y$.

\begin{defin}[fractional hitting set, cf.~{\cite[{Definition~2.1}]{HS:18}}]
Let $Q\colon \cA \to [0,1]$. For a set $B\in \cA$, we say that $Q$ is a \defn{fractional hitting set for $B$} if $Q(\emptyset)=0$ and $$\sum_{Y\In B}Q(Y)\ge 1.$$ For a complex event $\cB\In \cA$, we say that $Q$ is a \defn{fractional hitting set for $\cB$} if it is a fractional hitting set for all $B\in \cB$. Moreover, for an event $\cE\In \Omega$, we say that $Q$ is a \defn{fractional hitting set for $\cE$} if it is a fractional hitting set for $\cB_{\cE}$.
\end{defin}

The intention of this definition is to offer a more flexible notion of dependency within the framework, which finally eliminates the need for a dependency graph.
Suppose we want to find an assignment that avoids the complex events $\cB_1,\ldots,\cB_K$, for which we are given fractional hitting sets $Q_1,\dots,Q_K$, respectively.

The \defn{partial resampling algorithm (PRA)} starts with a random assignment of the~$X_i$. Then it repeats the following, as long as some bad event is true: choose any $B\in \cB_k$ which is true. Now, select randomly a set $Y\In B$, where the probability of selecting $Y$ is given by $\frac{Q_k(Y)}{\sum_{Y'\In B}Q_k(Y')}$. Then resample all $X_i$ with $Y\sim i$ (according to their distribution).
Clearly, if the algorithm terminates, it outputs an assignment of the $X_i$ which avoids all bad events.

The next definition is the final ingredient before we can state the result of Harris and Srinivasan.
Therein, the function $\lambda$ may be thought of as an `inflated' probability vector.

\begin{defin}
Let $Q\colon \cA\to [0,1]$ and $\lambda\colon \cX \to [0,\infty)$. We define
$$\Gamma(Q,\lambda):= \sum_{Y\in \cA}Q(Y)\lambda^Y,$$ and for each $i\in [N]$, we define $$\Gamma_i(Q,\lambda):= \sum_{Y\in \cA \colon Y\sim i} Q(Y)\lambda^Y.$$\COMMENT{Note that we do not have $\Gamma(Q,\lambda)=\sum_{i\in [N]}\Gamma_i(Q,\lambda)$.}
\end{defin}

\begin{theorem}[Harris and Srinivasan {\cite[{Theorem~3.8, see also Equation~(10)}]{HS:18}}] \label{thm:PRA}
\COMMENT{By Equation (10), we can use $\Gamma(Q_k,\lambda)$ instead of $S(\cB_k,Q_k,\lambda)$}
Let $\cB_1,\dots,\cB_K \In \cA$ be complex (bad) events.
Let $\lambda\colon \cX \to [0,\infty)$ be such that $p_{i,j}=\frac{\lambda_{i,j}}{\sum_{j\in L_i}\lambda_{i,j}}$. Suppose we are given fractional hitting sets $Q_1,\dots,Q_K$ for $\cB_1,\dots,\cB_K$, respectively. Assume that the following conditions are satisfied:
\begin{enumerate}[label=\rm{(\roman*)}]
  \item for all $k\in [K]$, we have $\Gamma(Q_k,\lambda)<1$;
  \item for all $i\in [N]$, we have $\sum_{k\in [K]}\frac{\Gamma_i(Q_k,\lambda)}{1-\Gamma(Q_k,\lambda)} +1 \le \sum_{j\in L_i}\lambda_{i,j}$.
\end{enumerate}
Then the expected number of resamplings of the PRA for these parameters is at most $\sum_{(i,j)\in \cX}\lambda_{i,j}$. In particular, the PRA terminates with probability~$1$.
\end{theorem}

\subsection{Proof of Corollary~\ref{cor:degree split advanced}}

Theorem~\ref{thm:PRA} guarantees that there exists an assignment which avoids all of $\cB_1,\dots,\cB_K$. It does however not state a lower bound on the probability (in the space $(\Omega,\mathbb{P})$) that no bad event happens. Fortunately, we can deduce the following corollary, which yields such a lower bound. For this, we use the PRA in an indirect proof. It would be interesting to know whether there exists an appropriately formulated generalized Lov\'asz local lemma with a direct proof.

\begin{cor}\label{cor:PRA explicit prob}
Let $\cE_1,\dots,\cE_K\In \Omega$ be (bad) events.
Let $\lambda\colon \cX \to [0,\infty)$ be defined as $\lambda_{i,j}:=(1+\kappa)p_{i,j}$, where $\kappa\in (0,1)$.\COMMENT{Could even allow another $\kappa_i$ for every $i$.} Suppose there are fractional hitting sets $Q_1,\dots,Q_K$ for $\cE_1,\dots,\cE_K$, respectively, which satisfy the following conditions:
\begin{enumerate}[label=\rm{(\roman*)}]
  \item for all $k\in [K]$, we have $\Gamma(Q_k,\lambda)< 1$;
  \item for all $i\in [N]$, we have $\sum_{k\in [K]}\frac{\Gamma_i(Q_k,\lambda)}{1-\Gamma(Q_k,\lambda)} \le \kappa/2$.
\end{enumerate}
Then $\prob{\bigcap_{k\in[K]}\overline{\cE_k}} \ge \frac{\kappa}{3}(1+\kappa)^{-N}$.
\end{cor}

\proof
For any event $\cE\In \Omega$,
we have that $\prob{\cE}=\sum_{\omega \in \cE}\prob{\Set{\omega}}=\sum_{\omega \in \cE}\prod_{(i,j)\in A_{\omega}}p_{i,j}$.
Therefore,
\begin{align}\label{eq:E*}
	\prob{\bigcap_{k\in[K]}\overline{\cE_k}}=\prob{\cE^\ast}=\sum_{A\in \cB_\ast}\prod_{(i,j)\in A} p_{i,j},
\end{align}
where $\cE^\ast:=\Omega\sm \bigcup_{k\in[K]}\cE_k$, $\cB_\ast:=\cB_{\cE^\ast}$, and $\cB_k:=\cB_{\cE_k}$ for all $k\in [K]$.
Observe that $\bigcup_{k\in [K]\cup \Set{\ast}}\cB_k = \cB_{\Omega}$ and hence there is no assignment which avoids all of $\cB_1,\dots,\cB_K,\cB_\ast$.

Suppose, for a contradiction, that $\prob{\cE^\ast}\le p^\ast:=\frac{\kappa}{3}(1+\kappa)^{-N}$.
Let $Q_\ast\colon \cA \to [0,1]$ be the trivial fractional hitting set for $\cB_\ast$, that is $Q_\ast(Y):=1$ if $Y\in \cB_\ast$ and $Q_\ast(Y):=0$ otherwise.
We emphasize that $\cB_\ast$ only contains assignments
and hence
\begin{align}
\Gamma(Q_\ast,\lambda)
=\sum_{Y\in \cB_\ast}\lambda^Y
= \sum_{Y\in \cB_\ast}\prod_{(i,j)\in Y} (1+\kappa)p_{i,j}
\stackrel{(\ref{eq:E*})}{=} (1+\kappa)^N\prob{\cE^\ast} \le p^\ast (1+\kappa)^N
= \kappa/3. \label{gamma upper bound}
\end{align}
Trivially, for all $i\in[N]$, we have $\Gamma_i(Q_\ast,\lambda)\le\Gamma(Q_\ast,\lambda)\le \kappa/3$, which implies $\frac{\Gamma_i(Q_\ast,\lambda)}{1-\Gamma(Q_\ast,\lambda)}\le \kappa/2$.
We now apply Theorem~\ref{thm:PRA} to the complex events $\cB_1,\dots,\cB_K,\cB_\ast$ with fractional hitting sets $Q_1,\dots,Q_K,Q_\ast$, respectively. For all $k\in[K]\cup \Set{\ast}$, we have $\Gamma(Q_k,\lambda)< 1$ by~(i) and~\eqref{gamma upper bound}. Moreover, for all $i\in [N]$, we have $$\sum_{k\in [K]\cup \Set{\ast}}\frac{\Gamma_i(Q_k,\lambda)}{1-\Gamma(Q_k,\lambda)} +1 \le \kappa/2+\kappa/2 +1 = \kappa +1 =\sum_{j\in[t]}\lambda_{i,j}.$$ Thus, by Theorem~\ref{thm:PRA}, there exists an assignment $A\In \cX$ which avoids all of $\cB_1,\dots,\cB_K,\cB_\ast$. This clearly is a contradiction, which completes the proof.
\endproof

We now use Corollary~\ref{cor:PRA explicit prob} to prove Corollary~\ref{cor:degree split advanced}. To this end, we need to introduce some more notation.

Let $\cS_\Delta^t:=\set{(s_1,\dots,s_t)}{\sum_{j=1}^t s_j=\Delta,s_j\in \bN_0}$ denote the $t$-dimensional discrete simplex.
For $\mathbf{s}=(s_1,\dots,s_t)\in \cS_\Delta^t$, let $\binom{\Delta}{\mathbf{s}}:=\frac{\Delta!}{s_1!\cdots s_t!}$ denote the multinomial coefficient. Recall that
\begin{align}
\sum_{\mathbf{s} \in \cS_\Delta^t}\binom{\Delta}{\mathbf{s}} = t^\Delta. \label{multinomial identity}
\end{align}

Let $G$ be a graph, $c\colon E(G)\to \binom{C}{\Delta}$ and suppose that $C_1,\ldots, C_t$ is a (random) partition of $C$.
For a vertex~$v$ and $\mathbf{s}=(s_1,\dots,s_t)\in \cS_\Delta^t$,
let $d_{\mathbf{s},v}$ be the (random) number of edges~$e$ incident to~$v$ with $|c(e)\cap C_j|=s_j$ for all $j\in[t]$.
For instance, for $\mathbf{s}=(\Delta,0,\dots,0)$, we have $d_{\mathbf{s},v}=d_{G_{C_1}}(v)$.
Observe that if each $i\in C$ is independently assigned to $C_j$ with probability $1/t$, then $\expn{d_{\mathbf{s},v}}=\binom{\Delta}{\mathbf{s}} t^{-\Delta} d_G (v)$.

\begin{lemma} \label{lem:degree split advanced}
Suppose $1/n\ll \mu \ll 1/t,d,1/\Delta$ and let $\kappa:=4\sqrt{t^\Delta d^{-1}\mu}$.
Let $G$ be a graph on $n$ vertices with $\delta(G)\ge dn$, and let $c\colon E(G)\to \binom{C}{\Delta}$ be $\mu n$-bounded.
Suppose we split $C$ randomly into $C_1,\dots,C_t$, by independently assigning $i\in C$ to $C_j$ with probability $1/t$.
Then with probability at least $\frac{\kappa}{3}(1+\kappa)^{-|C|}$, the following holds:
$$d_{\mathbf{s},v}=\binom{\Delta}{\mathbf{s}} t^{-\Delta} d_G (v) \pm  (2t)^\Delta \kappa n$$
for all $v\in V(G)$ and $\mathbf{s}\in \cS_\Delta^t$.
\end{lemma}

We remark that for any $1\le K\le 2^{-\Delta}\kappa^{-1}$, the above probability bound can be improved to $\frac{\kappa}{3K}(1+\kappa/K)^{-|C|}$ on the expense of an extra factor $K$ in the degree error term, with the same proof.\COMMENT{Define $\lambda=(1+\kappa/K)/t$ and apply Corollary~\ref{cor:PRA explicit prob} with $\kappa'=\kappa/K$}

\proof
We may assume that $C=[N]$.
Let $X_1,\dots,X_N$ be independent random variables taking values in $[t]$ uniformly at random. (Hence, the probability space is $\Omega=[t]^N$ and the set of elements is $\cX=[N]\times [t]$.) For $j\in[t]$, define the random sets $C_j:=\set{i\in [N]}{X_i=j}$. Thus, the $C_j$ are as in the statement of the lemma.

For every vertex $v\in V(G)$ and every $\mathbf{s}\in \cS_\Delta^t$, we define $\cE_{v,\mathbf{s}}$ as the (bad) event that $d_{\mathbf{s},v} \ge \binom{\Delta}{\mathbf{s}} t^{-\Delta} d_G (v) + 2^\Delta \kappa n$.
We will show that
\begin{align}
\prob{\bigcap_{v\in V(G),\mathbf{s}\in \cS_\Delta^t}\overline{\cE_{v,\mathbf{s}}}}\ge \frac{\kappa}{3}(1+\kappa)^{-N}.\label{target prob advanced}
\end{align}
 Observe that \eqref{target prob advanced} would complete the proof: if no bad event $\cE_{v,\mathbf{s}}$ happens, we have $d_{\mathbf{s},v} \le \binom{\Delta}{\mathbf{s}} t^{-\Delta} d_G (v) + 2^\Delta\kappa n$ for every vertex $v\in V(G)$ and every $\mathbf{s}\in \cS_\Delta^t$. We can then use \eqref{multinomial identity} to also establish the lower bound
\begin{align*}
d_{\mathbf{s},v} &= d_G(v)-\sum_{\mathbf{s}'\in \cS_\Delta^t\sm\Set{\mathbf{s}}}d_{\mathbf{s}',v}  \ge \left(1-\sum_{\mathbf{s}'\in \cS_\Delta^t\sm\Set{\mathbf{s}}} \binom{\Delta}{\mathbf{s}'} t^{-\Delta}\right)d_G(v) - 2^\Delta(|\cS_\Delta^t|-1)\kappa n \\
              &\ge \binom{\Delta}{\mathbf{s}} t^{-\Delta} d_G (v) - (2t)^\Delta\kappa n.
\end{align*}
We will obtain~\eqref{target prob advanced} by applying Corollary~\ref{cor:PRA explicit prob}.
To this end, we first define fractional hitting sets for our bad events.

\begin{NoHyper}
\begin{step}
Defining fractional hitting sets
\end{step}
\end{NoHyper}

Fix $v\in V(G)$ and $\mathbf{s}\in \cS_\Delta^t$ throughout Step~1.
For $I\in \binom{[N]}{\Delta}$, let $d_I(v)$ denote the number of edges $e$ incident to $v$ with $c(e)=I$.
We say that $Y\in \cA$ is of type $(I,\mathbf{s})$ if $Y=\Set{(i_1,j_1),\dots,(i_\Delta,j_\Delta)}$ with $\Set{i_1,\dots,i_\Delta}=I$ and for each $j\in [t]$,
the number of indices $k\in [\Delta]$ with $j_k=j$ is equal to $s_j$. Let $\cT(I,\mathbf{s})$ denote the set of all $Y\in \cA$ of type $(I,\mathbf{s})$.
This means that the edges~$e$ at~$v$ with $c(e)=I$ contribute to the value of the random variable $d_{\mathbf{s},v}$ if and only if $Y$ is a subset of the assignment
given by the outcome of $X_1,\ldots,X_N$.
More formally, given an assignment $A$, we have
\begin{align}\label{eq:dIAfix}
	\sum_{I\in \binom{[N]}{\Delta},Y\in \cT(I,\mathbf{s})\colon Y\In A}d_I(v)= d_{\mathbf{s},v}(A),
\end{align}
where $d_{\mathbf{s},v}(A)$ refers to the value of $d_{\mathbf{s},v}$ if $X_1,\ldots,X_N$ have outcome $A$.

Observe that
\begin{align}\label{eq:Is}
	|\cT(I,\mathbf{s})|=\binom{\Delta}{\mathbf{s}}.
\end{align}

We make another remark for later use.
We have
\begin{align}\label{eq:dI}
	\sum_{I\in \binom{[N]}{\Delta},Y\in \cT(I,\mathbf{s})}d_I(v)
	\stackrel{\eqref{eq:Is}}{=}\sum_{I\in \binom{[N]}{\Delta}}\binom{\Delta}{\mathbf{s}}d_I(v)
		= \binom{\Delta}{\mathbf{s}} d_G(v).
\end{align}

Next we define a fractional hitting set $Q_{v,\mathbf{s}}$ for the event $\cE_{v,\mathbf{s}}$.
Define $Q_{v,\mathbf{s}}\colon \cA \to [0,\infty)$ for all $I\in \binom{[N]}{\Delta}$ and all $Y\in \cT(I,\mathbf{s})$ as
\begin{align}
Q_{v,\mathbf{s}}(Y):=\frac{d_I(v)}{\binom{\Delta}{\mathbf{s}} t^{-\Delta} d_G (v) + 2^\Delta \kappa n}\label{def fractional advanced}
\end{align}
and set $Q_{v,\mathbf{s}}(Y):=0$ for all other $Y\in \cA$.
We claim that $Q_{v,\mathbf{s}}$ is indeed a fractional hitting set for $\cE_{v,\mathbf{s}}$.
Recall that $\cB_{v,\mathbf{s}}$ is the corresponding set of assignments.
Hence we need to check that $Q_{v,\mathbf{s}}$ is a fractional hitting set for all assignments $A\in \cB_{v,\mathbf{s}}$.

Consider first any assignment $A$.
Then
\begin{align*}
	\sum_{Y\In A}Q_{v,\mathbf{s}}(Y)
	=\sum_{I\in \binom{[N]}{\Delta},Y\in \cT(I,\mathbf{s})\colon Y\In A} Q_{v,\mathbf{s}}(Y)
	%&=\sum_{I\in \binom{[N]}{\Delta},Y \mbox{ of type }(I,\mathbf{s}),Y\In A}\frac{d_I(v)}{\binom{\Delta}{\mathbf{s}} t^{-\Delta} d_G (v) + 2^\Delta \kappa n}\\
	\stackrel{\eqref{eq:dIAfix},\eqref{def fractional advanced}}{=}\frac{d_{\mathbf{s},v}(A)}{\binom{\Delta}{\mathbf{s}} t^{-\Delta} d_G (v) + 2^\Delta \kappa n}.
\end{align*}
Since $d_{\mathbf{s},v}(A) \ge \binom{\Delta}{\mathbf{s}} t^{-\Delta} d_G (v) + 2^\Delta \kappa n$ for all $A\in\cB_{v,\mathbf{s}}$ by definition of $\cE_{v,\mathbf{s}}$,
we conclude that for any $A\in \cB_{\cE_{v,\mathbf{s}}}$, we have $\sum_{Y\In A}Q_{v,\mathbf{s}}(Y)\ge 1$, as required.

\begin{NoHyper}
\begin{step}
Verification of the conditions in Corollary~\ref{cor:PRA explicit prob}
\end{step}
\end{NoHyper}

Now, define $\lambda\colon \cX\to [0,\infty)$ as $\lambda_{i,j}:=(1+\kappa)/t$ for all $(i,j)\in \cX$. We need to check conditions (i) and~(ii) from Corollary~\ref{cor:PRA explicit prob}. In order to check~(i), consider $v\in V(G)$ and $\mathbf{s}\in \cS_\Delta^t$. We have
\begin{eqnarray*}
\Gamma(Q_{v,\mathbf{s}},\lambda)&=& \sum_{Y\in \cA}Q_{v,\mathbf{s}}(Y)\lambda^Y
=   \sum_{I\in \binom{[N]}{\Delta},Y\in \cT(I,\mathbf{s})} Q_{v,\mathbf{s}}(Y)\left(\frac{1+\kappa}{t}\right)^\Delta  \\
&\overset{\eqref{eq:dI},\eqref{def fractional advanced}}{=}& \left(\frac{1+\kappa}{t}\right)^\Delta \cdot \frac{\binom{\Delta}{\mathbf{s}}d_G(v)}{\binom{\Delta}{\mathbf{s}} t^{-\Delta} d_G (v) + 2^\Delta \kappa n} \\
           &\le& \frac{(1+\kappa)^\Delta}{1+2^\Delta \kappa} 
           \le \frac{1+(2^\Delta-1) \kappa}{1+2^\Delta\kappa} 
           = 1-\frac{\kappa}{1+2^\Delta\kappa}
\end{eqnarray*}
\COMMENT{$= \frac{\binom{\Delta}{\mathbf{s}}d_G(v)}{\binom{\Delta}{\mathbf{s}}  d_G (v) + 2^\Delta t^\Delta \kappa n} =  \frac{1}{1+\frac{2^\Delta \kappa n t^\Delta}{\binom{\Delta}{\mathbf{s}} d_G(v)}} $ and $n\ge d_G(v)$ and $t^\Delta\ge \binom{\Delta}{\mathbf{s}}$}
and thus (as $2^\Delta \kappa \le 1 $)
\begin{align}
\Gamma(Q_{v,\mathbf{s}},\lambda)&\le 1-\kappa/2 <1.\label{1st gamma condition new}
\end{align}

We now check (ii). Consider $i\in [N]$. First observe that for any $v\in V(G),\mathbf{s}\in \cS_\Delta^t$, we have
\begin{align}
\Gamma_i(Q_{v,\mathbf{s}},\lambda)
&=\sum_{Y\in \cA\colon Y\sim i }Q_{v,\mathbf{s}}(Y)\lambda^Y = \sum_{I\in \binom{[N]}{\Delta},Y\in \cT(I,\mathbf{s})\colon I\ni i} Q_{v,\mathbf{s}}(Y)\left(\frac{1+\kappa}{t}\right)^\Delta .\label{eq:Gammai}
\end{align}
Crucially, we have for all $i\in [N]$
\begin{align}
\sum_{v\in V(G)}\sum_{I\in \binom{[N]}{\Delta}\colon I\ni i} d_I(v) \le 2\mu n   \label{column sparse}
\end{align}
since $c$ is $\mu n$-bounded and every edge $e\in E(G)$ with $i\in c(e)$ is counted twice on the left hand side. Therefore, for any $\mathbf{s}\in \cS_\Delta^t$, we have
\begin{eqnarray*}\notag
\sum_{v\in V(G)} \sum_{I\in \binom{[N]}{\Delta},Y\in \cT(I,\mathbf{s})\colon I\ni i} Q_{v,\mathbf{s}}(Y)
&\overset{\eqref{eq:Is},\eqref{def fractional advanced},\eqref{column sparse}}{\le}& \frac{\binom{\Delta}{\mathbf{s}} \cdot 2\mu n}{\binom{\Delta}{\mathbf{s}} t^{-\Delta} d_G (v) + 2^\Delta \kappa n} \le  t^\Delta \frac{2\mu }{d}.
\end{eqnarray*}
Using~\eqref{multinomial identity} and \eqref{eq:Gammai}, we deduce that
$$\sum_{v\in V(G),\mathbf{s}\in \cS_\Delta^t}\Gamma_i(Q_{v,\mathbf{s}},\lambda) 
\le |\cS_\Delta^t|\left(\frac{1+\kappa}{t}\right)^\Delta \cdot t^{\Delta} \cdot \frac{2\mu }{d} 
\le \frac{3t^\Delta \mu}{d}.$$
Altogether, with \eqref{1st gamma condition new}, we conclude that
$$\sum_{v\in V(G),\mathbf{s}\in \cS_\Delta^t}\frac{\Gamma_i(Q_{v,\mathbf{s}},\lambda)}{1-\Gamma(Q_{v,\mathbf{s}},\lambda)} \le \frac{2}{\kappa}  \cdot \frac{3t^\Delta \mu}{d} \le \kappa/2$$ by definition of $\kappa$.
Thus, \eqref{target prob advanced} follows from Corollary~\ref{cor:PRA explicit prob}.
\endproof

For the proof of Lemma~\ref{lem:separate colours}, we need the additional flexibility that not only the total degree of a vertex is of interest, but its degree into each of the clusters of a regularity partition. By splitting every vertex we can deduce this easily from Lemma~\ref{lem:degree split advanced}.

\lateproof{Corollary~\ref{cor:degree split advanced}}
We apply Lemma~\ref{lem:degree split advanced} to the graph $G^*$, which is obtained from $G$ by splitting each vertex as follows: for each $v\in V(G)$, let $v_1,\dots,v_{m_v}$ be new vertices. Let $V(G^\ast)=\bigcup_{v\in V(G)}\Set{v_1,\dots,v_{m_v}}$.
For every edge $uv\in E(G)$, there are unique $i\in [m_u]$ and $j\in[m_v]$ such that $v\in S^u_i$ and $u\in S^v_j$. Add the edge $u_iv_j$ to $G^\ast$ and define $c^\ast(u_iv_j):=c(uv)$.
Clearly, $n^\ast:=|V(G^\ast)|\le rn$ and for all $v\in V(G)$ and $i\in[m_v]$, we have $|N_{G^\ast}(v_i)|=|S^v_i|$. In particular, $\delta(G^\ast)\ge drn\ge d n^\ast$. Moreover, $c^\ast\colon E(G^\ast)\to \binom{C}{\Delta}$ is $\mu n^\ast$-bounded.
By Lemma~\ref{lem:degree split advanced} (with $G^\ast$ playing the role of $G$), the following holds with probability at least $\frac{\kappa}{3}(1+\kappa)^{-|C|}$:
$$d_{G^\ast_{C_\ell}}(v_i) = t^{-\Delta} d_{G^\ast} (v_i) \pm  (2t)^\Delta \kappa n^\ast$$
for all $v\in V(G)$, $i\in[m_v]$ and $\ell\in[t]$.
Since $d_{C_\ell}(v,S_i^v)=d_{G^\ast_{C_\ell}}(v_i)$ and $|S_i^v| =d_{G^\ast}(v_i)$, this completes the proof.
\endproof

\section{Proof of the rainbow blow-up lemma} \label{sec:main proof}

In this section, we state and prove the full version of the rainbow blow-up lemma. As indicated in Section~\ref{subsec:blow-up intro}, we also include exceptional vertices and candidate graphs in the statement. 
In~\cite{KSS:97}, the authors provided a version of the blow-up lemma which also allowed for the assignment of certain `candidate sets' for some of the vertices. More precisely, for sufficiently small $\alpha$, it is possible to assign for every cluster $X_i$ to $\alpha |X_i|$ vertices $x$ of $X_i$ a candidate set $A_x\In V_i$, which needs to be of size at least $\beta |V_i|$, such that the blow-up lemma still holds with the additional property that the image of $x$ will lie inside $A_x$.
We keep track of candidate sets by storing the essential information in a `candidacy graph' $A^i$ with bipartition $(X_i,V_i)$, where $N_{A^i}(x)$ encodes the candidate set for~$x$. Thus, we allow candidate sets for all vertices, not only a small fraction. However, we require that the resulting candidacy graph for each cluster is lower super-regular. Clearly this includes the original framework.
More formally, a bipartite graph $A^i$ with bipartition $(X_i,V_i)$ is a \defn{candidacy graph for $(X_i,V_i)$}.
We say that the blow-up instance $(H,G,R,(X_i)_{i\in[r]_0},(V_i)_{i\in[r]_0})$ with candidacy graphs $(A^i)_{i\in[r]}$ is \defn{(lower) $(\eps,d)$-super-regular} if for all $ij\in E(R)$, the bipartite graph $G[V_i,V_j]$ is (lower) $(\eps,d)$-super-regular, and for all $i\in [r]$, the graph $A^i$ is (lower) $(\eps,d)$-super-regular.

Candidate sets are especially helpful if a part of $H$ is already embedded, for example if one has to deal with `exceptional vertices' before the application of the blow-up lemma. For instance, if $x_0$ is an exceptional vertex which already has been assigned its image $v_0$, and $x$ is an $H$-neighbour of $x_0$ to be embedded by the blow-up lemma, then the image of $x$ better lies $N_G(v_0)$. This can be achieved by assigning $x$ a candidate set which is a subset of $N_G(v_0)$. In the rainbow setting, we face the additional challenge that, depending on which of the candidates $v$ we pick as the image of $x$, the edge $vv_0$ will already use a colour which is then forbidden for the rest of the embedding.
In the full statement of our rainbow blow-up lemma, we thus already include the exceptional vertices. For the proof to work, we put the following restrictions on the partial embedding of the exceptional vertices. We will see in Section~\ref{sec:apps} that these criteria can be met in applications easily.

Given a blow-up instance $(H,G,R,(X_i)_{i\in[r]_0},(V_i)_{i\in[r]_0})$ with candidacy graphs $(A^i)_{i\in[r]}$ and an edge set colouring $c\colon E(G)\to 2^C$, we say that the bijection $\phi_0\colon X_0\to V_0$ is \defn{$D$-feasible} if the following hold:
\begin{enumerate}[label=\rm{(EXC\arabic*)}]
\item for all $x\in X_0$, $j\in[r]$ and $y\in N_H(x)\cap X_j$,
we have $N_{A^j}(y)\subseteq N_G(\phi_0(x))$; \label{exc condition:neighbourhoods}
\item for all $j\in [r]$, $x\in X_j$, $v\in N_{A^j}(x)$ and distinct $x_0,x_0'\in N_H(x)\cap X_0$, we have $c(\phi_0(x_0)v)\cap c(\phi_0(x_0')v)=\emptyset$ ; \label{exc condition:rainbow}
\item for all colours $\alpha \in C$, we have $\sum_{x\in X_0} d^{\alpha}_G (\phi_0(x))\cdot d_H(x) \le D$.\COMMENT{Make sure $d_G^{\alpha}(v)$ is defined.} \label{exc condition:high degrees}
\end{enumerate}

Condition \ref{exc condition:neighbourhoods} ensures that whenever we pick for $y\in V(H)\sm X_0$ an image~$v$ from its candidate set, then $v$ is appropriately connected to~$V_0$. Condition \ref{exc condition:rainbow} in turn ensures that no conflict arises from this, i.e.~the star with center $v$ and the images of the neighbours of $v$ in $X_0$ is rainbow.
Condition \ref{exc condition:high degrees} is designed for an application where the exceptional vertices are not required to have bounded degree.
Note that \ref{exc condition:neighbourhoods}--\ref{exc condition:high degrees} are trivially satisfied if $X_0$ is empty.\COMMENT{And $D\ge 0$} Moreover, \ref{exc condition:rainbow} is clearly satisfied if $X_0$ is $2$-independent. Note that if $c$ is $k$-bounded and $\Delta(H)\le \Delta$, then \ref{exc condition:high degrees} holds with $D=2\Delta k$.

Given an edge set colouring $c$ of $G$, we say that a subgraph $H$ is \defn{rainbow} if $c(e)\cap c(e')=\emptyset$ for all distinct $e,e'\in E(H)$. This allows to model slightly more general systems of conflicts. For instance, if $c_1,\dots,c_{\Delta}$ are edge colourings of $G$, we can define $c^\ast(e):=\Set{c_1(e),\dots,c_\Delta(e)}$ for all $e\in E(G)$. If $H$ is rainbow with respect to the edge set colouring $c^\ast$, then $H$ is simultaneously rainbow with respect to all the~$c_i$.

We also note that, even if the given colouring of $G$ is a `normal' edge colouring, the ability to handle edge set colourings is crucial when dealing with exceptional vertices (see Step~\ref*{step:colour splitting} in the proof of Lemma~\ref{lem:blow up matchings}).

We now state the full version of the rainbow blow-up lemma. It clearly implies Lemma~\ref{lem:blow-up simple}. 
\begin{lemma}[Rainbow blow-up lemma]\label{lem:blow up}
Let $1/n \ll \mu , \eps \ll d,1/\Delta$ and $\mu \ll 1/r$.
Let $(H,G,R,(X_i)_{i\in[r]_0},(V_i)_{i\in[r]_0})$ with candidacy graphs $(A^i)_{i\in[r]}$ be a lower $(\eps,d)$-super-regular blow-up instance and assume further that
%$\Delta(H),\Delta(R)\le \Delta$ and $|V_i|=(1\pm \eps)n/r$ for all $i\in[r]$.
\begin{enumerate}[label={\rm (\roman*)}]
	\item $\Delta(R)\le \Delta$ and $d_H(x)\le \Delta$ for all $x\in V(H)\sm X_0$;
	\item $|V_i|=(1\pm \eps)n/r$ for all $i\in[r]$;
	\item for all $i\in[r]$, at most $(2\Delta)^{-4}|X_i|$ vertices in $X_i$ have a neighbour in $X_0$.\label{connected to exceptional}
\end{enumerate}
Let $c\colon E(G)\to 2^C$ be $(\mu n,\Delta)$-bounded.
Suppose a $2\Delta\mu n$-feasible bijection $\phi_0 \colon X_0 \to V_0$ is given.
Then there exists a rainbow embedding $\phi$ of $H$ into $G$ which extends $\phi_0$
such that $\phi(x)\in N_{A^i}(x)$ for all $i\in [r]$ and $x\in X_i$.
\end{lemma}

In the remaining four subsections of this section, we prove Lemma~\ref{lem:blow up}. 
In the first subsection, we deduce Lemma~\ref{lem:blow up} from a similar statement (Lemma~\ref{lem:blow up matchings}),
where we impose considerably stronger assumptions on $G$ and $H$.
In Subsection~\ref{sec:4GL}, we introduce the so-called `Four graphs lemma' from~\cite{RR:99}, which we use in the proof of Lemma~\ref{lem:blow up matchings}.
In Subsection~\ref{sec:conflictM}, we discuss `conflict-free' matchings and see how a recent result of Coulson and Perarnau~\cite{CP:17} via the switching method can be used to find perfect conflict-free matchings in super-regular pairs.
This will be another important ingredient for the proof of Lemma~\ref{lem:blow up matchings}.
In the last subsection we finally prove Lemma~\ref{lem:blow up matchings}.

\subsection{Split into matchings} \label{subsec:matching split}

Our first step in proving Lemma~\ref{lem:blow up} is to reduce it to a similar statement where $H$ is highly structured in the sense that $H$ only induces (perfect) matchings between its partition classes. The main steps of this reduction are essentially the same as in~\cite{RR:99}: we apply the Hajnal--Szemer\'edi theorem to $H^2[X_i]$ for each vertex class $X_i$ to obtain a refined partition of $H$ where every vertex class is now $2$-independent. We refine the partition of $G$ randomly such that super-regularity is preserved. Our reduction is more intricate than the one in~\cite{RR:99} as we also consider candidacy graphs and exceptional vertices and allow the cluster sizes to be slighty unbalanced.

We now state the auxiliary lemma, which we will prove in Section~\ref{subsec:main proof}.
In this subsection, we deduce the rainbow blow-up lemma (Lemma~\ref{lem:blow up}) from Lemma~\ref{lem:blow up matchings}.

\begin{lemma} \label{lem:blow up matchings}
Suppose $1/n \ll \mu \ll \eps \ll d,1/\Delta$ and $\mu \ll 1/r$ such that $r$ divides~$n$.
Let $(H,G,R,(X_i)_{i\in[r]_0},(V_i)_{i\in[r]_0})$ with candidacy graphs $(A^i)_{i\in[r]}$ be an $(\eps,d)$-super-regular blow-up instance and assume further that
%$\Delta(H),\Delta(R)\le \Delta$ and $|V_i|=(1\pm \eps)n/r$ for all $i\in[r]$.
\begin{enumerate}[label={\rm (\roman*)}]
	\item $\Delta(R)\le \Delta$ and $d_H(x,X_0)\leq \Delta$ for all $x\in V(H)\sm X_0$;
	\item $|V_i|=n/r$ for all $i\in[r]$;
	\item for all $ij\in \binom{[r]}{2}$, the graph $H[X_i,X_j]$ is a perfect matching if $ij\in E(R)$ and empty otherwise.
\end{enumerate}
Let $c\colon E(G)\to 2^C$ be $(\mu n,\Delta)$-bounded.
Suppose a $\mu n$-feasible bijection $\phi_0 \colon X_0 \to V_0$ is given.
Then there exists a rainbow embedding $\phi$ of $H$ into $G$ which extends $\phi_0$
such that $\phi(x)\in N_{A^j}(x)$ for all $j\in [r]$ and $x\in X_j$.
\end{lemma}

In the reduction, we will need the classical Hajnal--Szemer\'edi theorem and two facts about regular pairs.

\begin{theorem}[\cite{HS:70}] \label{thm:HS}
Let $G$ be a graph on $n$ vertices with $\Delta(G)< k \le n$. Then $V(G)$ can be partitioned into $k$ independent sets of size $\lfloor \frac{n}{k}\rfloor$ or $\lceil \frac{n}{k}\rceil$.
\end{theorem}

\begin{prop}[{\cite[Fact~2]{RR:99}}] \label{prop:lower to super}
Let $1/n\ll \eps \ll \eps' \ll d$. Let $G$ be a bipartite graph with bipartition $(V_1,V_2)$, where $|V_1|=|V_2|=n$. If $G$ is lower $(\eps,d)$-super-regular, then $G$ contains a spanning subgraph $G'$ which is $(\eps',d^2/2)$-super-regular.
\end{prop}

\begin{fact} \label{fact:regularity}
Let $G$ be a bipartite graph with vertex partition $(A,B)$.
Suppose $G$ is lower $(\eps,d)$-regular and $Y\In B$ with $|Y|\ge \eps|B|$.
Then all but at most $\eps|A|$ vertices of $A$ have at least $(d-\eps)|Y|$ neighbours in~$Y$.
\end{fact}
\COMMENT{
Let $A_1\In A$ be the set of vertices which have less than $(d-\eps)|Y|$ neighbours in~$Y$. Suppose for a contradiction that $|A_1|\ge \eps |A|$. Then, since $G[A,B]$ is lower $(\eps,d)$-regular, $d_G(A_1,Y)\ge d-\eps$ and hence $e_G(A_1,Y)\ge (d-\eps)|A_1||Y|$. On the other hand, by definition of $A_1$, we have $e_G(A_1,Y)< |A_1|(d-\eps)|Y|$, a contradiction.}

\lateproof{Lemma~\ref{lem:blow up}}
The proof divides into four steps. Firstly, we modify the partition of $H$. From this, we obtain a few more exceptional vertices. In Step~2, we then extend $\phi_0$ to the new exceptional set. Subsequently, in Step~3 we refine the partition of $G$ accordingly, and finally we apply Lemma~\ref{lem:blow up matchings}.

Choose a new constant $\eps^\ast$ such that $1/n \ll \mu , \eps \ll \eps^\ast \ll d,1/\Delta$.
Let $n'$ be the largest integer not exceeding $(1-\eps)n/r$ which is divisible by $\Delta^2$. For $i\in[r]$, define $a_i:=|X_i|-n'$. Note that
\begin{align}
0\le a_i \le 3\eps n/r. \label{additional exceptionals}
\end{align}

Let $R^\ast$ be the graph with vertex set $[r]\times [\Delta^2]$ and edges $(i,j)(i',j')\in E(R^\ast)$ whenever $ii'\in E(R)$. Clearly, $\Delta(R^\ast) \le \Delta^3$. Let $r^\ast:=r\Delta^2$ and $n^\ast:=n'r$. So $(1-2\eps)n\le n^\ast \le n$. Later, we will apply Lemma~\ref{lem:blow up matchings} with $n^\ast,r^\ast,R^\ast$ playing the roles of $n,r,R$, respectively.

\begin{NoHyper}
\begin{step}
Refining $H$
\end{step}
\end{NoHyper}

First, we move $a_i$ vertices from each cluster $X_i$ to the exceptional set in order to adjust the sizes.
For vertex sets $U_1\In U_2 \In V(H)$, we say that $(U_2,U_1)$ is \defn{$(2,1)$-independent} if $U_2$ is independent in $H$ and whenever $u,u'\in U_2$ have a common neighbour in $H$, then $u,u'\in U_1$.
\begin{claim}
There is a set $B\In V(H)\sm X_0$ such that $|B\cap X_i|=a_i$ for all $i\in[r]$ and such that $(X_0\cup B,X_0)$ is $(2,1)$-independent in~$H$.
\end{claim}

\claimproof
The set $B$ can be constructed greedily.
Assume that for some $i\in[r-1]_0$, we have found a set $B_i$ such that $|B_i\cap X_j|=a_j$ for all $j\in[i]$ and $(X_0\cup B_i,X_0)$ is $(2,1)$-independent in~$H$. (Note that $B_0=\emptyset$ satisfies this for $i=0$.)
Using \eqref{additional exceptionals} and the fact that $\Delta(R)\le \Delta$ and $\Delta(H)\le \Delta$, it is not hard to see that at most $(2\Delta)^{-2} n/r$ vertices in $X_{i+1}$ are at distance at most $2$ from $X_0\cup B_i$.\COMMENT{Indeed, by~\ref{connected to exceptional}, at most $(2\Delta)^{-4}|X_{i+1}|$ vertices of $X_{i+1}$ are in distance $1$ to $X_0$, and using \eqref{additional exceptionals}, there are at most $\sum_{j\in N_R(i+1)}\Delta a_j \le 3\Delta^2 \eps n/r$ vertices in $X_{i+1}$ which are in distance $1$ of $B_i$. Similarly, there are at most $3\Delta^4 \eps n/r$ vertices in $X_{i+1}$ which are in distance $2$ of $B_i$ (but are not in distance $1$ of $X_0$, those are already counted before). Using \eqref{additional exceptionals} again, there are at most $\sum_{j\in N_R(i+1)}\Delta(2\Delta)^{-4}|X_{j}|  \le (4\Delta)^{-2}(1+\eps) n/r$ vertices in $X_{i+1}$ which are in distance $2$ of $X_0$.}

Now, since $\Delta(H^2[X_{i+1}])\le \Delta^2-1$, there exists a $2$-independent subset of $X_{i+1}$ of size at least $|X_{i+1}|/\Delta^2$. Thus, since $|X_{i+1}|/\Delta^2 - (2\Delta)^{-2} n/r \ge a_{i+1}$, we can pick $a_{i+1}$ $2$-independent vertices from $X_{i+1}$ which are at distance at least $3$ from $X_0\cup B_i$ and add them to $B_i$ to obtain~$B_{i+1}$. Observe that then $(X_0\cup B_{i+1},X_0)$ is $(2,1)$-independent in~$H$.
\endclaimproof

Define $X_0^\ast:=X_0\cup B$ and $X_i':=X_i\sm B$ for $i\in[r]$. Clearly, we have $|X_i'|=n'$.
For all $i\in[r]$, apply the Hajnal--Szemer\'edi theorem (Theorem~\ref{thm:HS}) to $H^2[X_i']$. Note that $\Delta(H^2[X_i'])\le \Delta(\Delta-1)<\Delta^2$. Recall that $n'$ is divisible by $\Delta^2$. Thus, there exists a partition of $X_i'$ into $2$-independent sets $X^\ast_{i,1},\dots,X^\ast_{i,\Delta^2}$ in $H$ of size exactly $n'/\Delta^2$ each.
For all $(i,j),(i',j')\in [r]\times [\Delta^2]$, we have that $H[X^\ast_{i,j},X^\ast_{i',j'}]$ is a matching if $(i,j)(i',j')\in E(R^\ast)$ and empty otherwise.
Thus, $R^\ast$ is a suitable reduced graph for the new partition $(X^\ast_{i,j})_{(i,j)\in [r]\times [\Delta^2]}$ of $V(H)\sm X_0^\ast$, and $X_0^\ast$ is the new exceptional set.
Surely, we can add edges to $H$ to obtain a supergraph $H^\ast$ such that $H^\ast[X^\ast_{i,j},X^\ast_{i',j'}]$ forms a perfect matching for all $(i,j)(i',j')\in E(R^\ast)$
Clearly, any rainbow embedding of $H^\ast$ induces a rainbow embedding of $H$.

We observe that by the $(2,1)$-independence of $(X_0^\ast,X_0)$, we have that for all $i\in[r]$ and $y\in X_i'$, exactly one of the following alternatives applies:
\begin{enumerate}[label=\rm{(\alph*)}]
\item $N_H(y)\cap X_0^\ast=N_H(y)\cap X_0$; \label{exceptional neighbourhood a}
\item $N_H(y)\cap X_0^\ast=\Set{x}$ for a unique $x\in X_j$ with $j\in N_R(i)$; \label{exceptional neighbourhood b}
\item $N_H(y)\cap X_0^\ast=\emptyset$.
\end{enumerate}
Let $W_a$ be the set of all vertices $y\in V(H)\sm X_0^\ast$ for which \ref{exceptional neighbourhood a} applies, and let $W_b$ be the set of all vertices $y\in V(H)\sm X_0^\ast$ for which \ref{exceptional neighbourhood b} applies. For $i\in[r]$, we have
\begin{align}
|W_b\cap X_i|\le \sum_{j\in N_R(i)}\Delta a_j \overset{\eqref{additional exceptionals}}{\le} \Delta^2 \cdot  3 \eps n/r. \label{new neighbours of exceptionals}
\end{align}

\begin{NoHyper}
\begin{step}
Extending $\phi_0$
\end{step}
\end{NoHyper}

We want to find a suitable set $V_0^\ast \supseteq V_0$ and extend $\phi_0$ to a bijection $\phi_0^\ast \colon X_0^\ast \to V_0^\ast$ such that $\phi_0^\ast(x) \in N_{A^j}(x)$ for all $j\in[r]$ and $x\in X_0^\ast\cap X_j$, and such that for all $x\in X_0^\ast\sm X_0$, all $i\in[r]$ and all $y\in N_H(x)\cap X_i$, we have
\begin{align}
|N_{A^i}(y) \cap N_G(\phi_0^\ast(x))| \ge (d-\eps)|N_{A^i}(y)|.\label{new exceptional nbhds}
\end{align}
We can find $V_0^\ast$ and $\phi_0^\ast$ by successively picking a suitable image for each $x\in X_0^\ast\sm X_0$.
Suppose we have already defined images for $Z\In X_0^\ast\sm X_0$ and now want to find a suitable image for $x\in X_0^\ast\sm (X_0\cup Z)$. Let $j\in[r]$ be such that $x\in X_j$. We want to pick $\phi_0^\ast(x)$ from $N_{A^j}(x)$. A vertex $v\in N_{A^j}(x)$ is a suitable image if for all $i\in[r]$ and all $y\in N_H(x)\cap X_i$, we have $d_G(v,N_{A^i}(y))\ge (d-\eps)|N_{A^i}(y)|$.
For fixed $i\in[r]$ and $y\in N_H(x)\cap X_i$, since $ij\in E(R)$ and $G[V_i,V_j]$ is lower $(\eps,d)$-regular, Fact~\ref{fact:regularity}\COMMENT{with $Y=N_{A^i}(y)$} implies that all but at most $\eps |V_j|$ vertices in $V_j$ have at least $(d-\eps)|N_{A^i}(y)|$ neighbours in $N_{A^i}(y)$. Thus, at most $\Delta \eps |V_j|$ vertices of $V_j$ are not suitable images for $x$. Moreover, at most $a_j$ vertices of $V_j$ have already been used as images for the vertices in $Z$. Thus, there exists a suitable image $\phi_0^\ast(x)$ for $x$.

Let $V_0^\ast:=\phi_0^\ast(X_0^\ast)$. Clearly, we have $|V_0^\ast \cap V_i|=a_i$ for all $i\in[r]$. Before we continue to refine the partition of $G$, we have to redefine the neighbourhoods of the vertices in $W_b$ in their respective candidate graphs in order to meet condition~\ref{exc condition:neighbourhoods} for $\phi_0^\ast$ (see Step~4).
Consider $i\in[r]$. Let $A'^i$ be the spanning subgraph of $A^i$ obtained by deleting for every vertex $y\in X_i\cap W_b$ all edges $yv$ for which $v\notin N_G(\phi_0^\ast(x))$, where $x$ is the unique $H$-neighbour of $y$ in $X_0^\ast$ (cf.~\ref{exceptional neighbourhood b}). Note that $A'^i$ still contains vertices from $X_0^\ast$ and $V_0^\ast$. We might as well remove them, but it is more convenient to leave them in for now.

Note that by~\eqref{new exceptional nbhds}, we still have that $d_{A'^i}(y) \ge (d-\eps)|N_{A^i}(y)| \ge d^2 |V_i|/2$ for all $y\in X_i\cap W_b$. Moreover, since $|X_i\cap W_b|\le 4\Delta^2\eps |X_i|$ by \eqref{new neighbours of exceptionals}, it is easy to see that $A'^i$ is still lower $(8 \Delta^2\eps,d^2/2)$-super-regular.\COMMENT{Clearly, the degrees are fine. Moreover, for sets $Y\In X_i$ and $Z\In V_i$ with $|Y|,|Z|\ge 8 \Delta^2\eps |X_i|$, we have $e(Y,Z)\ge e(Y\sm W_b,Z) \ge (d-\eps)|Y\sm W_b||Z| \ge (d-\eps)|Y||Z|/2$.}

\begin{NoHyper}
\begin{step}
Refining $G$
\end{step}
\end{NoHyper}

Let $d':=d^2/2$. For $i\in[r]$, let $V_i':=V_i\sm V_0^\ast$. Clearly, $|V_i'|=|X_i'|=n'$. We will partition each $V_i'$ randomly into $\Delta^2$ equal-sized parts and then match those with the refined parts $X^\ast_{i,j}$ of $X_i'$. In order to obtain super-regular new candidacy graphs, we need to take special care of some vertices.

Consider $i\in [r]$. We say that a vertex $v\in V_i'$ is \defn{good} if
\begin{align}
d_{A'^i}(v,X^\ast_{i,j}) \ge (d'-8 \Delta^2\eps)|X^\ast_{i,j}| \label{def good vertex}
\end{align}
for all $j\in[\Delta^2]$. By Fact~\ref{fact:regularity} and since $A'^i$ is lower $(8 \Delta^2\eps,d')$-regular, all but at most $8 \Delta^4 \eps |V_i|$ vertices of $V_i$ are good. %(Note that if $A'^i$ is complete bipartite, then all vertices are good.)
Let $V_i''$ be the set of good vertices of $V_i'$. Thus,
\begin{align}
|V_i\sm V_i''|\le 8 \Delta^4 \eps |V_i| + a_i \le \sqrt{\eps} |V_i|. \label{most suitable}
\end{align}

We now partition $V_i'$ into equal-sized parts $V^\ast_{i,1},\dots,V^\ast_{i,\Delta^2}$ of size exactly $n'/\Delta^2=n^\ast/r^\ast$ each. First, we take care of the vertices which are not good. For every $v\in V_i'\sm V_i''$, choose $j\in [\Delta^2]$ such that $d_{A'^i}(v,X^\ast_{i,j}) \ge (d'-\sqrt{\eps})|X^\ast_{i,j}|$. Clearly, such an index $j\in [\Delta^2]$ exists since $d_{A'^i}(v)\ge (d'-8 \Delta^2\eps)|X_i| \ge (d'-\sqrt{\eps})|X_i'|$. For $j\in [\Delta^2]$, let $V'_{i,j}$ be the set of all vertices $v\in V_i'\sm V_i''$ which have been assigned to $j$ in this way.

Now, for each $i\in [r]$, let $V''_{i,1},\dots,V''_{i,\Delta^2}$ be a partition of $V_i''$ such that $|V''_{i,j}|=n'/\Delta^2-|V'_{i,j}|$ and such that the following hold:
for all $i\in[r]$, $i'\in N_R(i)$, $v\in V_i$ and $j'\in[\Delta^2]$, we have
\begin{align}
d_G(v,V''_{i',j'}) \ge (d-3 \sqrt{\eps})|V''_{i',j'}|, \label{G partition random degree}
\end{align}
and for all $i\in[r]$, $x\in X_i$ and $j\in[\Delta^2]$, we have
\begin{align}
d_{A'^i}(x,V''_{i,j}) \ge (d'-3 \sqrt{\eps})|V''_{i,j}|. \label{G partition random degree candidates}
\end{align}

That such partitions exist can be seen using a probabilistic argument as follows:
For each $i\in[r]$, let $V''_{i,1},\dots,V''_{i,\Delta^2}$ be a partition of $V_i''$ such that $|V''_{i,j}|=n'/\Delta^2-|V'_{i,j}|$, chosen uniformly at random amongst all such partitions. (We may also assume that the partitions of $V_{i}''$ and $V_{i'}''$, say, are independent, but this is not even necessary.)
In particular, for all $(i,j)\in [r]\times [\Delta^2]$, the set $V''_{i,j}$ is a uniformly random subset of $V_i''$ of size $n'/\Delta^2-|V'_{i,j}|$. Note that for all $i\in[r]$, $i'\in N_R(i)$, $v\in V_i$, we have $$d_G(v,V''_{i'}) \overset{\eqref{most suitable}}{\ge} (d-\eps)|V_{i'}|-\sqrt{\eps} |V_{i'}| \ge (d-2\sqrt{\eps})|V_{i'}''|.$$ Thus, for all $j'\in[\Delta^2]$, we have $\expn{d_G(v,V''_{i',j'})} \ge  (d-2 \sqrt{\eps})|V''_{i',j'}|$. Similarly, for all $i\in[r]$, $x\in X_i$ and $j\in[\Delta^2]$, we have $\expn{d_{A'^i}(x,V''_{i,j})}\ge (d'-2 \sqrt{\eps})|V''_{i,j}|$.
Using a Chernoff-Hoeffding-type bound for the hypergeometric distribution and a union bound, we can see that \eqref{G partition random degree} and \eqref{G partition random degree candidates} are satisfied with positive probability.

Finally, for $(i,j)\in [r]\times [\Delta^2]$, let $V^\ast_{i,j}:= V'_{i,j} \cup V''_{i,j}$.
We claim that for all $(i,j)(i',j')\in E(R^\ast)$, the bipartite graph $G[V^\ast_{i,j},V^\ast_{i',j'}]$ is lower $(3\sqrt{\eps},d)$-super-regular. Indeed, it follows simply from the lower $(\eps,d)$-regularity of the pairs $G[V_i,V_{i'}]$ for $ii'\in E(R)$ that $G[V^\ast_{i,j},V^\ast_{i',j'}]$ is lower $(2\Delta^2\eps,d)$-regular, say.
That every vertex has large enough degree in the respective pair follows from~\eqref{G partition random degree}.

Similarly, for all $(i,j)\in [r]\times [\Delta^2]$, the new candidacy graph $A^{(i,j)}:=A'^i[X^\ast_{i,j},V^\ast_{i,j}]$ is lower $(3\sqrt{\eps},d')$-super-regular. Here, every vertex $x\in X^\ast_{i,j}$ has sufficiently high degree in $V^\ast_{i,j}$ by~\eqref{G partition random degree candidates}. Moreover, all good vertices of $V_i'$ have automatically sufficiently high degree by~\eqref{def good vertex}, and all vertices which are not good have sufficiently high degree in their new candidate graph because of their assignment to a set $V_{i,j}'$.

Finally, using Proposition~\ref{prop:lower to super}, we can transition to a spanning subgraph $G^\ast$ of $G$ such that $G^\ast[V^\ast_{i,j},V^\ast_{i,j'}]$ is $(\eps^\ast,d^2/2)$-super-regular for all $(i,j)(i',j')\in E(R^\ast)$, and for each $(i,j)\in [r]\times[\Delta^2]$, we can transition to a spanning subgraph $A^{\ast(i,j)}$ of $A^{(i,j)}$ such that $A^{\ast(i,j)}$ is $(\eps^\ast,d'^2/2)$-super-regular. Note that we do not delete any edges incident to $V_0^\ast$.

\begin{NoHyper}
\begin{step}
Applying Lemma~\ref{lem:blow up matchings}
\end{step}
\end{NoHyper}

We can now complete the proof. It remains to check that $\phi_0^\ast$ is feasible.

First, consider $x\in X_0^\ast$, $(i,j)\in[r]\times[\Delta^2]$ and $y\in N_{H^\ast}(x)\cap X^\ast_{i,j}$. If $y\in W_a$, then we must have $x\in X_0$ and thus $N_{A^i}(y) \In N_G(\phi_0(x))$ by \ref{exc condition:neighbourhoods} for $\phi_0$. If $y\in W_b$, then we have $N_{A'^i}(y)\In N_G(\phi_0^\ast(x))$ by the definition of $A'^{i}$.
In both cases, we conclude that $N_{A'^i}(y)\In N_{G^\ast}(\phi_0^\ast(x))$ since edges in $G$ are only removed between regular pairs. Since $A^{\ast(i,j)}$ is a subgraph of $A'^i$, we have $N_{A^{\ast(i,j)}}(y) \In N_{G^\ast}(\phi_0^\ast(x))$.
Thus, \ref{exc condition:neighbourhoods} holds for $\phi_0^\ast$.

Condition~\ref{exc condition:rainbow} also holds for $\phi_0^\ast$ because only the vertices in $W_{b}$ gained a new neighbour in~$X_{0}^\ast$. However, as each $y\in W_b$ only has one neighbour in $X_0^\ast$, the condition holds trivially.

Finally, consider \ref{exc condition:high degrees}. Let $\alpha\in C$. Note that $d_{H^\ast}(x)=d_H(x)\le \Delta$ for all $x\in X_0^\ast\sm X_0$. Thus, $$\sum_{x\in X_0^\ast\sm X_0} d_{G^\ast_\alpha} (\phi^\ast_0(x))\cdot d_{H^\ast}(x) \le 2\Delta \mu n$$
and hence $\sum_{x\in X_0^\ast} d_{G^\ast_\alpha} (\phi^\ast_0(x))\cdot d_{H^\ast}(x) \le 2\Delta \mu n + 2\Delta \mu n \le \mu^{0.9}n^\ast$.

Therefore, $\phi_0^\ast$ is $\mu^{0.9}n^\ast$-feasible.
Now apply Lemma~\ref{lem:blow up matchings} as follows:

\medskip
{
\noindent
{
\begin{tabular}{c|c|c|c|c|c|c|c|c|c|c|c|c|c}
 $n$ & $\mu$ & $\eps$ & $d$ & $\Delta$ & $r$ & $H$ & $G$ & $R$ & $X_i$ & $V_i$ & $A^i$ & $c$ & $\phi_0$   \\ \hline
$n^\ast$ & $\mu^{0.9}$ & $\eps^\ast$ & $d^4/8$ & $\Delta^3$ & $r^\ast$ & $H^\ast$ & $G^\ast$ & $R^\ast$ & $X_{i,j}^\ast$ & $V^\ast_{i,j}$ & $A^{\ast (i,j)}$ & $c{\restriction_{E(G^\ast)}}$ & $\phi^\ast_0$
\end{tabular}
}
}

\medskip
This yields a rainbow embedding $\phi$ of $H^\ast$ into $G^\ast$ which extends $\phi^\ast_0$
such that $\phi(x)\in N_{A^{\ast (i,j)}}(x)$ for all $(i,j)\in [r] \times [\Delta^2]$ and $x\in X_{i,j}^\ast$.
Clearly, $\phi$ is thus a rainbow embedding of $H$ into $G$ which extends $\phi_0$
such that $\phi(x)\in N_{A^i}(x)$ for all $i\in [r]$ and $x\in X_i$.
\endproof

\subsection{The Four graphs lemma}\label{sec:4GL}
In this subsection,
we state the so-called `Four graphs lemma' due to R\"odl and Ruci\'nski,
which is an important ingredient in the proof of the blow-up lemma.

Let $V_1,V_2,V_3$ be three disjoint sets of size $n$.
Suppose $G_{i,j}$, for $1\leq i<j\leq 3$,
is an $(\eps,d_{i,j})$-super-regular bipartite graph with vertex partition $(V_i,V_j)$
and for all $v_1v_2\in E(G_{1,2})$ (with $v_1\in V_1$), we have
\begin{align}
	|N_{G_{1,3}}(v_1)\cap N_{G_{2,3}}(v_2)| = (d_{1,3}d_{2,3}\pm \eps)n. \label{four graphs condition}
\end{align}
In such a case we say that the triple $(G_{1,2},G_{1,3},G_{2,3})$ is $(\eps,d_{1,2},d_{1,3},d_{2,3})$-regular.
For a perfect matching $\sigma\colon V_1\to V_2$ of $G_{1,2}$,
let $A_\sigma$ be the spanning subgraph of $G_{1,3}$
such that $v_1v_3\in E(A_\sigma)$ (with $v_i\in V_i$ for $i\in \{1,3\}$)
if $\sigma(v_1)v_3\in E(G_{2,3})$.

\begin{lemma}[Four graphs lemma, see~\cite{RR:99}] \label{lem:four graphs}
Suppose $1/n\ll a \ll \eps \ll \eps' \ll d_{1,2},d_{1,3},d_{2,3}$.
Suppose $(G_{1,2},G_{1,3},G_{2,3})$ is $(\eps,d_{1,2},d_{1,3},d_{2,3})$-regular
on the vertex set $V_1\cup V_2\cup V_3$.
Suppose $\sigma\colon V_1\to V_2$ is a perfect matching of $G_{1,2}$ drawn uniformly at random,
then
\begin{align*}
	\prob{A_\sigma \mbox{ is }(\eps',d_{1,3}d_{2,3})\mbox{-super-regular}}\ge 1-(1-a)^n.
\end{align*}
\end{lemma}

The following proposition is useful when the graphs $G_{i,j}$ are super-regular, but \eqref{four graphs condition} is not satisfied. One can always delete a small fraction of the edges of $G_{1,2}$ such that \eqref{four graphs condition} is then satisfied. The proof is based on standard regularity methods and thus omitted. (See also Fact~1 in~\cite{RR:99} for a very similar statement.)

\begin{prop}\label{prop:four graphs condition}
Let $1/n \ll \eps \ll d,1/k$. Let $V_1,\dots,V_{k+2}$ be disjoint vertex sets of size~$n$. Suppose that $G_{1,2}$ is an $(\eps,d_{1,2})$-super-regular bipartite graph with bipartition $(V_1,V_2)$, where $d_{1,2}\ge d$, and for all $i\in [2]$ and $j\in \Set{3,\dots,k+2}$, $G_{i,j}$ is an $(\eps,d_{i,j})$-super-regular bipartite graph with bipartition $(V_i,V_j)$, where $d_{i,j}\ge d$.
Let $G_{1,2}'$ be the spanning subgraph of $G_{1,2}$ consisting of those edges $v_1v_2\in E(G_{1,2})$ (with $v_1\in V_1,v_2\in V_2$) which satisfy the following: $|N_{G_{1,j}}(v_1)\cap N_{G_{2,j}}(v_2)| = (d_{1,j}d_{2,j}\pm \eps)n$ for all $j\in \Set{3,\dots,k+2}$.
Then $G_{1,2}'$ is still $(2k\sqrt{\eps},d_{1,2})$-super-regular.
\end{prop}

We close this subsection by giving some intuition as to how the Four graphs lemma is applied. Suppose that we embed cluster $X_i$ into $V_i$ by choosing a perfect matching $\sigma$ in a suitable candidacy graph $A^i$ with bipartition $(X_i,V_i)$. Suppose we have not embedded cluster $X_j$ yet. In order to proceed with the embedding, we need to update the candidacy graph $A^j$. This involves the graphs $H[X_i,X_j]$, $A^i$, $A^j$, and $G[V_i,V_j]$. More precisely, suppose that $x\in X_j$ and that $v\in V_j$ is a candidate for $x$ before the embedding of $X_i$. After embedding $X_i$, vertex $v$ remains a valid candidate for $x$ if and only if it is suitably connected to the images of the neighbours of $x$ which are already embedded. Now, if $H[X_i,X_j]$ is a perfect matching, we can define the bijection $\pi\colon X_j\to X_i$ where $\pi(x)$ is the unique $H$-neighbour of $x$ in $X_i$. With this notation, $v$ remains a valid candidate for $x$ if and only if $\sigma(\pi(x))v\in E(G[V_i,V_j])$. We can now identify $X_j$ with $X_i$ according to~$\pi$. More precisely, define the graph $P$ with bipartition $(X_i,V_j)$ which is isomorphic to $A^j$, where $\pi(x)$ plays the role of $x$. We are now left with three vertex sets $X_i,V_i,V_j$ and the three graphs $A^i,P,G[V_i,V_j]$. In this setting, we can apply the Four graphs lemma, which yields a super-regular spanning subgraph $P_\sigma$ of $P$, which, via $\pi$, translates back to the updated candidacy graph for $(X_j,V_j)$.

\subsection{Conflict-free perfect matchings}\label{sec:conflictM}

A \defn{system of conflicts} in a graph $G$ is a set $\cF$ of unordered pairs of edges. If $\Set{e,f}\in \cF$, we say that $e$ and $f$ \defn{conflict}. We say that $\cF$ is $k$-bounded if every edge is contained in at most $k$ conflicts. A subgraph $H$ is \defn{conflict-free} if no two edges of $H$ conflict.

Let $G$ be a bipartite graph with vertex classes $A,B$ such that $|A|=|B|=n$. Suppose $M$ is a perfect matching of $G$ and $e=a_1b_1\in M$. An edge $ab\in E(G)\sm M$ is \defn{$(e,M)$-switchable} if $a_1b_2,a_2 b_1\in E(G)$, where $a_2$ and $b_2$ are matched to $b$ and $a$ by $M$, respectively, i.e.~$a_2b,ab_2\in M$.

The following lemma from~\cite{CP:17} is another important tool in our proof. Its proof is based on the Lopsided Lov\'asz local lemma. In~\cite{CP:17}, it is used to show the existence of conflict-free perfect matchings in Dirac bipartite graphs.

\begin{lemma}[Coulson and Perarnau~{\cite[{Lemma~6}]{CP:17}}] \label{lem:switchings}
Suppose $1/n\ll \mu \ll \gamma$. Let $G$ be a bipartite graph with bipartition $(A,B)$ such that $|A|=|B|=n$. Assume that $G$ has at least one perfect matching, and for every perfect matching $M$ of $G$ and for every $e\in M$ there are at least $\gamma n^2$ edges in $G$ that are $(e,M)$-switchable. Then, given any $\mu n$-bounded system of conflicts for $G$, a uniformly chosen perfect matching of $G$ is conflict-free with probability at least $\eul^{-\mu^{1/2} n}$.
\end{lemma}
%\COMMENT{in old version only stated as existence result, but probability bound follows from proof. Every bad event has probability at most $p:=\frac{4}{\gamma^2 n^2}$. The number of bad events $|\cE|=|Q|$ is at most $\mu n^3$ ($n^2$ choices for the first edge, then at most $\mu n$ choices for a conflict edge). Thus, the probability that no bad event happens is at least $(1-2p)^{|\cE|}\ge (1-\frac{8}{\gamma^2 n^2})^{\mu n^3}\to \eul^{-\frac{8\mu n}{\gamma^2}}$ }

We will apply it to find conflict-free perfect matchings of our candidacy graphs.
It is easy to see that there are many switchings in an $(\eps,d)$-super-regular graph.

\begin{prop}\label{prop:many switchings}
Let $1/n \ll \eps \ll d$.\COMMENT{moderate $\eps$ enough here} Let $G$ be a bipartite graph with bipartition $(A,B)$ such that $|A|=|B|=n$ and $G$ is lower $(\eps,d)$-super-regular. Then $G$ has a perfect matching. Moreover, for every perfect matching $M$ of $G$ and every $e\in M$, there are at least $(d-2\eps)^3 n^2$ edges in $G$ that are $(e,M)$-switchable.
\end{prop}

\proof
It is well known that $G$ has a perfect matching. Let $M$ be a perfect matching of $G$ and suppose $e\in M$.
Let $a\in A$, $b\in B$ with $e=ab$. Define $N_a:=N_G(a)$ and $N_b:=N_G(b)$. Moreover, let $N_a'\In A$ be the set of vertices which are matched to the vertices in $N_a$ by $M$, and let $N_b'\In B$ be the set of vertices which are matched to the vertices in $N_b$ by $M$. Clearly, $|N_a'|=|N_a|$ and $|N_b'|=|N_b|$. Note that all edges in $G[N_a',N_b']$ are $(e,M)$-switchable, except those which already belong to~$M$. Since $d_G(a),d_G(b)\ge (d-\eps)n\ge \eps n$, we have $d_G(N_a',N_b') \ge d-\eps$ and hence $e_G(N_a',N_b')-|M| \ge (d-\eps)|N_a'||N_b'|-n \ge (d-2\eps)^3 n^2$.
\endproof

\subsection{Proof of Lemma~\ref{lem:blow up matchings}} \label{subsec:main proof}
We are now ready to prove the auxiliary blow-up lemma. We split the proof into four steps.
In Step 1,
we show that we may assume that $|C|$ is not too large. This is needed for the application of Lemma~\ref{lem:separate colours}.
We embed $H$ into $G$ in a number of rounds, which depends only on $\Delta$;
in particular, all vertices that belong to the same cluster are embedded simultaneously.
As $r$ may be much larger than $\Delta$, we even have to embed several clusters in a single round (cf.~Section~\ref{sec:sketch}).
In Step 2, we reserve for each round a set of colours.
During the embedding procedure, we will in each round only use colours of $G$ which were assigned to this round.
In Step 3, we set up the induction statement
and in Step 4, we perform the induction step.

\lateproof{Lemma~\ref{lem:blow up matchings}}
Let $n':=|V(G)|=n+|V_0|$. By removing isolated vertices from $X_0$ and their images determined by $\phi_0$ from $V_0$,
we can assume that $|V_0|\le \Delta n$ and hence $n'\le (\Delta+1)n$.\COMMENT{Do we need this? $n'$ doesn't really appear again.}
We may clearly assume that $V_0,\dots,V_r$ are independent in $G$.
Let $$T:=\Delta^2+1,\quad d':=\frac{d}{\binom{T+1}{2}^{\Delta^2}}\quad \mbox{and }\eps_0:=2\eps.$$
Choose new constants $a,\eps_1,\dots,\eps_T$ such that $$1/n \ll \mu \ll a \ll \eps \ll \eps_1 \ll \dots \ll \eps_T \ll d,1/\Delta.$$

\begin{NoHyper}
\begin{step}\label{step:linear colours}
Modifying the colouring assignment
\end{step}
\end{NoHyper}

We claim that we may assume that $|C|\leq 7\mu^{-1}\Delta^2 n$.
Roughly speaking, if two colours $\alpha,\beta$ appear on at most $\mu n/2$ edges each, say, then we wish to merge them to a single colour.
Clearly, the new colouring will still be $(\mu n,\Delta)$-bounded and any rainbow embedding of $H$ in the new colouring is also a rainbow embedding in the original colouring. However, we have to be a bit careful not to violate the feasibility of $\phi_0$.
We start with a simple observation:
\begin{align*}
	\sum_{\alpha\in C}\sum_{x\in X_0} d^{\alpha}_G (\phi_0(x))\cdot d_H(x) \le \sum_{x\in X_0} d_H(x) \Delta d_G (\phi_0(x)) \le \Delta n \sum_{x\in X_0}  d_H(x) \leq \Delta^2 n^2.
\end{align*}\COMMENT{Using $\sum_{\alpha\in C}d^{\alpha}_G (v_0) \le \Delta d_G(v_0) \le 2\Delta n$. Here I thought we get an extra $\Delta+1$ if $V_0$ is not independent. Of course we can assume that $V_0$ is independent. $\sum_{x\in X_0}  d_H(x)\le \sum_{x\in V(H)\sm X_0}d_H(x,X_0) \le n\Delta$ }
Hence there are at most $2\mu^{-1}\Delta^2 n$ colours $\alpha \in C$
for which $\sum_{x\in X_0} d^{\alpha}_G (\phi_0(x))\cdot d_H(x)\geq \mu n/2$ holds.
Let us call these colours \emph{critical}. In order to ensure that~\ref{exc condition:high degrees} still holds, we will not merge colours where one is critical.

We say $\alpha,\beta\in C$ \emph{block} each other if there exist distinct $x,x'\in X_0, x''\in N_H(x)\cap N_H(x')$ and $v\in N_{A^i}(x'')$ for some $i\in [r]$
such that $\alpha\in c(\phi_0(x)v)$ and $\beta\in c(\phi_0(x')v)$. For~\ref{exc condition:rainbow} to be preserved, we must not merge colours which block each other.
Next we seek an upper bound on the number of pairs $\alpha,\beta$ that block each other.
For any $x,x'\in X_0, x''\in N_H(x)\cap N_H(x')$,
there are clearly at most $n$ vertices $v\in N_{A^i}(x'')$
and for each such $v$, there are at most $\Delta^2$ such pairs.
We claim that there are at most $\Delta^2 n$ choices for $x,x',x''$ as above.
Indeed, there are clearly at most $n$ choices for $x''$ as $x''\in V(H)\setminus X_0$ and
as $d_H(x'',X_0)\leq \Delta$, we have at most $\Delta^2$ choices for the pair $x,x'$.
Hence there are at most $n \cdot \Delta^2 \cdot n \Delta^2=\Delta^4n^2$
pairs $\alpha,\beta$ that block each other.
Therefore, any set of colours of size, say, $10 \Delta^2 n$ contains a pair $\alpha,\beta$ such that $\alpha,\beta$ do not block each other.

The following observation motivates our discussion above.
Given $\alpha,\beta$ that do not block each other, are both not critical,
and there are at most $\mu n$ edges on which $\alpha$ or $\beta$ appear,
then we may replace every appearance of $\beta$ by $\alpha$ (some colour sets may become smaller)
and the new edge set colouring is still $(\mu n,\Delta)$-bounded and $\phi_0$ is still $\mu n$-feasible.
In addition, any rainbow embedding of $H$ in the new colouring is also a rainbow embedding in the original colouring.

Therefore, we may assume that there are at most $2\mu^{-1}\Delta^2 n+ 10\Delta^2 n \leq 3\mu^{-1}\Delta^2 n$
colours $\alpha\in C$
such that $\alpha$ appears on at most $\mu n/2$ edges.
Hence $(|C|-3\mu^{-1}\Delta^2 n)\cdot\mu n/2 \leq \Delta e(G)\le 2\Delta^2 n^2$,
which implies $|C|\leq 7\mu^{-1}\Delta^2 n$.

\begin{NoHyper}
\begin{step} \label{step:colour splitting}
Colour splitting
\end{step}
\end{NoHyper}

Let $\psi\colon V(R) \to [T]$ be a proper vertex colouring of $R^2$.\COMMENT{$\Delta(R^2)\le \Delta^2$, so $\chi(R^2)\le \Delta^2+1$.} Moreover, set $\psi(0):=0$. For $t \in[T]_0$, let
$$J_t:=\psi^{-1}(t) \text{ and } J_t^\ast:=\bigcup_{\ell\in [t]_0}J_{\ell}.$$
Note that the sets $(J_t)_{t\in [T]}$ are $2$-independent in $R$, and $J_0=\Set{0}$.

In round~$t\in [T]$, we will embed all vertices from clusters $X_i$ with $i\in J_t$.
In order to reserve colours for the respective rounds, we first partition $E(R)$ according to the colouring~$\psi$. For $t_1t_2\in \binom{[T]}{2}$, let $E_{t_1t_2}:=\set{ij\in E(R)}{\psi(i)=t_1,\psi(j)=t_2}$.
Clearly, $(E_{t_1t_2})_{t_1t_2\in \binom{[T]}{2}}$ is a partition of $E(R)$ since $\psi(i)\neq \psi(j)$ for all $ij\in E(R)$. We will also reserve colours for the edges that have an endpoint in $V_0$. For $t\in [T]$, define $E_{0t}:=\set{\Set{j,j+r}}{j \in J_t}$.\COMMENT{Explanation?}

Our aim is now to use Lemma~\ref{lem:separate colours} to reserve for every set $(E_{t_1t_2})_{t_1t_2\in \binom{[T]_0}{2}}$ an exclusive set of colours, and sparsify $G$ accordingly. This will also let the neighbourhoods of exceptional vertices shrink, and thus we need to update the candidacy graphs $A^j$. To ensure that the new candidacy graphs are again super-regular, we define an auxiliary colouring of the candidacy graphs and apply Lemma~\ref{lem:separate colours} to the somewhat artificial graph $G^{exc}:=G[V(G)\sm V_0]\cup \bigcup_{j\in[r]}A^j$ which is the union of $G[V(G)\sm V_0]$ and all the candidacy graphs. For $j\in[r]$, let $V_{j+r}:=X_j$. Let $R^{exc}$ be the graph on $[2r]$ which is the union of $R$ and the perfect matching $\set{\Set{j,j+r}}{j\in [r]}$. Note that $R^{exc}$ is a reduced graph for $G^{exc}$, and $(E_{t_1t_2})_{t_1t_2\in \binom{[T]_0}{2}}$ is a partition of $E(R^{exc})$.

We now transfer the colouring of $G[V_0,V(G)\sm V_0]$ onto the candidacy graphs in a natural way: Consider $i\in [r]$, $x\in X_i$ and $v\in V_i$ with $xv\in E(A^i)$. Define
\begin{align}
	c^{exc}(xv):=\bigcup_{x'\in N_H(x)\cap X_0}c(\phi_0(x')v).\label{exceptional colouring lifted}
\end{align}
Note that $|c^{exc}(xv)|\le \Delta^2$ since $c$ is $(\mu n, \Delta)$-bounded and $d_H(x,X_0)\le \Delta$.
We also define $c^{exc}(e):=c(e)$ for all $e\in E(G[V(G)\sm V_0])$. We claim that $c^{exc}$ is $(2\mu n,\Delta^2)$-bounded. Consider any $\alpha \in C$. Since $\phi_0$ is $\mu n$-feasible, we can deduce from~\ref{exc condition:high degrees} that $\alpha$ appears on at most $\mu n$ edges $xv\in \bigcup_{j\in[r]}E(A^j)$. Moreover, there are at most $\mu n$ edges $e\in E(G[V(G)\sm V_0])$ on which $\alpha$ appears.

Recall from Step~\ref*{step:linear colours} that we can assume that $|C|\leq 7\mu^{-1}\Delta^2 n$. We now apply Lemma~\ref{lem:separate colours} (with $G^{exc}$, $c^{exc}$, $R^{exc}$, $\Delta^2$, $\binom{T+1}{2}$ playing the roles of $G,c,R_S,\Delta,t$) to obtain a partition $(C_{t_1t_2})_{t_1t_2\in \binom{[T]_0}{2}}$ of $C$ such that for all $t_1t_2\in \binom{[T]_0}{2}$ and all $ij\in E(R^{exc})$,
there is a subgraph $G^{ij}_{t_1t_2}$ of $G^{exc}_{C_{t_1t_2}}[V_i,V_j]$ which is $(\eps_0,d')$-super-regular. We will only keep those for which $\psi(ij)=t_1t_2$. 
More precisely, for all $t_1t_2\in \binom{[T]}{2}$ and all $ij\in E_{t_1t_2}$, we let $G^{ij} := G^{ij}_{t_1t_2}\In G_{C_{t_1t_2}}[V_i,V_j]$, and for all $1\le t\le T$ and all $j\in J_t$, we let $A_0^j:=G^{0j}_{0t} \In A^j_{C_{0t}}$.\footnote{The `0' in $A^j_0$ refers to the fact that these graphs are candidacy graphs after round 0. Later we will also define $A_t^j$ for $t\geq 1$.}
This means that we reserve the colours in $C_{t_1t_2}$ for the edges in $G[V_i,V_j]$ with $\psi(ij)=t_1t_2$. The colours in the sets $C_{0t}$ are `reserved' for the candidacy graphs, which we now transfer back to the exceptional edges. For all $t\in [T]$ and all $j\in J_t$, define the bipartite graph $G^{0j}$ with bipartition $(V_0,V_j)$ and edge set $E(G^{0j})=\set{v_0v\in E(G)}{v_0\in V_0,v\in V_j,c(v_0v)\In C_{0,t}}$.

We now transition from $G$ to the spanning subgraph $G^\ast$ defined as
 $$G^\ast:= \bigcup_{j\in[r]} G^{0j} \cup  \bigcup_{ij\in E(R)} G^{ij} $$        %\bigcup_{0\le t_1 < t_2 \le T}\bigcup_{ij\in E_{t_1,t_2}} G_{C_{t_1,t_2}}[V_i,V_j], $$
which is the spanning subgraph of $G$ containing those edges of $G$ which are `admissibly' coloured. More precisely, for all $t_1t_2\in \binom{[T]_0}{2}$ and all $ij\in E_{t_1t_2}$, we have
\begin{align}
	\text{$c(vw)\In C_{t_1t_2}$ for all $vw\in G^\ast[V_i,V_j]$.} \label{admissible colours}
\end{align}
The following is the reason why we transferred the colouring of the exceptional edges onto the candidacy graphs before applying Lemma~\ref{lem:separate colours}. Observe that for all $j\in [r]$, $x_0\in X_0$ and $x\in N_H(x_0)\cap X_j$,
we have
\begin{align}
N_{A_0^j}(x)\subseteq N_{G^\ast}(\phi_0(x_0)),\label{exceptional containment updated}
\end{align}
i.e.~the new candidacy graphs satisfy~\ref{exc condition:neighbourhoods} with respect to the sparsified graph $G^{\ast}$.
Indeed, given $j,x_0,x$ as above and $v\in N_{A_0^j}(x)$, we have $xv\in A^j$ and $c^{exc}(xv)\In C_{0,t}$, where $t:=\psi(j)$.
By~\ref{exc condition:neighbourhoods}, we have $\phi_0(x_0)v\in E(G)$.
Moreover, by definition of $c^{exc}$, we have $$c(\phi_0(x_0)v) \overset{\eqref{exceptional colouring lifted}}{\In} c^{exc}(xv)\In C_{0,t}$$ and thus $\phi_0(x_0)v\in E(G^{0j})\In E(G^\ast)$, as claimed.
From now on, we do not need the colouring $c^{exc}$ anymore.

\begin{NoHyper}
\begin{step}
Candidacy graphs and the inductive statement
\end{step}
\end{NoHyper}

We will embed $H$ into $G^\ast$.
For brevity, define for $t\in[T]_0$:
\begin{align*}
\bX_t  &:=\bigcup_{i\in J_t}X_i,            &           \bV_t  &:=\bigcup_{i\in J_t}V_i,              &      C_t &:=\bigcup_{k\in[t-1]_0}C_{kt},       \\
\bX_t^\ast &:=\bigcup_{i\in J_t^\ast}X_i,   &           \bV_t^\ast &:=\bigcup_{i\in J_t^\ast}V_i,       &        C_t^\ast &:= \bigcup_{\ell \in[t]_0} C_\ell,        \\
H_t &:= H[\bX_t^\ast], 										 &           G_t &:= G^\ast[\bV_t^\ast].                      &
\end{align*}
Note that $\bX_t^\ast$ contains $X_0$ and $\bV_t^\ast$ contains $V_0$ for all $t\in [T]_0$.
Moreover, $\bX_0^\ast=\bX_0=X_0$, $\bV_0^\ast=\bV_0=V_0$, and $C_0^\ast=C_0=\emptyset$.
After round $t$, we want to have embedded $H_t$ into $G_t$, only using colours from $C_t^\ast$. Given a partial embedding of $H$ into $G^\ast$, an edge of $G^\ast$ is called \defn{used} if it is the image of an edge of $H$, and a colour is called \defn{used} if some edge is used on which this colour appears.

A bijection $\phi\colon \bX_t^\ast \to \bV_t^\ast$ is \defn{valid} if $\phi{\restriction_{X_0}}=\phi_0$ and for all $j\in J_t^\ast$ and all $x\in X_j$, we have $\phi(x)\in N_{A_0^j}(x)$. In particular, this implies that $\phi(X_j)=V_j$ for all $j\in J_t^\ast$.

The following claim ensures that edges which are embedded in different rounds have automatically distinct colours. Moreover, it also ensures that the edges embedded at one vertex in the same round have distinct colours.
\begin{NoHyper}
\begin{claim} \label{claim:colour separated}
Let $t\in [T]_0$ and suppose that $\phi\colon H_t \to G_t$ is a valid embedding. Then $\phi$ uses only colours from $C_t^\ast$.
Moreover, for all $z\in X_k$, $v\in V_k$ with $zv\in E(A^k)$ and $\ell:=\psi(k) >t$ and all distinct $x,y\in N_H(z)\cap \bX_t^\ast$ such that $\phi(x)v,\phi(y)v\in E(G^\ast)$,
we have $c(\phi(x)v),c(\phi(y)v)\In C_\ell$ and $c(\phi(x)v)\cap c(\phi(y)v) = \emptyset$. 
\end{claim}
\end{NoHyper}

\claimproof
To prove the first part of the claim, assume that $vw\in E(G_t)$. Thus, there are unique $i,j\in J_{t}^\ast$ with $ij\in E_{t_1t_2}$ such that $v\in V_i$ and $w\in V_j$. By~\eqref{admissible colours}, we have $c(vw)\In C_{t_1t_2}$. Since $i,j\in J_{t}^\ast$, we have $t_1,t_2\le t$ and thus $C_{t_1t_2}\In C_t^\ast$, as claimed.

To prove the second part of the claim, suppose that $z,v,\ell,k,x,y$ are as in the statement. There are unique $i,j\in J_t^\ast$ with $x\in X_i$ and $y\in X_j$. Let $\ell_1:=\psi(i)$ and $\ell_2:=\psi(j)$. Without loss of generality, assume that $\ell_1\le \ell_2\le t<\ell$.
By~\eqref{admissible colours}, we have that $c(\phi(x)v)\In C_{\ell_1\ell}$ and $c(\phi(y)v)\In C_{\ell_2\ell}$.
If $\ell_1\neq \ell_2$, then $C_{\ell_1\ell}$ and $C_{\ell_2\ell}$ are disjoint subsets of $C_\ell$, and the claim follows.
If $i,j>0$, this implies that $ik,jk\in E(R)$ and we can conclude that $\ell_1\neq \ell_2$ by definition of $\psi$.
If $i=0$ and $j>0$, we have $\ell_1=0$ and $\ell_2>0$.

The remaining case is when $x,y\in X_0$. But then, since $zv\in E(A^k)$, we have $c(\phi(x)v)\cap c(\phi(y)v) = \emptyset$ by \ref{exc condition:rainbow} since $\phi_0$ feasible.
\endclaimproof

Given $t\in [T]_0$ and any bijection $\phi\colon \bX_t^\ast \to \bV_t^\ast$ with $\phi{\restriction_{X_0}}=\phi_0$ (e.g.~a valid embedding as above), we define (for round $t+1$) updated candidacy graphs for the clusters still to be embedded. Let $j\in [r]\sm J_t^\ast$. For vertices $x\in X_j$ and $v\in V_j$, we say that $v$ is a \defn{candidate for $x$} if $xv \in E(A_0^j)$ and for all $y\in N_{H}(x)\cap \bX_t^\ast$ we have that $\phi(y)v\in E(G^\ast)$.

Let $A_t^j$ be the bipartite auxiliary graph with bipartition $(X_j,V_j)$ and edge set
\begin{align}
E(A_t^j):=\set{xv}{x\in X_j,v\in V_j\mbox{ and }v \mbox{ is a candidate for }x}. \label{def candidate set}
\end{align}
Note that $A_t^j$ depends on the specific choice of $\phi$. We might write $A_t^j(\phi)$ to indicate this dependency, but will just write $A_t^j$ if $\phi$ is clear from the context.
Note that $A_0^j(\phi_0)=A_0^j$ by~\eqref{exceptional containment updated}.

We will prove by induction that the following statement \ind{t} holds for all $t\in[T]_0$.
\begin{itemize}
\item[\ind{t}.] There exists a valid rainbow embedding $\phi_t\colon H_t \to G_t$ such that for all $j\in [r]\sm J_t^\ast$, the candidacy graph $A_t^j(\phi_t)$ is $(\eps_t,d_j)$-super-regular for some $d_j\ge d'^{t+1}$.
\end{itemize}
The instance \ind{0} holds since the graphs $(A_0^j)_{j\in[r]}$ are $(\eps_0,d')$-super-regular.
The instance \ind{T} completes the proof since then $\phi_T$ is a valid rainbow embedding of $H$ into $G$.

\begin{NoHyper}
\begin{step}
The inductive step
\end{step}
\end{NoHyper}

Now, assume the truth of \ind{t} for some $t\in[T-1]_0$, and let $\phi_t$ be as in \ind{t}. We will extend $\phi_t$ to $\phi_{t+1}$ such that \ind{t+1} holds. Any bijection $\sigma\colon \bX_{t+1} \to \bV_{t+1}$ induces a bijection $\phi_{t+1}\colon \bX_{t+1}^\ast \to \bV_{t+1}^\ast$ which extends $\phi_t$ as follows:
\begin{align}
\phi_{t+1}(x):=\begin{cases}
\phi_t(x) & \mbox{if } x\in \bX_{t}^\ast, \\
\sigma(x) & \mbox{if } x\in \bX_{t+1}.
\end{cases} \label{def new phi}
\end{align}
We will choose $\sigma$ as a perfect matching in a suitably defined bipartite graph with bipartition $(\bX_{t+1}, \bV_{t+1})$. The natural choice for this graph is $\bigcup_{i\in J_{t+1}} A_t^i$. Indeed, if we pick for every $x\in \bX_{t+1}$ its image $\sigma(x)$ among the neighbours of $x$ in $A_t^i$, then, since every such neighbour $v$ is a candidate for $x$, we will obtain a valid embedding of $H_{t+1}$ into $G_{t+1}$. However, not every perfect matching $\sigma$ in $\bigcup_{i\in J_{t+1}} A_t^i$ induces a $\phi_{t+1}$ which satisfies \ind{t+1}. For this, we also need to ensure that the embedding is rainbow, and that the new candidacy graphs are again super-regular. In order to show that a suitable $\sigma$ exists, we pick $\sigma$ randomly and show that the desired properties hold with positive probability. We will use Lemma~\ref{lem:switchings} to show that $\phi_{t+1}$ is rainbow (see~Claim~\ref*{claim:conflict free}), and the Four graphs lemma (Lemma~\ref{lem:four graphs}) to show that the new candidacy graphs are again super-regular (see~Claim~\ref*{claim:candidate update}).

We now prepare the application of the Four graphs lemma by finding a suitable spanning subgraph $\tilde{A}_t$ of $\bigcup_{i\in J_{t+1}} A_t^i$ of which we will pick the perfect matching $\sigma$ uniformly at random. (The reason for the transition to $\tilde{A}_t$ is to satisfy condition~\eqref{four graphs condition}.)

Recall that for every $ij\in E(R)$, the graph $H[X_i,X_j]$ is a perfect matching. For $i,j\in V(R)$ with $ij\in E(R)$, we define the bijection $\pi_{j\to i}\colon X_j\to X_i$ where $\pi_{j\to i}(x)$ is the unique $H$-neighbour of $x$ in $X_i$. Clearly, $\pi_{j\to i}^{-1}=\pi_{i \to j}$.

\begin{figure}[t]%
\begin{center}
\begin{tikzpicture}[scale = 0.7,every text node part/.style={align=center}]
\def\ver{0.1} %size of a vertex

\draw[black!20,fill=black!20] (8,-1) rectangle (14,1);
\draw (11,0) node {$G^*[V_i,V_{j}]$};

\draw (16,6) node {$H$};
\draw (16,0) node {$G^*$};

\draw (1,7.5) node {$X_{i'}$};
\draw (5,7.5) node {$X_{i''}$};
\draw (9,7.5) node {$X_{i}$};
\draw (15,7.5) node {$X_{j}$};

\draw (1,-1.5) node {$V_{i'}$};
\draw (5,-1.5) node {$V_{i''}$};
\draw (9,-1.5) node {$V_{i}$};
\draw (15,-1.5) node {$V_{j}$};

\filldraw[fill=black!20, draw=black!50]
(7.3,6)--(8.7,6) -- (8.7,0) -- (7.3,0) -- cycle;

\draw[]
(8,5.8) -- (14,5.8)
(8,6) -- (14,6)
(8,6.2) -- node[anchor=south] {$\pi_{i\rightarrow j}$}(14,6.2);

\draw[fill=black!50,black!50]
(13.5,0) rectangle (14.5,6)
;
\draw (16,3) node {$A_t^{j}$ \color{black!70}{$\to A_{t+1}^j$}};
\draw (6.8,4) node {$\tilde{A}_t^{i}$};
\draw (12.6,3) node {$P_i^{j}$};

\filldraw[fill=black!50, draw=black!50]
(13.5,0) -- (14.5,0.5) -- (8.5,6.5) -- (7.5,6) --cycle;

\draw[thick,->] (14,4)  arc  (90:130:4cm);
\draw (12.5,4.2) node {$\pi_{j \rightarrow i}$};

\draw[very thick, fill=white] 
(0,0) ellipse (1 and 1.5)
(4,0) ellipse (1 and 1.5)
(8,0) ellipse (1 and 1.5)
(14,0) ellipse (1 and 1.5)
(0,6) ellipse (1 and 1.5)
(4,6) ellipse (1 and 1.5)
(8,6) ellipse (1 and 1.5)
(14,6) ellipse (1 and 1.5);

\filldraw[fill=black!50, draw=black!60]
(8,6) -- (7.6,0) -- (8.4,0) -- cycle;
\draw[fill=black!50, thick] (8,0) ellipse (0.4 and 0.6);

%\draw[] (8,6) node {$x$};
\draw[fill] (8,6) circle (\ver)
node[anchor=south west] {$x$}
;

\draw[fill] 
(4,6) circle (\ver)
node[anchor=south east] {$x''$}
(0,6) circle (\ver)
node[anchor=south east] {$x'$}
(4,0) circle (\ver)
node[anchor=south] {$\phi(x'')$}
(0,0) circle (\ver)
node[anchor=south] {$\phi(x')$}
(8,0) circle (\ver)
node[anchor=north] {$v$}
;

\draw[thick]
(0,6) .. controls (2,7) and (6,7) .. (8,6)
(4,6) .. controls (5,6.5) and (7,6.5) .. (8,6)
(0,0) .. controls (2,-1) and (6,-1) .. node[near end,anchor=north] {$\alpha$} (8,0)
(4,0) .. controls (5,-0.5) and (7,-0.5) .. node[anchor=south] {$\beta$} (8,0)

(8,0) -- node[anchor=west] {$\{\alpha,\beta\}$} (8,6)  
;

\draw[ultra thick,->]
(0,4.5) -- node[anchor=west] {$\phi_t$}(0,1.5);
\draw[ultra thick,->]
(4,4.5) -- node[anchor=west] {$\phi_t$} (4,1.5)
;

\draw[ultra thick,->,dashed]
(7.5,4.7) -- node[anchor=east] {$\sigma_i$} (7.5,1.3)
;

\end{tikzpicture}
\end{center}
\caption{A sketch of the induction step. Here, $i\in J_{t+1}$, that is, the cluster $X_i$ in the middle is currently embedded into $V_i$ by choosing a random perfect matching $\sigma_i$ in the candidacy graph $\tilde{A}_t^{i}$. To the left, we have clusters that are already embedded, so $i',i''\in J_t^\ast$. The edge $xv$ receives all colours from edges from $v$ to already embedded neighbours of~$x$. This produces a conflict system for $\tilde{A}_t^{i}$. To the right, we have a cluster which will be embedded in a later round, so $j\in [r]\sm J_{t+1}^*$. The current candidacy graph $A_t^j$ is projected using the perfect matching $H[X_i,X_j]$. The random choice of $\sigma_i$ produces a subgraph $P^j_i(\sigma_i)$ of this projected candidacy graph $P_i^j$, which translates back to the new candidacy graph $A_{t+1}^j$ for $(X_j,V_j)$.}
\label{fig:induction step}
\end{figure}
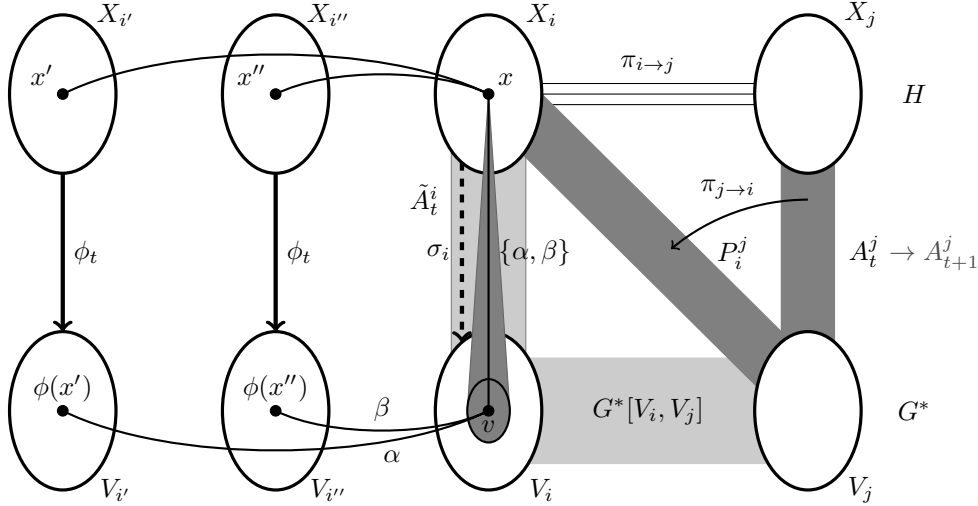

Let $i\in J_{t+1}$. For $j\in N_R(i)\sm J_{t+1}^\ast$, we define the new bipartite auxiliary graph $P^j_i$ with bipartition $(X_i,V_j)$, where $xv\in P^j_i$ if and only if $\pi_{i\to j}(x)v\in E(A_t^j)$. Clearly, $P^j_i$ is isomorphic to $A_{t}^j$ by construction. The reason for the definition of $P^j_i$ is that it forms a super-regular triple with $A_t^i$ and $G^\ast[V_i,V_j]$ on vertex set $X_i\cupdot V_i \cupdot V_j$ (cf.~Figure~\ref{fig:induction step}). We would like to apply the Four graphs lemma with $A_t^i$, $P^j_i$, and $G^\ast[V_i,V_j]$ playing the roles of $G_{1,2}$, $G_{1,3}$, and $G_{2,3}$, respectively.
However, in order to satisfy condition~\eqref{four graphs condition}, we define $\tilde{A}_t^i$ as the spanning subgraph of $A_t^i$ which contains only those edges $xv\in E(A_t^i)$, $x\in X_i,v\in V_i$, that satisfy the following: for all $j\in N_R(i)\sm J_{t+1}^\ast$, we have
\begin{align}
|N_{P^j_i}(x)\cap N_{G^\ast[V_i,V_j]}(v)|=(d_jd'\pm \eps_t)n/r. \label{four graphs condition satisfied}
\end{align}

By Proposition~\ref{prop:four graphs condition}, we still have that for every $i\in J_{t+1}$,
\begin{align}
	\text{$\tilde{A}_t^i$ is $(2\Delta \sqrt{\eps_t},d_i)$-super-regular.} \label{sparsified regular}
\end{align}

Now, define
\begin{align}
\tilde{A}_t:=\bigcup_{i\in J_{t+1}} \tilde{A}_t^i.
\end{align}

Choose $\sigma$ uniformly at random among all perfect matchings of $\tilde{A}_t$. We will show that with positive probability, $\phi_{t+1}$ as defined in~\eqref{def new phi} satisfies \ind{t+1}, which then completes the proof.

Note first that the graphs $(\tilde{A}_t^i)_{i\in J_{t+1}}$ are all vertex-disjoint. In particular, every perfect matching $\sigma$ of $\tilde{A}_t$ is the union of $|J_{t+1}|$ perfect matchings $\sigma_i\colon X_i \to V_i$, where $\sigma_i$ is a perfect matching in $\tilde{A}_t^i$. In particular, every choice of $\sigma$ yields a $\phi_{t+1}$ which is valid. Moreover, a simple but crucial observation is that since $\sigma$ is chosen uniformly at random, each $\sigma_i$ is a uniformly random perfect matching of $\tilde{A}_t^i$.

In order to make sure that $\phi_{t+1}$ is a rainbow embedding, we define a suitable system of conflicts for $\tilde{A}_t$ and then use Lemma~\ref{lem:switchings} to obtain a lower bound on the probability that $\sigma$ is conflict-free. We will subsequently use this lower bound in a union bound, hence a mere existence result would not be sufficient here.
Consider $xv\in E(\tilde{A}_t)$ with $x\in \bX_{t+1},v\in \bV_{t+1}$. There is a unique $i\in J_{t+1}$ with $xv\in E(\tilde{A}_{t}^i)$, i.e.~$x\in X_i$, $v\in V_i$. Let $$F_{xv}:=\set{\phi_t(y)v}{y\in N_{H}(x)\cap \bX_t^\ast}.$$ %Note that $N_{H_{t+1}}(x)\In V(H_t)$ since $J_{t+1}$ is an independent set in $R$, hence $\phi_t(y)$ is defined.
By definition of $A_{t}^i$, we have that $F_{xv}\In E(G^\ast)$.
In particular, $F_{xv}$ is the set of edges of $G^\ast$ which is used in the embedding of $H$ if $x$ is mapped to $v$.
For our embedding to be rainbow, we clearly need that the colour sets\COMMENT{Check terminology in final draft} $(c(e))_{e\in F_{xv}}$ of the edges in $F_{xv}$ are pairwise disjoint, and disjoint from the set of colours used to embed $H_t$.
That this is the case follows from Claim~\ref*{claim:colour separated}. Recall from Claim~\ref*{claim:colour separated} that the embedding of $H_t$ only uses colours from $C_t^\ast$. Moreover, since $\psi(i)=t+1>t$, the second part of Claim~\ref*{claim:colour separated} implies that the sets $(c(e))_{e\in F_{xv}}$ are pairwise disjoint subsets of $C_{t+1}$. Recall that $C_{t+1}\cap C_t^\ast=\emptyset$.
Let $$Z_{xv}:=\bigcup_{e\in F_{xv}} c(e)$$ be the set of all colours which appear on the edges in $F_{xv}$ (cf.~Figure~\ref{fig:induction step}). By the above, $Z_{xv}\In C_{t+1}$.
We define a system of conflicts $\cF$ for $E(\tilde{A}_t)$ as follows: two edges $xv,x'v'\in E(\tilde{A}_t)$ conflict, i.e.~$\Set{xv,x'v'}\in \cF$, if and only if $Z_{xv}\cap Z_{x'v'}\neq \emptyset$.
Hence, if $\sigma$ is a conflict-free perfect matching of $\tilde{A}_t$, then $\phi_{t+1}\colon H_{t+1}\to G_{t+1}$ is a valid rainbow embedding.

\begin{NoHyper}
\begin{claim}\label{claim:conflict free}
The probability that $\sigma$ is conflict-free is at least $\eul^{-\mu^{1/3} n}$.
\end{claim}
\end{NoHyper}

\claimproof
We first claim that $\cF$ is $4\mu \Delta^2 n$-bounded.
Fix any edge $xv\in E(\tilde{A}_t)$.
Clearly, $|Z_{xv}|\le 2\Delta^2$ since $d_H(x,\bX_t^\ast)\le 2\Delta$ and $c$ is $(\mu n,\Delta)$-bounded.
Fix a colour $\alpha\in Z_{xv}$.
Note that if $\alpha \in Z_{x'v'}$, where $x'\in \bX_{t+1}$, $v'\in \bV_{t+1}$ and $x'v'\in E(\tilde{A}_t)$, then $\alpha$ was \defn{forced on $x'v'$} by some edge $v'v''\in E(G^\ast)$ with $v''\in \bV_{t}^*$, $\phi_t^{-1}(v'')x'\in E(H)$ and $\alpha\in c(v'v'') $.
Note that if $v''\in \bV_{t}^*\sm V_0$, then $x'':=\phi_t^{-1}(v'')\in V(H)\sm X_0$ satisfies $d_H(x'',\bX_{t+1})=1$ by definition of $\psi$. Thus, $v'v''$ forces $\alpha$ on at most one edge $x'v'\in E(\tilde{A}_t)$. Since there are at most $\mu n$ edges $v'v''$ with $\alpha\in c(v'v'')$, we conclude that $\alpha$ is forced on at most $\mu n$ edges $x'v'\in \tilde{A}_t$ in this way.
On the other hand, if $v''\in V_0$, then $v'v''$ forces $\alpha$ on at most $d_H(\phi_0^{-1}(v''),\bX_{t+1})\le d_H(\phi_t^{-1}(v''))$ edges $x'v'\in \tilde{A}_t$. Since there are at most $d_G^\alpha (v'')$ vertices $v'$ in $\bV_{t+1}$ such that $\alpha\in c(v'v'')$, the total number of edges $x'v'\in \tilde{A}_t$ on which $\alpha$ is forced by some $v'v''$ with $v''\in V_0$, can be bounded from above by
$$\sum_{v''\in V_0} d^{\alpha}_G (v'') \cdot d_H(\phi_0^{-1}(v''),\bX_{t+1}) \le  \sum_{x''\in X_0} d^{\alpha}_G (\phi_0(x''))\cdot d_H(x'') \overset{\ref{exc condition:high degrees}}{\leq} \mu n.$$
Therefore, $\alpha$ appears on at most $2\mu n$ edges in $\tilde{A}_t$
and hence $\cF$ is $4\mu \Delta^2 n$-bounded.

By~\eqref{sparsified regular} and Proposition~\ref{prop:many switchings} (applied for each $i\in J_{t+1}$), there exists a perfect matching of $\tilde{A}_t$.
Let $M$ be any perfect matching of $\tilde{A}_t$, and $e\in M$. We need to show that there are many $(e,M)$-switchable edges.
Suppose that $e\in E(\tilde{A}_t^i)$ for $i\in J_{t+1}$. Let $M_i$ be the perfect matching of $\tilde{A}_t^i$ induced by $M$. Clearly, every $(e,M_i)$-switchable edge is also $(e,M)$-switchable. By~\eqref{sparsified regular}, we have that $\tilde{A}_t^i$ is $(2\Delta \sqrt{\eps_t},d_i)$-super-regular. By Proposition~\ref{prop:many switchings} (with $n/r,2\Delta \sqrt{\eps_t},d'^{t+1}$ playing the roles of $n,\eps,d$), there are at least $\frac{d'^{3(t+1)}}{2} (n/r)^2$ edges in $\tilde{A}_t^i$ that are $(e,M_i)$-switchable, and thus $(e,M)$-switchable.

Hence, by Lemma~\ref{lem:switchings} (with $|J_{t+1}|n/r, 4\mu r\Delta^2/|J_{t+1}| , d'^{3(t+1)}/2|J_{t+1}|^2$ playing the roles of $n,\mu,\gamma$),\COMMENT{$n^\ast=|J_{t+1}|n/r$. Need that $4\mu \Delta^2 n \le \mu^\ast n^\ast$, i.e. $\mu^\ast\ge 4\mu r\Delta^2/|J_{t+1}|$. Need that $ \frac{d'^{3(t+1)}}{2} (n/r)^2 \ge \gamma^\ast n^{*2}$ and thus $\gamma^\ast \ge \frac{d'^{3(t+1)}}{2|J_{t+1}|^2}$} the probability that $\sigma$ is conflict-free is at least $\eul^{-\mu^{1/3} n}$.
\endclaimproof

It remains to show that with high enough probability, the updated candidacy graphs $A_{t+1}^j$ are super-regular. For this, we use the Four graphs lemma.

\begin{NoHyper}
\begin{claim}\label{claim:candidate update}
For all $j\in [r]\sm J_{t+1}^\ast$, with probability at least $1-(1-a)^{n/r}$, the candidacy graph $A_{t+1}^j$ is $(\eps_{t+1},d'')$-super-regular for some $d''\ge d'^{t+2}$.
\end{claim}
\end{NoHyper}

\claimproof
Crucially, $A_{t+1}^j$ only depends on (at most) one of the $\sigma_i$.
Consider $x\in X_j$ and $v\in V_j$. Observe first that by~\eqref{def candidate set} and~\eqref{def new phi}, we have $xv\in E(A_{t+1}^j (\phi_{t+1}))$ if and only if $xv\in E(A_{t}^j (\phi_{t}))$ and $\sigma(y)v\in E(G^\ast)$ for all $y\in N_H(x)\cap \bX_{t+1}$.
Recall that $j$ has at most one $R$-neighbour in $J_{t+1}$. If $j$ has no such neighbour, then $A_{t+1}^j=A_{t}^j$ and thus there is nothing to prove. Assume now that $i$ is the unique $R$-neighbour of $j$ in $J_{t+1}$. Thus, $\pi_{j\to i}(x)$ is the unique $H$-neighbour of $x$ in $\bX_{t+1}$. We conclude that $xv\in E(A_{t+1}^j)$ if and only if $xv\in E(A_{t}^j)$ and $\sigma_i(\pi_{j\to i}(x))v\in E(G^\ast)$. Thus, $A_{t+1}^j$ only depends on $\sigma_i$. Moreover, recalling the definition of $P^j_i$, we have that
\begin{align}
xv\in E(A_{t+1}^j) \quad \Leftrightarrow \quad \pi_{j\to i}(x)v \in E(P^j_i) \mbox{ and } \sigma_i(\pi_{j\to i}(x))v\in E(G^\ast). \label{four graphs identification}
\end{align}

Consider the triple $(\tilde{A}_t^i,P^j_i,G^\ast[V_i,V_j])$ on vertex set $X_i\cupdot V_i \cupdot V_j$. Recall that $\sigma_i\colon X_i\to V_i$ is a uniformly random perfect matching of $\tilde{A}_t^i$.
Define $P^j_i(\sigma_i)$ as the spanning subgraph of $P^j_i$ which contains all edges $x'v\in E(P^j_i)$, with $x'\in X_i,v\in V_j$, for which $\sigma_i(x')v\in E(G^\ast[V_i,V_j])$. (Thus, $P^j_i(\sigma_i)$ is the graph $A_\sigma$ defined in the Four graphs lemma.) By~\eqref{four graphs identification}, $P^j_i(\sigma_i)$ is isomorphic to $A_{t+1}^j$.

By~\eqref{sparsified regular}, $\tilde{A}_t^i$ is $(2\Delta \sqrt{\eps_t},d_i)$-super-regular. Moreover, by \ind{t}, $P^j_i$ (which is isomorphic to $A_t^j$) is $(\eps_t,d_j)$-super-regular. Also, recall that $G^\ast[V_i,V_j]$ is $(\eps_0,d')$-super-regular. Finally, condition~\eqref{four graphs condition} is satisfied by~\eqref{four graphs condition satisfied}. Therefore, we can apply the Four graphs lemma (Lemma~\ref{lem:four graphs}) as follows:

\medskip
{
\noindent
{
\begin{tabular}{c|c|c|c|c|c|c|c|c|c|c|c}
 $n/r$ & $a$ & $\eps_t^{1/3}$ & $\eps_{t+1}$ & $\tilde{A}_t^i$ & $P^j_i$ & $G^\ast[V_i,V_j]$ & $d_i$ & $d_j$ & $d'$ & $\sigma_i$ & $P^j_i(\sigma_i)$   \\ \hline
 $n$ & $a$ & $\eps$ & $\eps'$ & $G_{1,2}$ & $G_{1,3}$ & $G_{2,3}$ & $d_{1,2}$ & $d_{1,3}$ & $d_{2,3}$ & $\sigma$ & $A_\sigma$
\end{tabular}
}
}
\medskip

With probability at least $1-(1-a)^{n/r}$, the graph $P^j_i(\sigma_i)$ (and thus $A_{t+1}^j$) is $(\eps_{t+1},d_jd')$-super-regular. Note that $d_jd'\ge d'^{t+2}$, as required.
\endclaimproof

Using a union bound, it follows from Claims~\ref*{claim:conflict free} and~\ref*{claim:candidate update} that with probability at least $$\eul^{-\mu^{1/3} n}-r(1-a)^{n/r}>0,$$ $\sigma$ has the property that $\phi_{t+1}$ satisfies \ind{t+1}. This completes the proof.
\endproof

\section{Applications} \label{sec:apps}

\subsection{Rainbow spanning trees in Dirac graphs}\label{sec:appstree}

Although the following result is implied by the rainbow bandwidth theorem which we prove in the next subsection, we state and prove it here to illustrate an easy application of our blow-up lemma.

\begin{theorem}\label{thm:diractrees}
Suppose $1/n\ll\mu\ll \delta,1/\Delta$.
Suppose $G$ is a graph on $n$ vertices such that $\delta(G)\geq (1/2+\delta)n$.
Then, given any $\mu n$-bounded edge colouring of $G$,
the graph $G$ contains any tree $T$ on $n$ vertices with $\Delta(T)\leq \Delta$ as a rainbow subgraph.
\end{theorem}

Koml\'os, S\'ark\"ozy and Szemer\'edi \cite{KSS:95} proved this result in the non-rainbow setting even before the development of the blow-up lemma (some of the ideas which later led to the blow-up lemma were already present there). With a blow-up lemma at hand, the proof becomes much shorter. It essentially boils down to distributing the vertices of $T$ evenly among the clusters of the regularity partition, with the reduced graph being a Hamilton cycle.
For this, we need the following result, which essentially appears in \cite{JKKO:ta} and uses the fact that the symmetric Markov chain on an odd cycle mixes rapidly.

\begin{lemma}[\cite{JKKO:ta}]\label{lem:randomwalk}
Suppose $1/n\ll 1/r, 1/\Delta$ and $r$ is odd.
Suppose $T$ is a tree on $n$ vertices with $\Delta(T)\leq \Delta$.
Then there is a partition $(X_i)_{i\in[r]}$ of $V(T)$
such that
\begin{enumerate}[label={\rm (\roman*)}]
	\item the endvertices of any edge of $T$ lie in two consecutive sets $X_{i},X_{i+1}$ (where $X_{r+1}$ is identified with $X_{1}$),
	\item $|X_i|=n/r\pm n/\log^2 n$ for all $i\in [r]$, and
	\item for every $i\in[r]$,
	the vertex set $X_{i+1}$ contains at least $2^{-\Delta-3}n/r$ vertices $v\in V(T)$ such that $N_T(v)\subseteq X_{i}$.
\end{enumerate}
\end{lemma}

We remark here that for a cycle on $n$ vertices
it is easy to find a vertex decomposition as described in Lemma~\ref{lem:randomwalk}.
As the following proof only uses these properties of $T$ and the fact that $\Delta(T)\leq \Delta$,
(the proof of) Theorem~\ref{thm:diractrees} also implies the existence of a rainbow Hamilton cycle under the same assumptions on $G$.%
\COMMENT{Hence it implies the result of Cano, Perarnau and Serra. Guillem said they leave it with the Eurocomb abstract. There's not gonna be a paper.}

\lateproof{Theorem~\ref{thm:diractrees}}
Let $G$ and $T$ be as in the statement; in particular, we assume that some $\mu n$-bounded edge colouring of $G$ is given.
Choose new constants $\eps,d$ such that $\mu\ll \eps \ll d \ll \delta,1/\Delta$.
We apply the regularity lemma to $G$
and obtain a vertex partition $(V_i)_{i\in[r]_0}$
such that
\begin{itemize}
	\item $|V_1|=\ldots = |V_r|$ for some odd $r$ (we may assume $\mu \ll 1/r \ll \eps$),\COMMENT{So $|V_i|\le n/r$}
	\item $|V_0| \leq \eps n$, and
	\item for every $i\in [r]$ and for all but at most $\eps r$ integers $j\in [r]\sm\{i\}$,
	the graph $G[V_i,V_j]$ is $\eps$-regular.
\end{itemize}
Let $R$ be the graph with vertex set $[r]$
where two vertices $i,j\in[r]$ are joined by an edge if $G[V_i,V_j]$ is lower $(\eps,d)$-regular.
As $\delta(G)\geq (1/2+\delta)n$, it is not hard to see
that $\delta(R)\geq (1/2+\delta/2)r$.
Hence $R$ contains a Hamilton cycle $C$. Without loss of generality, we may assume that $1,2,\ldots,r$ appear in this order on~$C$. Below we always identify vertex sets $S_{r+1}$ with $S_1$.

By Fact~\ref{fact:regularity}, one can remove at most $2\eps n/r$ vertices from every $V_i$  (call the new vertex sets $\widetilde{V}_i'$)
such that $G[\widetilde{V}_i',\widetilde{V}_{i+1}']$ is lower $(4\eps,d)$-super-regular for all $i\in [r]$.
For purposes that will become clear below, we remove even more vertices.
So let us remove from every $V_i$ exactly $\eps^{1/2}n/r$ vertices (in total) and add them to $V_0$ (call the new vertex sets $V_i'$)
such that  $G[V_i',V_{i+1}']$ is lower $(2\eps^{1/2},d)$-super-regular for all $i\in [r]$ (where $V_{r+1}':=V_1'$).
We claim that
\begin{enumerate}[label=(\roman*)]
	\item $\eps^{1/2}n \leq |V_0'|\leq (\eps^{1/2}+\eps)n$,
	\item for every vertex $v\in V_0'$, there are at least $(1/2+\delta/2)r$ integers $i\in [r]$
	such that $v$ has at least $d n/2r$ neighbours in $V_i'$,
	and
	\item for every $i\in [r]$, there are at least $(1/2+\delta/3)|V_0'|$ vertices in $V_0'$ that have at least $dn/2r$ neighbours in $V_i'$.
\end{enumerate}
Indeed, (i) is trivial.
To see (ii), consider $v\in V_0'$.
Clearly, $v$ has at most $|V_0'|\leq 2\eps^{1/2} n$ neighbours in $V_0'$ and at most $dn$ neighbours in clusters $V_i'$ in which $v$ has less than $d n/2r$ neighbours.
As $d\ll \delta$ and $\delta(G)\geq (1/2+\delta)n$, statement (ii) follows.
To see (iii),
we fix some $i\in [r]$ and let $j\in N_R(i)$.
As $G[V_i,V_j]$ is lower $(\eps,d)$-regular,
among the $\eps^{1/2}n/r$ vertices we moved from $V_j$ to $V_0$,
at least $(\eps^{1/2}-2\eps)n/r$ have at least $(d-\eps)|V_i|\ge dn/2r$ neighbours in $V_i$ by Fact~\ref{fact:regularity}.
Now (iii) follows as $d_R(i)\geq (1/2+\delta/2)r$ and by (i).

Next we apply Lemma~\ref{lem:randomwalk} to $T$
and obtain a vertex partition $(X_i)_{i\in[r]}$ of $V(T)$ such that the endvertices of every edge lie in consecutive parts (modulo $r$)
and $|X_i|=n/r\pm n/\log^2 n$ for all $i\in[r]$.
Moreover, $X_{i+1}$ contains at least $2^{-\Delta-3}n/r$ vertices such that all their neighbours lie in $X_{i}$;
among those vertices we select a $2$-independent set $\widehat{X}_{i+1}$ of size exactly $|X_{i+1}|-|V_{i+1}'|$. For all $i\in[r]$, since $|X_i|= n/r \pm n/\log^2 n$ and $|V_i'|= (1- \eps^{1/2} \pm \eps)n/r$, we have that
\begin{align}\label{eq:size hat X}
	|\widehat{X}_{i}|=(1\pm \delta/10)\eps^{1/2}n/r.
\end{align}
Let $X_i':=X_i\setminus \widehat{X}_{i}$.
Hence $|X_i'|=|V_i'|$.

Let $X_0':=\widehat{X}_1\cup \dots \cup \widehat{X}_r$.
It is easy to see that $X_0'$ is $2$-independent and $|V_0'|=|X_0'|$. We wish to apply the rainbow blow-up lemma to the blow-up instance $(T,G,C,(X_i')_{i\in[r]_0},(V_i')_{i\in[r]_0})$. Before, we need to find a feasible partial embedding $\phi_0\colon X_0'\to V_0'$ and define appropriate candidacy graphs for the remaining vertices.
Next, we want to match the vertices in $V_0'$ and $X_0'$
such that whenever $v\in V_0'$ is matched to $x\in \widehat{X}_{i+1}$,
the vertex $v$ has at least $dn/2r$ neighbours in $V_i'$.
To this end, let $H$ be the bipartite graph with bipartition $(V_0',X_0')$
and join $v\in V_0',x\in \widehat{X}_{i+1}$ whenever $d_G(v,V_i')\ge dn/2r$.
By (iii), we have $d_H(x)\geq (1/2+\delta/3)|V_0'|$ for all $x\in X_0'$.
Moreover, for all $v\in V_0'$, we have
$$d_H(v) \overset{{\rm (ii)},\eqref{eq:size hat X}}{\geq} (1/2+\delta/2)r\cdot (1-\delta/10)\eps^{1/2}n/r \overset{\rm(i)}{\geq} (1/2+\delta/3)|V_0'|.$$
Clearly, by Hall's theorem, $H$ contains a perfect matching $\phi_0\colon X_0'\to V_0'$.

Now, consider $i\in [r]$. We need to define a suitable candidacy graph $A^i$. Note that at most $2\Delta \eps^{1/2}|X_i'|$ vertices of $X_i'$ have a neighbour in $X_0'$.\COMMENT{In particular, condition (iii) from the rainbow blow-up lemma is satisfied.} For these vertices we have to restrict their candidate set. Now suppose that $x\in X_i'$ has a (unique) neighbour $x'$ in $X_0'$. Observe that we must have $x'\in  \widehat{X}_{i+1}$. We define $N_{A^i}(x):=N_G(\phi_0(x'),V_i')$. By definition of $H$, we have that $|N_{A^i}(x)|\ge dn/2r$. For the vertices $x\in X_i'\setminus N_T(X_0')$, we can define $N_{A^i}(x):=V_i'$. Thus, $A^i$ is lower $(\eps^{1/3},d/2)$-super-regular.
Moreover, we have met condition~\ref{exc condition:neighbourhoods}, and since $X_0'$ is $2$-independent, $\phi_0$ is $2\Delta\mu n$-feasible.
Thus, since the blow-up instance $(T,G,C,(X_i')_{i\in[r]_0},(V_i')_{i\in[r]_0})$ with candidacy graphs $(A^i)_{i\in[r]}$ is lower $(\eps^{1/3},d/2)$-super-regular, an application of the rainbow blow-up lemma (Lemma~\ref{lem:blow up}) yields a rainbow embedding of $T$ into $G$.
\endproof

\subsection{Rainbow bandwidth theorem} \label{subsec:bandwidth}

%One of the fundamental results of extremal combinatorics is the Erd\H{o}s--Stone theorem, stating that for a fixed graph $H$, every large graph $G$ on $n$ vertices with average degree at least $(1-1/(\chi(H)-1)+o(1))n$ contains $H$ as a subgraph.
%The bandwidth theorem can be viewed as an analogue of the Erd\H{o}s--Stone theorem when $H$ is a spanning subgraph of~$G$. Clearly, in this setting, one has to replace the average degree condition by a minimum degree condition in order to obtain sensible results. 
%Let $\ell:=\chi(H)$ and assume that $H$ has bounded degree.
%A long line of research confirmed for various cases that $\delta(G)\geq (1-1/\ell+o(1))n$ suffices to find $H$ as a spanning subgraph in $G$, for instance when $H$ is a spanning tree, a (power of a) Hamilton cycle, or a clique factor.
%A conjecture of Bollob\'as and Koml\'os, which became known as the bandwidth conjecture, attempted to
%generalize all the mentioned results by claiming that $\delta(G)\geq (1-1/\ell+o(1))n$ suffices whenever $H$ has not too strong expansion properties (in this case measured by the parameter bandwidth).
%B\"ottcher, Schacht and Taraz proved this conjecture roughly ten years ago~\cite{BST:09}. We refer to their paper for more information on the history of the conjecture.
%By using our rainbow blow-up lemma, we extend the bandwidth theorem to the rainbow setting (see Theorem~\ref{thm:bandwsimple}).

In this subsection, we apply the rainbow blow-up lemma to prove a rainbow bandwidth theorem (Theorem~\ref{thm:bandwsimple}). Due to the proof in~\cite{BST:09} being carried out in a very modular way, we can directly make use of important parts of the original proof.
In fact, B\"ottcher, Schacht and Taraz~\cite{BST:09} proved a slightly more general result (and so do we). As mentioned earlier, the minimum degree threshold in the bandwidth theorem is in general not optimal. For example, if $H$ is a Hamilton cycle and $n$ is odd, then $\chi(H)=3$, yet $\delta(G)\ge n/2$ already suffices to find~$H$. Roughly speaking, this is because $H$ is essentially $2$-colourable, except that we need an additional colour for one vertex. This motivates the following definition.
Suppose the vertex set of $H$ is $[n]$.
For $\ell,x,y\in \mathbb{N}$,
we say an $(\ell + 1)$-colouring $\sigma\colon V (H) \to [\ell]_0$ of $H$ is \defn{$(x, y)$-zero free}
if for each $t \in [n]$,
there exists a $t'$ with $t \leq  t' \leq t+x$ such that $\sigma(u) \neq 0$ for all $u \in [t', t'+y]$.

The following theorem can provide the optimal threshold for such graphs $H$, for example if $H$ is the $(\ell-1)$-st power of a Hamilton cycle. We remark however that it does not always yield the optimal threshold, e.g.~for $C_5$-factors the optimal threshold is $3/5$ (cf.~\cite{KO:09}).
\begin{theorem}\label{thm:bandwidth}
Suppose $1/n\ll \mu , \beta \ll \delta,1/\Delta,1/\ell$.
Let $H$ be a graph on vertex set $[n]$ with $\Delta(H)\leq \Delta$. Assume that $H$ has an $(\ell+1)$-colouring that is $(8\ell \beta n,4\ell \beta n)$-zero free
and uses colour $0$ at most $\beta n$ times. Also assume that $\max \{|i-j|\colon ij\in E(G)\}\le \beta n$.
Suppose $G$ is a graph on $n$ vertices with $\delta(G)\geq (1-1/\ell+\delta)n$. Then, given any $\mu n$-bounded edge colouring of $G$, there is a rainbow copy of $H$ in~$G$.
\end{theorem}

Obviously Theorem~\ref{thm:bandwidth} implies Theorem~\ref{thm:bandwsimple}, and the remaining part of this section is devoted to the proof of Theorem~\ref{thm:bandwidth}. Its proof in~\cite{BST:09} without the rainbow condition is based on three results, the so-called `Lemma for $G$', `Lemma for $H$' and `partial embedding lemma' (see Lemmas~6, 8 and~9 in~\cite{BST:09}). The first two are strong structural decomposition results for $G$ and $H$, respectively,
which prepare the application of the blow-up lemma and require the main part of the work done in~\cite{BST:09}.
Fortunately, for the proof of Theorem~\ref{thm:bandwidth} we can use the Lemmas for $G$ and $H$ without any alterations.

The partial embedding lemma deals with a very small set of exceptional vertices that cannot be handled by the blow-up lemma directly.
In~\cite{BST:09}, this is by far the easiest part and follows by standard techniques.
However, it does not easily translate to the rainbow setting and hence we give a proof later in this section using the methods developed in this paper.
We now state the partial embedding lemma for the rainbow setting, which replaces Lemma 9 from \cite{BST:09}.
Roughly speaking, it yields a rainbow embedding $\phi$ of a small exceptional set $X$ such that every outside neighbour $y\in Y:=N_H(X)\sm X$ of vertices in $X$ still has a large candidate set (and will later be embedded by the blow-up lemma).

\begin{lemma}\label{lem:partemb}
Suppose $1/n\ll\mu, \eps \ll d' \ll  d,1/\Delta$ and $\mu \ll 1/r$.
Suppose $R$ is a graph on~$[r]$.
Suppose $G$ is a graph on $n$ vertices with vertex partition $(V_i)_{i\in [r]}$
such that $G[V_i,V_j]$ is lower $(\eps,d)$-regular whenever $ij\in E(R)$
and $|V_i|=(1\pm \eps)n/r$.
Suppose $H$ is a graph on at most $\eps n/r$ vertices with vertex partition $(X_1,\ldots,X_r,Y_1,\ldots,Y_r)$
such that $\Delta(H)\leq \Delta$ and $Z_i:=X_i\cup Y_i$ is an independent set for all $i\in [r]$.
Moreover,
whenever $z_iz_{j}\in E(H)$ with $z_i\in Z_i$, $z_{j}\in Z_{j}$, then $ij\in E(R)$.
Let $X:= X_1\cup \ldots \cup X_r$ and $Y:=Y_1\cup \ldots \cup Y_r$.
Given any $\mu n$-bounded edge colouring $c\colon E(G)\to C$ of $G$,
there exists a set of colours $C'\In C$, a rainbow embedding $\phi\colon H[X]\to G_{C\sm C'}$  
and candidate sets $S_y\subseteq V_i$ for all $y\in Y_i$ and $i\in [r]$
such that the following hold:
\begin{enumerate}[label={\rm (\roman*)}]
\item $\phi(X_i)\subseteq V_i$ for all $i\in [r]$;
\item for all $y\in Y$, we have $|S_y|\geq d' n/r$ and for all $x\in N_H(y)\cap X$ we have $S_y\In N_{G_{C'}}(\phi(x))$;
\item for all $y\in Y$, all distinct $x,x'\in N_H(y)\cap X$ and all $v\in S_y$, we have $c(\phi(x)v)\neq c(\phi(x')v)$.
\end{enumerate}
\end{lemma}

With the Lemmas for $G$ and $H$ at hand,
the proof of Theorem~\ref{thm:bandwidth} without the rainbow condition
works verbatim also for Theorem~\ref{thm:bandwidth}, by replacing the partial embedding lemma with the rainbow partial embedding lemma and the blow-up lemma with the rainbow blow-up lemma.
We provide a proof for completeness, but focus on the differences to the original proof. 
For integers $r,\ell$,
let $K^r_\ell$ denote the vertex-disjoint union of $r$ cliques of order~$\ell$. This will be a spanning subgraph of the reduced graph $R$ indicating super-regular pairs, and the blow-up lemma will be applied with reduced graph $K^r_\ell$. In~\cite{BST:09}, $K^r_\ell$ is found within $R$ inside a so-called `backbone', that is, the $\ell$-cliques in $K^r_\ell$ form a sequence where two consecutive cliques are joined almost completely in~$R$. The purpose of this backbone is to allow changing the cluster sizes $V_i$ in $G$ slightly to match the cluster sizes of the $H$-partition given by the Lemma for~$H$.

\lateproof{Theorem~\ref{thm:bandwidth}}
Choose new constants $d,d',\eps,\xi,r',r''$ such that we have the following hierarchy of parameters
\begin{align*}
	1/n\ll \mu , \beta \ll \xi \ll 1/r'' \ll 1/r' \ll \eps \ll d'\ll d \ll \delta, 1/\Delta, 1/\ell.
\end{align*}
Suppose now that $G$ and $H$ are as in the statement. In exactly the same way as in~\cite{BST:09}, we first use the structural results for $G$ (Lemma~6 in~\cite{BST:09}) and $H$ (Lemma~8 in~\cite{BST:09}), respectively,
to obtain a partition $(V_i)_{i\in [\ell r]}$ of $V(G)$, a partition $(X_i\cup X_i')_{i\in [\ell r]}$ of $V(H)$ and a graph $R$ on $[\ell r]$ for some $r'\leq r\leq r''$ such that the following hold:
\begin{enumerate}[label=(\roman*)]
	\item $|\bigcup_{i\in [\ell r]}X_i'|\leq \xi n$;
	\item $|X_i\cup X_i'|= |V_i|= (1\pm \xi)\frac{n}{\ell r}$ for all $i\in [\ell r]$;
	\item $R$ contains a $K_\ell$-factor $K^r_\ell$ such that $G[V_i,V_j]$ is lower $(\eps,d)$-regular for all $ij\in E(R)$ and even lower $(\eps,d)$-super-regular if $ij\in E(K^r_\ell)$;
	\item for all $xy\in E(H)$ with $x\in X_i\cup X_i'$ and $y\in X_j\cup X_j'$, 
	we have $ij\in E(R)$, 
	and moreover we have $ij\in K_\ell^r$ if $x\in X_i$ and $y\in X_j$.
\end{enumerate}

For the first time, we need to make an easy amendment for the rainbow setting. Let $c\colon E(G) \to C$ be any $\mu n$-bounded edge colouring. Note that we may assume that $|C|\le 3\mu^{-1}n$, say. (Otherwise, there would be two colour classes with combined size at most $\mu n$ and we could merge them.)
We split the colour set $C$ into two sets $C_1,C_2$ such that $G_{C_i}$ has still $R$ as a reduced graph
and all lower (super-)regular pairs halve their density.
This can be achieved with Lemma~\ref{lem:separate colours}.

Let $X_0:=\bigcup_{i\in [\ell r]}X_i'$.
Next we use the (rainbow) partial embedding lemma (Lemma~\ref{lem:partemb} instead of Lemma~9 in~\cite{BST:09}), with $X_0,N_H(X_0)\sm X_0,H[X_0\cup N_H(X_0)]$ playing the roles of $X,Y,H$, to obtain a set $C'\In C_1$, a rainbow embedding $\phi_0\colon H[X_0] \to G_{C_1\sm C'}$ and a candidate set $S_y$ for each $y\in N_H(X_0)\sm X_0$ such that the following hold:
\begin{enumerate}[label={\rm (\alph*)}]
\item $\phi_0(X_i')\subseteq V_i$ for all $i\in [\ell r]$; \label{partial embedding location}
\item for all $y\in N_H(X_0)\sm X_0$ with $y\in X_i$ for some $i\in[\ell r]$, we have $S_y\In V_i$ and $|S_y|\geq d' n/\ell r$, and for all $x\in N_H(y)\cap X_0$, we have $S_y\In N_{G_{C'}}(\phi_0(x))$; \label{partial embedding candidates}
\item for all $y\in N_H(X_0)\sm X_0$, all distinct $x,x'\in N_H(y)\cap X_0$ and all $v\in S_y$, we have $c(v\phi(x))\neq c(v\phi(x'))$. \label{partial embedding rainbow}
\end{enumerate}

Let $V_0':=\phi_0(X_0)$ and $V_i':=V_i\setminus V_0$ for all $i\in [\ell r]$. By~\ref{partial embedding location}, we have $|V_i'|=|X_i|$ for all $i\in [\ell r]$. Let $H':=H-H[X_0]$.

Now we want to use the rainbow blow-up lemma to embed $H'$ into $G_{C_2\cup C'}$.\COMMENT{Adding edges with colours in $C'$ doesn't destroy the lower regularity of $G_{C_2}$}
As a blow-up instance we use $(H',G_{C_2\cup C'},K^r_\ell,(X_i)_{i\in [\ell r]_0},(V_i')_{i\in [\ell r]_0})$ with the obvious candidacy graphs, where the candidate set of $y\in N_H(X_0)\sm X_0$ is~$S_y$ and there are no restrictions for other vertices.
From~\ref{partial embedding candidates},~\ref{partial embedding rainbow} and the fact that $\Delta(H)\leq \Delta$, we conclude that $\phi_0$ is $2\Delta \mu n$-feasible. Thus, Lemma~\ref{lem:blow up} yields a rainbow embedding $\phi$ of $H'$ into $G_{C_2\cup C'}$ which extends~$\phi_0$. Since $\phi_0$ and $\phi$ use distinct colours, this is a rainbow embedding of~$H$, completing the proof.
\endproof

%\subsubsection{Proof of Lemma~\ref{lem:partemb}}

It remains to prove Lemma~\ref{lem:partemb}, for which we need the following straightforward consequence of Lemma~\ref{lem:switchings} and Proposition~\ref{prop:many switchings}.

\begin{cor}\label{cor:manyswitch}
Suppose $1/n\ll \mu, \eps \ll d$ and $\mu\ll 1/r$.
Let $V_1,\ldots,V_r$ be disjoint vertex sets such that $|V|=n$, where $V:=\bigcup_{i\in[r]}V_i$, and $n/(2r)\leq |V_i|\leq 2n/r$ for all $i\in [r]$.
Suppose $X_1,\ldots,X_r$ are disjoint vertex sets such that $|X_i|\leq \eps n/r$ for all $i\in [r]$.
Let $X:=\bigcup_{i\in[r]}X_i$.
Suppose $G$ is a bipartite graph with vertex partition $(X,V)$
and suppose every edge in $G$ joins $X_i$ and $V_i$ for some $i\in [r]$.
Suppose that $d_G(x)\geq dn/r$ for all $x\in X$.
Then, given any $\mu n$-bounded system of conflicts for $G$,
there is a conflict-free matching covering $X$.
\end{cor}
\proof
Define a supergraph $G'$ of $G$ on $2n$ vertices as follows.
First we add to every set $X_i$ exactly $|V_i|-|X_i|$ new vertices and call the new set $X_i'$. Let $X':=\bigcup_{i\in[r]}X_i'$.
Next, we join every vertex $x\in X_i'\sm X_i$ to all vertices in $V_i$.
It is easy to see that $G'[X_i',V_i]$ is lower $(4\eps,d/2)$-super-regular for all $i\in [r]$.
Every conflict system of $G$ easily transfers to a conflict system of $G'$.

Observe that by Proposition~\ref{prop:many switchings},
for every $i\in [r]$, the graph $G'[X_i',V_i]$ has a perfect matching, and thus $G'$ has one.
In addition, for every perfect matching $M$ of $G'$ and every edge $e\in  M$,
there are at least $d^3/9\cdot n^2/(4r^2)$ edges that are $(e,M)$-switchable.
Hence, Lemma~\ref{lem:switchings} implies the existence of a conflict-free perfect matching $M'$ in $G'$.
Let $M\subseteq M'$ be the matching obtained by deleting all edges incident to $X'\setminus X$.
Consequently, $M$ is a conflict-free matching of $G$ covering $X$.
\endproof

The following proof is reminiscent of our main proof in Section~\ref{subsec:main proof} in that we proceed in rounds and reserve exclusive colours for each round. It is much simpler though since we only need to embed a small fraction of all vertices.

\lateproof{Lemma~\ref{lem:partemb}}
As $\Delta(H)\leq \Delta$, we can colour $H^2$ with $T:=\Delta^2+1$ colours.
Hence, there exists a partition of $X$ into sets $X^1,\ldots,X^{T}$ which are $2$-independent in~$H$.\COMMENT{independent partition of $X_1,\dots,X_r$}
We define $X^{T+1}:=Y$. We will proceed in $T$ rounds and embed in round $t\in [T]$ all vertices of $X^t$, whilst keeping track of candidate sets for the remaining vertices.

Beforehand, we apply Lemma~\ref{lem:separate colours}\COMMENT{actually just McDiarmid here} to partition $C$ into sets $\{C_{t_1t_2}\colon t_1t_2\in \binom{[T+1]}{2}\}$ such that $G_{C_{t_1t_2}}[V_i,V_j]$ is lower $(2\eps,\hat{d})$-regular for all $t_1t_2\in \binom{[T+1]}{2}$ and all $ij\in E(R)$, where $\hat{d}:=d/\binom{T+1}{2}$.
Define $G^{t_1t_2}:=G_{C_{t_1t_2}}$.
%
%to partition $G$ into $\binom{\Delta^2+2}{2}$ edge-disjoint graphs $\{G^{ij}\colon ij\in \binom{[\Delta^2+2]}{2}\}$
%such that $G^{ij}$ is lower $(2\eps,\hat{d})$-regular on $R$ for every $i\in [\Delta^2+1]$, where we define $\hat{d}:=d/\binom{\Delta^2+2}{2}$,
%and two distinct graphs do not share a common colour.
%
%We proceed in $\Delta^2+1$ rounds and embed in round $t$ the vertices in $X^t$.

In the beginning we set $S_z(0):=V_i$ for all $i\in [r]$ and $z\in Z_i$.
We claim that after round~$t$, we have
\begin{enumerate}[label=(\alph*)]
	\item a rainbow embedding $\phi^t$ of $H^t:=H[X^1\cup \dots \cup X^t]$ such that for every edge $xy\in E(H^t)$ with $x\in X^{t_1}$ and $y\in X^{t_2}$, we have $c(\phi^t(xy))\in C_{t_1t_2}$;
	\item $|S_z(t)|\geq (\hat{d}/2)^t\cdot n/2r$ and
	$S_z(t)\subseteq V_{i}\setminus \phi^{t}(V(H^{t}))$ for all $i\in[r]$ and $z\in Z_i \sm V(H^{t})$;
	\item for all $t'\in [t]$, $t''\in [T+1]\setminus [t]$, and all $x'\in X^{t'},x''\in N_H(x')\cap X^{t''}$,
	we have $S_{x''}(t)\subseteq N_{G^{t't''}}(\phi^t(x'))$.
\end{enumerate}

Indeed, our claim is true for $t=0$.
Now let $t>0$ and suppose it is true for $t-1$. We first restrict the candidate sets for the vertices in $x\in X^t$ slightly.
Consider any $x\in X^t$.
By Fact~\ref{fact:regularity}, for all but at most $\eps^{1/2}n/r$ vertices $v\in S_x(t-1)$, we have 
\begin{align}
d_{G^{tt'}}(v,S_{x'}(t-1))\ge 3\hat{d}/4 \cdot |S_{x'}(t-1)| \overset{\rm{(b)}}{\ge} 3\hat{d}/4  \cdot (\hat{d}/2)^{t-1} \cdot  n/2r. \label{candidate shrink}
\end{align}
for all $t'>t$ and $x'\in N_H(x)\cap X^{t'}$.
Let $S_x(t)$ be obtained from $S_x(t-1)$ by removing all those vertices $v\in S_x(t-1)$ for which~\eqref{candidate shrink} does not hold for some $t',x'$. This will ensure that (b) and (c) hold again for the next step.

We now embed $X^t$ by picking for every $x\in X^t$ a suitable image from $S_x(t)$.
Let $V:=V(G)$.
We define a bipartite graph $J$ with bipartition $(X^t,V\setminus \phi^{t-1}(V(H^{t-1})))$
and we join every $x\in X^t$ to all vertices in $S_x(t)$. To ensure that the new embedding is again rainbow, we define a system of conflicts for $E(J)$. For an edge $xv\in E(J)$ with $x\in X^t$ and $v\in S_x(t)$, let $F_{xv}$ be the set of edges $\phi^{t-1}(x')v$ for $x'\in N_H(x)\cap V(H^{t-1})$. By (c), we have $F_{xv}\In E(G)$. Moreover, the colours of the edges in $F_{xv}$ are distinct as the sets $X^{t'}$ are $2$-independent.
We let two distinct edges $xv,x'v'\in E(J)$ conflict if $c(F_{xv})\cap c(F_{x'v'})\neq \emptyset$.
It is easy to see that this gives a $\mu n$-bounded conflict system.\COMMENT{}
We apply Corollary~\ref{cor:manyswitch} and obtain a conflict-free matching $\sigma \colon X^t\to V\setminus \phi^{t-1}(V(H^{t-1}))$.
Now, let $\phi^t$ be the extension of $\phi^{t-1}$ defined by $\phi^t(x):=\sigma(x)$ for all $x\in X^t$.
Observe that $\phi^t$ is a rainbow embedding of $H^t$ with the properties required for~(a).

Next we update the candidate sets $S_x$ for all $x\in X^{t'}$ with $t'>t$. 
If $x$ has no neighbour in $X^t$ we only remove $\phi^t(X^t)$ from $S_x(t-1)$ to obtain $S_x(t)$.
If $x'\in X^t$ is the neighbour of $x$ (there is at most one for every $x$), we define $$S_x(t):=(S_x(t-1)\cap N_{G^{tt'}}(\sigma(x')))\sm \phi^t(X^t).$$
 This automatically ensures that (c) holds again for~$t$. Moreover, by~\eqref{candidate shrink} and since $|X^t|\leq \eps n/r$, we deduce that $|S_x(t)|\ge 3\hat{d}/4 \cdot |S_{x'}(t-1)|- \eps n/r \ge (\hat{d}/2)^t\cdot n/2r$, hence (b) also holds for $t$.

Assume now that the claim holds for $T$. Let $C':=\bigcup_{t\in [T]}C_{t,T+1}$.
By (a), $\phi:=\phi^T$ is a rainbow embedding of $H[X]$ into $G_{C\sm C'}$.
For all $y\in Y$, let $S_y:=S_y(T)$. It remains to check that (i)--(iii) hold. Clearly, (i) holds as we have $S_z\In V_i$ for all $z\in Z_i$ and all $i\in[r]$ throughout the embedding.
Since $d'\leq (\hat{d}/2)^{T}/2$,
we have $|S_y|\geq d' n/r$ for all $y\in Y$ by~(b). Moreover, from (c) we deduce that for all $y\in Y$, $t\in [T]$ and $x\in N_H(y)\cap X^t$, we have $S_y\In N_{G^{t,T+1}}(\phi(x))$. Since $G^{t,T+1}\In G_{C'}$, this establishes (ii). Finally, for all $y\in Y$, all distinct $x,x'\in N_H(y)\cap X$ and all $v\in S_y$, we have $c(\phi(x)v)\neq c(\phi(x')v)$, i.e.~(iii) holds. Indeed, by the above, we have $c(\phi(x)v)\in C_{t,T+1}$ and $c(\phi(x')v)\in C_{t',T+1}$, where $x\in X^t,x'\in X^{t'}$. Since we partitioned $X$ into $2$-independent sets, $t$ and $t'$ are distinct and thus $C_{t,T+1}\cap C_{t',T+1}=\emptyset$.
This completes the proof.
\endproof

\section{Concluding remarks} \label{sec:conclusion}

We have proved a rainbow blow-up lemma for $\mu n$-bounded edge colourings and have used it to prove a rainbow bandwidth theorem in this setting and to show that every bounded degree spanning graph exists as a rainbow subgraph in quasi-random graphs of arbitrarily small fixed density. We conclude this paper with the following remarks:

\begin{itemize}
\item In fact, our blow-up lemma applies to slightly more general systems of conflicts (see beginning of Section~\ref{sec:main proof}), allowing for instance to obtain embeddings which are simultaneously rainbow with respect to several given edge colourings. A natural question is whether a blow-up lemma still holds for arbitrary $\mu n$-bounded conflict systems (as defined in Section~\ref{sec:conflictM}). The bottleneck in the current proof is that it relies on the colour splitting technique, which seems to be limited to `highly transitive' conflict systems.

\item Note that in order to guarantee a rainbow copy of a graph $H$ with maximum degree $\Delta$ in a graph $G$ of density $d$, 
the given edge-colouring of $G$ needs to be $dn/\Delta $-bounded, 
as otherwise there may be less than $\Delta n/2$ different colours available.
In particular, our theorems are optimal up to the value of the constant~$\mu$.
As noted before, the constant $\mu$ in our theorems is very small. 
In particular, in an embedding obtained with the rainbow blow-up lemma, only a small fraction of the colours available in $G$ is used for the embedding. 
On the contrary, affirmative answers to many open rainbow conjectures would imply that (almost) all colours need to be used. 
As mentioned in Section~\ref{sec:intro}, 
Kim, K\"uhn, Kupavskii and Osthus
used the rainbow blow-up lemma on a small random subset (of vertices and colours) to complete a partial embedding (even a partial approximate decomposition), thus effectively using almost all colours.
We expect further applications in this direction.

\item There has been some exciting progress towards rainbow decompositions of properly edge-coloured complete graphs, for instance in~\cite{BLM:17,PS:17} it is shown that there are $\Omega(n)$ edge-disjoint rainbow spanning trees in every properly edge-coloured~$K_n$. For $\mu n$-bounded edge colourings, it might be possible to achieve (approximate) rainbow decompositions for any prescribed collection of bounded degree graphs.
For instance, it was proved in~\cite{AST:95}
that for every $\mu n$-bounded edge colouring of $K_{n,n}$,
there is a decomposition into rainbow perfect matchings provided $\mu$ is small enough and $n$ is a power of~$2$.
We conjecture that for any $\mu n$-bounded edge colouring of~$K_n$, there exists a decomposition into rainbow Hamilton cycles (provided $n$ is odd). Similarly, 
we conjecture that for any collection of $n$-vertex graphs $H_1,\dots,H_t$ with bounded degree and $\sum_{i=1}^t e(H_i) \le (1-\alpha)e(K_n)$ and any $\mu n$-bounded edge colouring of~$K_n$, where $\mu \ll \alpha$, the graphs $H_1,\dots,H_t$ pack edge-disjointly into $K_n$ such that each subgraph is rainbow. In the uncoloured setting, this was proved in~\cite{KKOT:ta}.

\end{itemize}

\section*{Acknowledgement}

We are grateful to David Harris as well as Matthew Coulson and Guillem Perarnau for helpful discussions on their papers.

\bibliographystyle{amsplain_v2.0customized}
\bibliography{References}

\vfill

\small
\vskip2mm plus 1fill
\noindent
Version \today{}
\bigbreak

\noindent
Stefan Glock
{\tt <s.glock@bham.ac.uk>}\\
Felix Joos
{\tt <f.joos@bham.ac.uk>}\\
School of Mathematics\\ 
University of Birmingham\\
United Kingdom

\end{document}